\documentclass[a4paper,10pt]{article}
\usepackage{hyperref}

\usepackage[utf8]{inputenc}

\usepackage[english ]{babel}
 
\usepackage{xcolor}
\usepackage{bbding}

\usepackage[utf8]{inputenc}

\usepackage{tikz-cd}

\usepackage{amsfonts}
\usepackage{amssymb}
\usepackage{graphicx}
\usepackage{mathtools}
\usepackage{skak} 
\usepackage{foekfont}

\usepackage[T1]{fontenc}

\usepackage{amsmath,mathdots}

\usepackage[utf8]{inputenc}
\usepackage{devanagari}

\usepackage{mathtools}
\usepackage{mathdots}
\usepackage{tikz}

\definecolor{amber(sae/ece)}{rgb}{1.0, 0.49, 0.0}
\newfont{\rsfsten}{rsfs10 scaled 1200}

\usepackage{MnSymbol}
\usepackage{tikz}
\usepackage{amscd}
\usepackage{MnSymbol}
\usepackage{wasysym}

\usepackage{mathrsfs}

\newcommand*{\rom}[1]{\expandafter\@slowromancap\romannumeral #1@}

\usepackage{xcolor}

\usepackage{wrapfig}
\newcommand{\tightunderset}[2]{%
  \mathop{#2}\limits_{\vbox to .3ex{\kern-0.95ex\hbox{$#1$}\vss}}}

\usepackage[utf8]{inputenc}

\usepackage{CJKutf8}

\usepackage{pmboxdraw}

\usepackage{stmaryrd}
\usepackage{url}
\usepackage{harmony}
\usepackage{ifsym}

\usepackage{eufrak}

\usepackage{rotating}

\usepackage{dingbat}

\title {Lectures on  Immersions with Controlled Curvatures}

\author{Misha Gromov\footnote{Part of the material presented in these lectures were  preparered during author's stay at the  Isaac Newton Institute on the programm Operators, Graphs, Groups in July-August 2025.}}

\begin{document}

\maketitle

{\bf Historical Preamble:} {\color {blue}Heinz Hopf}: Selected Topics in Geometry, New York University 1946, Notes by {\color {blue} Peter Lax.}\footnote{\url{https://link.springer.com/book/10.1007/3-540-39482-6} Among many other  things, there is   a proof of Legendre-Cauchy-A. Schur "Arms-Bow-Lemma" on pp 31-32 in these lecture (attributed by Hopf to E. Schmidt),   which  has been reproduced in all further publication concerning this theorem. e.g. in \url {https://www.scribd.com/document/759520702/Chern-Curves-and-surfaces-in-Euclidean-spaces}}
\vspace {1mm}

{\color {blue}Heinz Hopf} (19 November 1894 – 3 June 1971)\vspace {1mm}

In 1925, he proved that any simply connected complete Riemannian 3-manifold of constant sectional curvature is globally isometric to Euclidean, spherical, or hyperbolic space.

In 1931, Hopf discovered the Hopf invariant of maps $S^3\to S^2$  (“element of the architecture of our world” in the words of Penrose)
 and proved that the Hopf fibration has invariant 1.
 This: 
 
 (1) disproved the then standing  intuitive conjecture that the continuous maps between spheres $S^N\to S^n$,  $N>n$, are contractible;  

(2) Opened the door to the world of vector bundles and the topology of spinors, where the curvature of the Hopf bundle is 1/2 
 curvature of the 2-sphere.

(Hopf bundle and Dirac Monopole \url{https://personal.math.ubc.ca/~mihmar/HopfDirac.pdf},\url{https://www.sciencedirect.com/science/article/abs/pii/S0393044002001213} \url{https://ncatlab.org/nlab/show/Hopf\%20fibration})
\vspace {1mm}

{\color {blue} Peter David Lax} (1 May 1926 – 16 May 2025) \vspace {1mm}

After the war ended, Lax remained with the Army at Los Alamos for another year and eventually  returned to NYU for the 1946–1947 academic year.

\begin{abstract}

When does a smooth $n$-manifold $X$ admit an immersion to the unit $N$-ball $Y=B^N=B^N(1)\subset \mathbb R^N$,
 such that   the  (normal) curvature of this immersion is {\it bounded  by a 
given constant} $c$? 

If $N$ is significantly greater than $n$,   then 
we know  what the critical constant $c$ is: 

(a) {\it all} $n$-manifolds $X$ {\it immerse} to $B^{20n^2}(1)$ with curvatures $<c=\sqrt 3$; 

(b) if $c=\sqrt 3-\varepsilon$ then, for all $n>{1\over 2\varepsilon}$ and all $N$,  {\it there are $n$-manifolds}, which     admit {\it no immersions }
 to $B^N(1)$ with curvatures $<c$.
 
 We explain this and  recapitulate what little is known and much of what is 
  unknown for $N$ comparable  with $n$.\footnote {Most our conjectures, as it is common in mathematics, are not based on a positive  evidence, but rather on  the absence  of  evidence to the contrary.}

\end {abstract}
  \tableofcontents
  

  \section{Definitions, Problems and First  Examples}


  It is amazing how little is known on how much a {\it bound on the normal curvature} constrains 
   the {\it geometry and topology}  of a submanifold in the Euclidean space.

   For instance, let 
$$X\overset {f}\hookrightarrow B^N(1)\subset \mathbb R^N$$
be a smooth {\it immersion}\footnote {
        "Immersion" signifies a 
$C^1$-map $f:X\to Y$ between smooth manifolds, such that
 the differentials $df:T(X)\to T(Y)$ nowhere vanishes, 
 $df(\tau)=0\implies \tau=0, \tau\in T(X).$
 
 Immersions are locally one-to-one maps, but  globally  they may have self intersections;  immersions 
 {\it without self intersections} are called 
{\it embeddings}, where, for  non-compact $X$, one usually  require   
  the induced topology in $X$ to be equal the original one.} 
   of a {\it closed $n$-dimensonal manifold} to the unit ball in the Euclidean N-space.\vspace {1mm}
 
 {\color {blue!24!black}Does a bound on the normal curvature,
   {\color {blue}$$curv^\perp(X)=curv^\perp(f)= curv^\perp\big(f( X))=curv^\perp\big( X\overset {f}\hookrightarrow B^N(1)\big )\leq C,$$}
impose a nontrivial constrain on the topology of $X$?}
   
   \vspace {1mm}


{\it \bf Notation:  $\bf \overrightarrow{\bf curv}_\tau$, $\bf curv_x^\perp$ and $\bf curv^\perp(X)$.} Let $f:X\to \mathbb R^N$ be a smooth immersion, let 
$\tau =\tau_x\in T_x(X)$ be a tangent vector and let $\gamma_\tau\subset X$ be a geodesic in $X$ issuing  from $x$ with the speed $\tau$. 
Then  the {\it normal curvature vector} 
$$\overrightarrow{ curv}_\tau(X)\in \mathbb T_{f(x)}(X)\in\mathbb R^N)=\mathbb R^N$$
is equal the acceleration  (the second derivative) at $f(x)$  of a point moving along  the curve 
$f(\gamma)$ in    $\mathbb R^N$.

Granted this,  define  
$$curv_x^\perp(X)=\sup_{\|\tau_x\|=1} { curv}_{\tau_x}(X)\mbox { and } curv^\perp(X)=\sup_{x\in X}{ curv}_{x}(X).$$

If  $dim(X)=1$, say $X=[0,1]$ and the curve $X\overset {f}\hookrightarrow \mathbb R^N$ is parametrized by 
arc length, that is 
$$\left\|\frac {d f(x)}{dx}\right\|=1,$$ then     
  this is  the usual curvature of a curve, 
$$\overrightarrow {curv}(X, x)=\frac {d^2 f(x)}{dx^2}  \mbox {and $curv^\perp(X)=\sup_{x\in X}\left \|\frac {d^2 f(x)}{dx^2}\right\|$}.$$ 

Thus, \vspace {1mm}

{\it the normal curvature  $curv^\perp (X\overset {f}\hookrightarrow  \mathbb R^N)$ is equal to  the supremum of the normal curvatures  of the $f$-images in $\mathbb R^N$ of the  geodesics from $X$.} 

\vspace {1mm}

{\it \textbf{Locality of the Curvature  and Curvature of Submanifolds.}}  Since curvature of an immersion  at a point $x\in X$  is a local invariant and since immersions locally are embeddings, the definition and many properties of curvatures of immersions  formally follow from those for submanifolds $X\subset \mathbb X^N$. In this in mind, we may often (but not always)  speak   of curvatures  of  "immersed submanifolds", and accordingly  simplify our notation.

\vspace {1mm}

{\it \textbf {Curvatures of Spheres.}}
Spheres $S^n(R)$  of  radius $R$  of all dimensions $n$  in the $N$-space  $\mathbb R^N$, $N>n$,   satisfy 
$$ curv^\perp(S^n(R))=\|\overrightarrow {curv}_\tau^\perp(S^n(R))\|=1/R\mbox {  for all unit tangent vectors } \tau\in T(S^n(R)).$$

 \SunshineOpenCircled \hspace {1mm} {\sl The unit  $n$-spheres $S^n(R=1)\subset  B^N(1)$, are {\it the only}  {\it closed connected immersed} $n$-sub-manifolds
 with curvatures $\leq 1$ in the Euclidean space $\mathbb R^N$, which are contained  in  the unit   $N$-ball $B^N(1)$, except for $n=1$, where   multiple covering of the unit circle are also such manifolds. }

   This follows by the  {\it maximum principle}   applied to the    distance function from  $X$ to the boundary $\partial B^N(1))$ or equivalently to  the squared distance to the centre of the ball $B^N$
 denoted $r^2(x)$. 
 
 In fact, since  $curv^\perp(X)\leq 1$,  the second derivatives of   $r^2$ along   geodesics parametrized by the arc length  satisfy:
 $$\mbox { $\|r''r\|\leq 1$ and $(r^2)''=2(r''r+ \|r'\|^2)\geq 0$, since  $\|r'\|^2=1$.}$$ 
  
  This says that $r^2$
    is { \it  a convex}, hence constant=1 function on $X$. 
    Thus, $X$ is contained in the unit sphere $S^{N-1}(1)=\partial B(1)$, where it has zero normal curvature by Pythagorean formula (1.1.C); hence, totally geodesic,  (compare with 3.B).

 {\it \textbf {Compact, Closed, Complete}.} Curvature has a limited effect on immersion of {\it open manifolds}, i.e. those which contain{\it  no compact   connected components without boundaries}, called  
 {\it closed manifolds}.
 
 For instance, according to  the generalized {\it Smale-Hirsch h-principle}\footnote{See  [C-E-M]  and references therein.},  an arbitrary  immersion $f$ of open manifold $X$  to an open subset  $U\subset \mathbb R^N$ admits a homotopy  (even a {regular homotopy}\footnote
 {A regular homotopy is a path in the space of $C^1$ immersions with the usual $C^1$-topology, that is  $f_t: X\to Y$, $t\in [0,1]$, where the differential $df_t$ of $f_t$ in  $x$-variables, $x\in X$, is continuous in $t$.} to an immersion $f_\varepsilon$,  such that 
 $$curv^\perp (X\overset {f_\varepsilon} \hookrightarrow U)\leq\varepsilon\mbox {  for a given $\varepsilon>0.$}$$
 
In what follows, we focus on  immersions of closed manifold $X$, where much of what we do equally applies to {\it complete immersed} manifolds $X\hookrightarrow \mathbb R^N$, i.e. where  the induced Riemannian metrics in $X$, sometimes called {\it inner metrics}, are {\it geodesically  complete}: geodesics starting at all  point $ x\in X$  extend infinitely in all directions $\tau_x\in T_x(X)$.

{\it Exercise.} Generalize   \SunshineOpenCircled \hspace {1mm} to complete immersed $X\hookrightarrow B^N(1)$.

\vspace {1mm}

{\it \textbf {Extremal Immersions.}}  We are much interested in $curv^\perp$-{\it extremal immersions} between Riemannian  manifolds, $f:X\hookrightarrow Y$,  especially for $Y=\mathbb R^N$, which minimize some geometric size  invariant of the image $f(X)\subset Y$, such as $diam_Y(f(X))$, among all immersions with $curv^\perp\leq c$\footnote {Our definitions of $curv^\perp$ naturally  generalize to all Riemannian manifold $Y$ receiving immersions from $X$.} or among  all such  immersion {\it regularly homotopic} to a given one.
Beside the diameter, it may be, some kind of {\it width}, the {\it radius  of the minimal ball}  which contained $f(X)$, etc.

If we don't 	specify any  invariant,  we call an immersion $f_0:X\hookrightarrow  Y$ {\it simple extremal} if it admits {\it  no regular homotopy}  
 $f_t :X\hookrightarrow Y$, such that $curv^\perp (f_1)<curv^\perp(f_0)$, where the local version of this says that    
 all  regular homotopies $f_t$, satisfy  $curv^\perp (f_t)\geq curv^\perp(f_0)$ for $t>0$.
 
If $Y=\mathbb R^N$, then this  may be applied to the convex hull  $Y_0=conv(f(X))\supset f(X)$ and then  
an immersion $f _0:X\hookrightarrow \mathbb R^N$ is called  {\it $conv$-$curv^\perp$-extremal}  if one can't decrease the  
normal curvature of $f_0$  by a regular homotopy of immersions $f_t:X \hookrightarrow conv(f_0(X))$.

{\bf 1.A.} {\it \textbf {Basic Spherical  Example.}}
By \SunshineOpenCircled,   spheres $S^n(1/c)\subset \mathbb R^N$ are extremal  with respect to all above criteria.\vspace {1mm}

{\bf 1.B.} {\it \textbf {Piecewise $C^2$ Circular Example.}} Some naturally arising  submanifolds with
 bounded normal curvatures, e.g. many extremal ones    
are $C^1$-smooth and only {\it piecewise $C^2$}.\footnote {This is a well know phenomenon in the {\it optimal control theory}, where one is predominantly concerned with $n=1$, [Feld 1965].}

For instance,  immersed closed 
curves, which  go around several circles in the plane, e.g.  {\huge$\Circle$}\hspace {-3.5mm}$\Circle${\large $^{\Circle}$},  by switching their tracks from  one circle to another at the points where  two circle  "touch" one another.

Such curves are $C^1$-smooth but they are not $C^2$: their curvatures  jump  when   a curve switches the track    from  one circle to another  at the contact points between the  circles. The $curv^\perp$ curvature of such a curve is equal the  reciprocal of the  radius of the smallest circle involved.

{\bf 1.C.}{\it \textbf{ {\small \bf $\Circle$\hspace {-0.35mm}$\Circle$}-Subexample.}}  Let $f:S^1\hookrightarrow B^2(1)\subset 
\mathbb R^2$ be  a $ C^1$ immersion with curvature
$$curv^\perp  (S^1\overset {f}\hookrightarrow \mathbb R^2)\leq 2.$$
  
{\sl If the corresponding oriented Gauss map to the unit circle
$$\overrightarrow G_f=\frac{df}{\|f\|}:S^1\to S^1\subset \mathbb R^2$$
has degree zero (hence contractible), then 
the image of $f$ is equal to the union of two circles of radii 1/2, which meet  at the center of the disc $B^2(1)$, where they are tangent one to another. } Thus the figure $\infty$ immersion is "radially extremal": it {\it minimizes  the  
radius  of the  2-ball}. around it.  
(We shall explain why this is so  in section 10).
\vspace {1mm}

{\it \textbf {Bi-invariants  $\bf curv^\perp_{min}(X,Y)$ and  $\bf Imm_{\perp\leq c}(X,Y)$}}.
Let  $X$ be a smooth  closed manifold and 
and $Y$   a Riemannian manifold and let  $ curv^\perp_{min}(X,Y)$ be the infimum of normal curvatures of smooth immersions $X\hookrightarrow Y$.

Now, if we choose and fix  a particular  $Y$, e.g. the unit ball in $\mathbb R^N$, the number $ curv^\perp_{min}(X,Y)$ becomes     a {\it topological invariant} 
of $X$, the value of
 which is unknown for most $n$-manifolds and $N>n$ .
 
 Dually, given a topological  $n$-manifold, e.g.
  (homeomorphic to) the product 
of spheres, the minimal $curv_{\min}(X,Y)$ of immersions $X\hookrightarrow Y$ appears as a {\it metric invariant} 
of $Y$, which is unknown in most cases,  for instance, for the $N$-balls and cubes $Y\subset  \mathbb R^N$.

The number $\bf curv^\perp_{min}(X,Y)$  carries only a small part of the information 
about immersions $f:X\hookrightarrow Y$  with  curvatures $curv^\perp(f)\leq c$.

A more comprehensive information is contained in the {\it homotopy types} of the spaces of immersions with 
$curv^\perp(f)\leq c$, denoted $Imm_{\perp\leq c}(X,Y)$ and the {\it homotopy classes} of the inclusion maps
$$Imm_{\perp\leq c_1}(X,Y)\subset Imm_{\perp\leq c_2}(X,Y), \mbox { } c_1\leq c_2,$$ 
where much of  this information is encoded by the diagram of the natural (co)homology 
homomorphisms between these spaces.


\subsection{Alternative Definitions  of Normal Curvature}


   
  {\bf 1.1.A.} The  full {\it second order infinitesimal information} of a smooth submanifold  $X$ in a Riemannian manifold $Y$, e.g. in the Euclidean $N$-space, at a point $x\in X$ is algebraically represented by the  {\it second fundamental form} 
  that is a {\it symmetric bilinear}  form on $X$ with values in the normal vector space $T^\perp(X)\subset T(Y)$,
  denoted 
    $${\rm II}(X,x) ={\rm II}(X,x,\tau_1,\tau_2) = {\rm II}_x(\tau_1,\tau_2), $$
    where $\tau_1,\tau_2\in T_x(X)$ are tangent vectors to $X$ and where  the value $\rm II(\tau_1,\tau_2)$ is a vector in $T_x(Y)$ normal to the tangent (sub)space  $T_x(X)\subset T_x(Y)$.     
 This form in the case $Y=\mathbb R^N$ is defined  as the {\it second differential} of a vector  function, say 
   $\Phi:  T_x(X)\to T^\perp_x(X)$, such that 
  the graph of $\Phi$ in a neighbourhood of $x\in\mathbb R^N\supset X$ is equal to $X\subset \mathbb R^N=T_x(\mathbb R^N)=T_x(X)\oplus T_x^\perp(X)$,
  $${\rm II}_x(\tau_1,\tau_2)= \partial_{\tau_1} \partial_{\tau_2}\Phi(x)$$
   
   In the general case, this definition applies by equating  $T_x(Y))$ with a small neighbourhood in $Y$ via the exponential map $\exp_x: T_x(Y)\to Y$.

   {\it  Exercises.}  {\bf  1.1.B.} Show that 
      ${\rm II}(\tau,\tau)$ is
    equal 
     to the second   (covariant) derivative in $Y$ of the geodesic in $X$ issuing from $x$ with the velocity  $\tau$,
    and that 
   $$\|{\rm II}_x(\tau_1,\tau_2)\|\leq curv^\perp_x(X)\leqno {[\tau_1\tau_2]_\leq}$$
 for all $x\in X$ and all unit tangent vectors $\tau_1,\tau_2\in T_X(X)$.
 
 
  {\bf  1.1.C. Pythagorean Curvature Composition} Let  $X\hookrightarrow  Y \hookrightarrow Z$ be isometric embeddings (or immersions) between Riemannian manifolds, i.e 
the Riemannian metrics in $Y$ and in $X$ are induced from a Riemannian metric in $Z$.
Show that 
{\color {blue}$$curv^\perp_\tau (X\hookrightarrow Z) =\sqrt {(curv^\perp_\tau (X\hookrightarrow Y))^2 +(curv^\perp_\tau (Y\hookrightarrow Z))^2 }$$}
for all tangent vectors $\tau \in T(X)\hookrightarrow T(Y) \hookrightarrow T(Z). $

For instance, if   $X\hookrightarrow  Y=S^{N-1}(1)  \hookrightarrow Z=\mathbb R^N$ 
  then 
{\color {blue}$$curv^\perp_\tau (X\hookrightarrow \mathbb R^N) =\sqrt {(curv^\perp_\tau (X\hookrightarrow S^{N-1}))^2 +1}.$$}

{ \textbf {1.1.D. Geodesic free definition of $curv^\perp$.}  Show that 
the normal curvature  $curv^\perp_\tau (X\hookrightarrow Y)=\|\overrightarrow{curv}_\tau \| $, $\|\tau\|=1$,     is equal to  the infimum  of the $Y$-curvatures  $curv_Y^\perp$-curvatures of the curves 
in $X$ tangent to $\tau$.}  \vspace {1mm}

This description of $curv^\perp$,  which doesn't refer to geodesics, has an advantage of being  applicable to mechanical systems with {\it non-holonomic} constrains that are submanifolds in the tangent bundle  of $Y$, say
$\mathcal  X \subset T(Y)$ 
rather  $X\subset $.

  {\bf 1.1.E. Metric Definition of II.} Let $Y=(Y,g)$ be a Riemannian manifold, e.g. $Y=(\mathbb R^N, g=\sum_{j=1}^N dy_j^2$), let $X\subset Y$ be a smooth submanifold, let $\nu\in T^\perp_x(X)$ be a normal vector to $X$ at $x$ and $\tilde \nu$ be a smooth vector field on $Y$, which extend $\nu_x$. 
  
  Let $g_{|X}$ be the restriction of the Riemannian quadratic form $g$ to $X$ and let $\tilde g'_{|X}$ be the restriction of (Lie) {\it derivative of $g$ by the field $\tilde \nu$ to $X$.}

Show that the value $\tilde g'_{|X}(\tau_1,\tau_2)$ for $\tau_1,\tau_2\in T_x(X)$ depends only on $\nu$ but not on the extension $\tilde \nu$  of $\nu$.

Moreoever, show that 
$$\tilde g'_{|X}(\tau_1,\tau_2)=\langle\nu,{\rm II}_x(\tau_1,\tau_2)\rangle_g, $$ 
and that the second fundamental form II is {\it uniquely determined}  by this identity.

(The  definition of the second fundamental form II as the derivative   $\tilde g'_{|X}$ of the induced Riemannian form uses no  covariant  derivatives or   geodesics  either in $X$ or in   $Y$.)

{\bf 1.1.F. Normal Curvature Defined via the Gauss Map.} Let $\mathcal H=Gr_n(N)$ be the space of $n$-dimensional linear  subspaces 
$H\subset \mathbb R^N$ and natuarally identify  the tangent space  $T_H(\mathcal H$  with the space of linear maps from $H$ to the normal space $H^\perp\subset \mathbb R^N$,
$$T_H(\mathcal H)=hom(H, H^\perp).$$

Let $f:X \hookrightarrow \mathbb R^N$ be a smooth immersed submanifold   and 
$ \overleftrightarrow G: X\to Gr_n(N)$, $n=dim(X)$, be the  (non-oriented) Gauss map where $\overleftrightarrow G(x)$ is the linear subspace parallel to  tangent  subspace  of $X$ in  $R^N $ (regarded as an affine   subspace) at $x$.

Let $D_x\overleftrightarrow G: T_x(X)\to T_x(X)^\perp$ be the differential of the map  $\overleftrightarrow G$ at $x\in X$ regarded as a linear operator  $T_x(X)\to T_x^\perp(X)^\perp$.
 
Show that 

{\it the normal curvature of $X$ at $x$ is equal to the norm of the operator $D_x\overleftrightarrow G$,
$$curv_x^\perp(X)=\sup_{\tau\in T_x(X), \|\tau\}=1}  \|D_x\overleftrightarrow G(\tau)\|\leqno {[D\overleftrightarrow G]^\perp}$$}
and derive from this the following corollary.

{\bf 1.1.G.  Angular  Arc Inequality.} If the (inner) distance   between two points $x_1,x_1\in X$ satisfies 
  $$dist_X(x_1,\underline x)\leq {\alpha} (curv^\perp(X))^{-1},\mbox {   } \alpha\leq \pi/2,  $$ 
then the angles between   vectors $\tau\in T_{x_1}(X)$ and their images $\bar \tau$ under the   normal projection  
$T_{x_1}(X)\to T_{x_1}(X)$ satisfy
$$\angle (\tau, \bar \tau)\leq \alpha,$$    
where  the equality holds if and only if  there exists  a 

{\it planar $\alpha$-arc of radius $1\over curv^\perp(X)$, which is  contained in $X$, which   join $x_1$ with $\underline x$
and such that  $\tau$ is tangent to this arc at its  $x_1$-end.}

Conversely,  the inequality $\angle (\tau, \bar \tau)\leq \epsilon/c+o(\epsilon)$, $c\geq 0$,  for all pairs of $\epsilon$-infinitesimally closed points 
 implies that  $curv^\perp(X)\leq c.$

\hspace {10mm}{\it  no non-zero tangent vector $\tau_1\in T_{x_1}(X)$ is normal to $T_{\underline x}(X)$}. 

Moreover the same non-normality conclusion holds if $$dist_X(x_1,\underline x)\leq {\pi\over 2} (curv^\perp(X))^{-1}, $$ 
unless   there exists  a 

{\it planar semicircle of radius $1\over curv^\perp(X)$ contained in $X$ and  joining $x_1$ with $\underline x$.}

\hspace {1mm}

{\bf  1.1.H.  Polygonal Curves}. Given a spacial  polygonal curve $P$ with vertices $p_i$ let $c_i=c_i(P)=c(P,p_i)$ denote the "external" angles  of $P$, 
$$c_i=c(P,p_i)=\pi-\angle(P,p_i)$$
where $0\leq \angle(P,p_i)\leq \pi$ is the angle between two segments of the curve   adjacent to $p_i$ that are $[p_{i-1},p_i]$ and $[p_i,  p_{i+1}]$.\footnote{According to our present  convention  even non-convex planar polygons have all there angles measured between zero and $\pi$.}

{\bf\octagon }\hspace {1mm} Let $P\subset \mathbb R^N $  be a {\it closed. connected}   spacial polygonal curve with $k$ vertices $p_i$ (where $p_1=p_k$ by the cyclic convention) and let us decompose $P$ to triangles $\triangle _j$, e.g. by drawing $ k-3$ segments $[p_{1} ,p_i]\subset \mathbb R^N$. for all $i\neq 1, 2,  k-1$.\footnote{One can decompose $P$   to  $\approx \log_2 k$ triangles.} 

Observe that the angles of triangles $\triangle_{j_i}$  adjacent to  $p_i$ satisfy the following "trianagle kind  inequaity"
$$\sum_{j_i} \angle (\triangle_{j_i},p_i)\geq \angle (P,p_i) \mbox{  for all } i,$$
where the equality implies that

$ \bullet $ the triangles $\triangle_{j_i}$ lie in the same 2-plane (which depends on $i$) 
 
$ \bullet $ the triangles $\triangle_{j_i}$  do not overlap

$ \bullet $  the union of these trianles is convex.

It follows that the sum of the "external angles" of  $P$ satisfies 
$$\sum_{i=1}^kc_i\geq 2\pi,$$
 where the equality  $\sum_{i=1}^kc_i\geq 2\pi$ holds if  and only if 
 $P$ is a {\it planar convex} curve.

{\it Exercise.} Let $P\hookrightarrow \mathbb R^2$ be a planar  connected, immersed {\it locally convex} polygonal curve. 

Show that $\sum_i c_i(P)=2\pi d$, where $d$ is a positive integer and that two such curves $P_0$ and $P_1$
can be joint by a homotopy $P_t$, $t\in [0,1]$ of immersed locally convex curves if and only if
$\sum c_i(P_0)=\sum c_i(P_1)$.

{\bf \bf  1.1.I.  Polygonal  Approximation  Exercises.}
Let a spacial curve $X$ be represented by a continuous map $x:[0,l]\to \mathbb R^N$
and let $s_i =\varepsilon i\in [0,l]$, for $i=1,2,...,k$ and $\varepsilon=l/k$ and let $P_\varepsilon\hookrightarrow \mathbb R^N$ be the polygonal curve with  vertices $p_i=x(s_i)\in \mathbb R^N$  and segments 
$[p_i,p_{i+1}]
\subset \mathbb R^N$

(a) Show that if $X$ is a  smooth immersed curve, then the sums  of the "external angles" of $P_\varepsilon$ approximate   the total curvature of $X$
$$\sum_i c_i(P_\varepsilon)\to \int_ 0^l curv^\perp(x(s))ds\mbox { for } \varepsilon\to 0,$$ 
$s\in[0,l]$  is    arc length  parameter on $X$.   

(b) Use this  approximation for the definition of the curvature of $X$ as the density of the weak limit of the  measures  
$\sum_i c_i\delta(p_i)$, where $\delta(p_i)$ are Dirac's $\delta$-measures at the points $p_i.$

(c) Similary define curvature measures of  {\it piecewise smooth} immersed curves and also localy convex planar curves.

(d) Prove  {\color {blue}\it  Fenchel total curvature $\geq 2\pi$-Inequality} for  closed smooth immersed curves 
$$\int_ {S^1} curv^\perp(x(s))ds\geq 2\pi$$ 
 where $X$ is paramerized by the unit circle $S^1$
and then generalise this to all closed curves.  
  \footnote {See [Chern] and also sections 7-9 for related matters.}
  
 (e) Construct  smooth embedded curves $X_{d, \varepsilon}\subset \mathbb R^3$ for all $d=1,2,...$ and $\varepsilon>0$ with total curvatures $2\pi+\varepsilon $, which have linking numbers $d$ with a straight  line in  $\mathbb R^3$.

(f) Let $r$ be the normal  projection from the 3-space to the  line  $\mathbb R=\mathbb R_r$ in space   defined  by a unit vector $r\in S^2\subset \mathbb R^3$.    
Let $N_r(X)$ be the number of  critical points of the composed function $S^1\mathbb R_r$  for  $s\mapsto r\circ x(s)$.

Show that 
$$\int _{S^2}N_r(X)dr= 4\int_ {S^1} curv^\perp(x(s))ds$$

{\it Hint}. Apply Crofton's formula to the spherical curve $s\mapsto x'(s)\in S^2$.

(e) Prove {\color {blue} \it Fáry–Milnor theorem.}  If the total  curvature a closed   embedded  $X\subset \mathbb R^3$ is strictly less than $2\pi$ then $X$ is unknotted.


\section   { Products of Spheres, Clifford's sub-Tori with Small Curvatures and Petrunin Inequality}


The product $X$ of spheres
 $S^{n_i}(R_i)\subset \mathbb R^{N_i=n_i+1}$, $i=1,...m$,
  $$X=S^{n_1}(R_1)\times S^{n_2}(R_2)\times  ...\times S^{n_m}(R_m)\subset  \mathbb R^{N=(n_1+n_2+...+n_m)+m}, $$
has  the
 curvature equal to the maximum of $ 1/R_i$, $i=1....,m,$ and if 
  $$R_1^2+R_2^2+...+R_m^2\leq 1,$$
   then  $X$ is contained in the unit ball in $\mathbb R^N.$ (If $R_1^2+R_2^2+...+R_m^2=1,$ then 
   $X$ is contained in the unit sphere $S^{N-1}(1)=\partial B^N(1)\subset \mathbb R^N.$)
 
 For example, the product of $m$-copies of  $S^{n}$ admits an embedding to the unit ball in $\mathbb R^{mn+m}$, where 
   {\color {blue}$$curv^\perp\big ( (S^n)^m\subset B^{mn+m}(1)\big)=\sqrt  m$$}

 The main instance of this   is the {\it \color {blue} Clifford $n$-torus}, that 
   the product of $n$ circles imbedded to the unit $2n$ ball, such that     
   {\color {blue}$$curv^\perp \big ( \mathbb T^n\subset  \partial B^{2n}(1)\big)=\sqrt n.$$}

  \vspace {1mm}
 
 It is {\color{magenta} conceivable} that the  above (Clifford's)  products of spheres $S^{n_1}(R_1)\times S^{n_2}(R_2)\times  ...\times S^{n_m}(R_m)\subset  \mathbb R^{N}$  are {\it $conv$-$curv^\perp$-extremal,}
 where this {\color{magenta} seems  realistic}   
 for $m<min_i n_i$,
  but we have {\color{magenta} no idea}, for instance,   if     there are  immersions of $n$-tori  to 
 $B^{2n}(1)$ with   $curv^\perp<\sqrt n$. 
 
 Yet, if  $N>>n$, then the $n$-torus can be immersed  to the unit ball $B^N(1)$ with {\it unexpectdly small} curvature.


  {\bf  2.A.  $\sqrt 3$-Clifford Sub-Torus  Theorem.} (Section ?) {\bf[a]} If $N$ is much greater than $n$,  then  the   Clifford torus 
  $$\mathbb T^N\subset S^{2N-1}\subset  B^{2N}(1),$$
  contains an $n$-subtorus $\mathbb T_o ^n\subset \mathbb T^N,$
   such that the   normal curvature of this $n$-torus  the ambient Euclidean space 
  $\mathbb R^{2N}\supset B^{2N}\supset \mathbb T_0^n$ satisfies
  {\color {blue} $$curv^\perp\big (\mathbb T_o^n\subset  B^{2N}(1)\big )\leq\sqrt \frac {3n}{n+2}.\leqno { \left [\frac {3n}{n+2}\right]_{\mathbb T^n}}$$}

  One has  a poor bound on the best (i.e. the smallest) $N$ for this purpose, (something like $10^{10^n}$, see section 13) but

 {\bf [b]}  {\sl if  $N\geq 8n^2+8$, then, there exists a a locally isometric 
  (with respect to the  Euclidean metrics in $\mathbb R^n$ and  $\mathbb T^N$)
    map, that is 
a group homomorphism 
  $$g: \mathbb R^n\hookrightarrow \mathbb T^N\subset B^{2N}(1),$$
  such that 
  {\color {blue} $$curv^\perp\big (\mathbb R^n\hookrightarrow  B^{2N}(1)\big )\leq\sqrt \frac {3n}{n+2}.\leqno { \left [\frac {3n}{n+2}\right]_{\mathbb R^n}}.$$}}
  
 {\bf[c]}  It follows that 
  {\sl for all $\varepsilon>0$,    there exists  
  a sub-torus 	 
  $$\mathbb T_\varepsilon^n\subset \mathbb T^N\subset B^{2N}(1),$$
  such that 
{\color {blue} $$curv^\perp\big (\mathbb T_\varepsilon^n\subset  B^{2N}(1)\big )\leq\sqrt \frac {3n}{n+2}+\varepsilon.\leqno { \left [\frac {3n}{n+2}+\varepsilon\right]_{\mathbb T^n}}$$}}

 \vspace {1mm}

{\bf  {\bf  2.B.   $\sqrt 3$-Immersion Corollary.} {\sl  Let  $f:X\hookrightarrow \mathbb R^m$ be an immersion  then, for all $\epsilon>0$,  there exist   an immersion (actually an embedding) $f_\epsilon$  to the unit  ball 
$B^{16m^2+16m}$ with curvature 
  {\color {blue} $$curv^\perp\big (X\overset{ f_\varepsilon}\subset  B^{{16m^2+16m}}(1)\big )\leq\sqrt \frac {3m}{m+2}+\varepsilon.$$}}}
  {\it Proof.}  Let   $ \lambda$ be a large constant, $\lambda>>1/\epsilon$, {\it scale the manifold} 
  $X\overset {f}\hookrightarrow \mathbb R^m$ by $\lambda$ and 
  compose the scaled map $\lambda\cdot f:X\hookrightarrow \mathbb R^m$
  with the map $g: \mathbb R^m \hookrightarrow\mathbb T^N\subset B^{2N}(1)$
  from  the above {\bf [b]}.
  
  Then, if one one wishes, one slightly perturbs the resulting immersion 
  $X:\to \hookrightarrow\mathbb T^N $. and makes it an embedding.

  {\it On sharpness of $\left [\frac {3n}{n+2}\right]$. }  It is not hard to show that  the Euclidean  curvatures of all  Clifford subtori $\mathbb T^n\subset \mathbb T^N \subset \mathbb R^{2N}$ (these $\mathbb T^n$  are very special  submanifolds in $B^{2N}(1)\supset \mathbb T^N)$)  satisfy
  $curv_{\mathbb R^{2n}}^\perp(\mathbb T^n)\geq  \sqrt \frac {3n}{n+2}$,  but the following is not so obvious.

{\bf  2.C.  Petrunin's $\sqrt 3$-Inequality.} (Section 14.2) 
  {\it All immersions  $\mathbb T^n\hookrightarrow  B^N(1)$  satisfy  
  {\color {blue}  $$curv^\perp\big( \mathbb T^n\hookrightarrow  B^N(1)\big)\geq \sqrt \frac {3n}{n+2}\mbox {  for all $n\geq 1$ and all $N$.}$$} }

It is  {\color {magenta} unclear} what is, in general,  the geometry  of
immersions $\mathbb T^n\hookrightarrow  B^N(1)$ with 
 $curv ^\perp\approx\sqrt 3n $\footnote {Anton Petrunin told me that there exist extremal tori 
in $   B^N(1).$
which are not contained in $S^{N-1}$.}   depending on the ambient dimension $N$.
 {\color{magenta} Conceivably}   the $n$-tori admit no immersions $\mathbb T^n\hookrightarrow  B^{N}(1)$ with $curv^\perp\leq\sqrt{3}$ for $N<<n^2,$
 but we have  {\color{magenta} no means} to rule out such immersions, say for  
 for $N\leq {3}n$ and $n\geq 4$.

 


\section  {  Focal Radius and $+\rho$-Encircling}.


 Let $Y$ be  a complete Riemannian manifold, let   $X\hookrightarrow Y$ be  a smooth embedded or  immersed  submanifold,
  let $x_0\in X$, 
 let $\nu_0\in T_{x_0}^\perp(X)$  be a unit normal vector at the point $x_0$  and $\gamma_\nu\hookrightarrow  Y$   be a geodesic ray issuing from $x_0$
 in the $\nu_0$-direction. 

Define {\it $\nu_0$-focal radius} $rad^\perp_{\nu_0}(X)$ as the supremum of $r\geq 0$,   such that 
the the segment  $[x_0,y]\subset \gamma_0$ {\it locally  minimises} the length of  curves 
in $Y$ between $y$ and $X$, that is all curves,  which are sufficiently  close to the segment $[x_0,y]$ in $C^0$-topology 
 and which join $y$ and $X$,  have $length>r$.

Then let 
$$rad^\perp_{x_0}(X)=\inf_{\nu_0\in T_{x_0}(X)} rad^\perp_{\nu_0}(X)\mbox { and } rad^\perp(X)=
\inf_{x_0\in X} rad^\perp_{x_0}(X).$$ 

Equivalently, 
 {\sl the focal radius of $X\hookrightarrow Y$ is equal to the supremum of $r$, such that the normal exponential map $\exp T^\perp(X)\to Y$ is an {\it immersion} on the $r$-ball subbundle $B_X^\perp(r)\subset T^\perp(X).$ }

Observe the following.

{\bf 3.A. Eclidean Reciprocity.} The focal radius of a submanifold in a  Euclidean space is  equal to reciprocal of its normal curvature:
{\color {blue}$$ rad_x^\perp(X\hookrightarrow\mathbb R^N)=  {1\over curv_x^\perp(X\hookrightarrow\mathbb R^N)}.$$}

{\bf 3.B.   Focal Radius in $S^n$.}  Focal radii of  {\it submanifolds in the $R$-spheres} and  in the Euclidean spaces  satisfy the following relation: 
  $$rad^\perp_{\mathbb R^N}(X)=2R\sin {1\over 2}rad^\perp_{S^{N-1}}(R)(X),$$
which agrees with the Pythagorean formula for the  curvature of immersions $X\hookrightarrow S^N(R)$
    from 1.1.C:
  $$
  \Big (curv_{S^{n-1}}\big(X\hookrightarrow S^N(R)\big)\Big)^2 =\Big (curv^\perp_{\mathbb R^N}(X\hookrightarrow \mathbb R^N)\Big)^2 -1/R^2.$$

   For instance, 
   
  $\bullet$ the spherical  focal  radii of an equatorial subsphere (with zero spheriacal curvature) in the unit sphere $ S^{N-1}(1)$ is   equal to $\pi/2$, while their Euclidean focal radius  is equal to one;

    $\bullet$  the spherical  focal  radii of a  subsphere
    with spherical   radius  $\pi/4$ is also  $\pi/4$, and  the
        spherical curvature is  equal to one, while 
        the  Euclidean curvature  is   $\sqrt 2$ with agreement with he identity $\frac{\sin\pi}{4}=1/\sqrt 2$.

{\it  Exercises.} 
(i) Let $x_0\in X$ be a {\it local maximum} point in  $X\subset Y$ for the distance function $x\mapsto dist_Y(x,y_0)$  
for some $y_0\in Y$. Show that 
$$dist(x_0, y_0)\geq rad^\perp_{x_0}(X).$$

(ii). Let $rad^\perp_{x_0}(X)\geq r$ and let $B(R)\subset Y$ be an $R$-ball, which  contains a (small) neighbourhood $V_+0\subset X$ of $x_0$ and such that the boundary sphere $S(R)=\partial B(R)$ contains $x_0$. 
Show that:

$\bullet$ $R\geq r $, 

$\bullet$ if $R=r+\varepsilon$ for a small $\varepsilon\geq 0$, then  the sphere $S(R)$ is {\it smooth} at the point $x_0$, 

 $\bullet$ if  $S(R)$ is smooth at $x_0$, then   the {\it radial component} of the second fundamental form of $X$ at $x_0$  is greater than that of $S(R)$, 
$$\langle {\rm II}_X(\tau, \tau),\nu\rangle \geq \langle {\rm II_{S(R)}}(\tau, \tau)\nu\rangle, $$
where $\nu $ is the inward looking  unit normal vector to $S(R)$ at $x_0$ and $\tau \in T_{x_0}(X)\subset T_{x_0}(S(R).$
 
 (If the sphere  $S(R)$ is convex at $x_0$, then $0\leq \langle {\rm II_{S(R)}}_X(\tau, \tau),\nu|\rangle =\sqrt {{\|\rm II_{S(R)}}_X(\tau, \tau)\|}$.)

{\bf 3.C. Defintion of $+\rho$-Encircling $\bf X_{+\rho}=T^\perp_{\rho} (X)$  }
Given an immersed $X\hookrightarrow\mathbb R^N$
 let $T^\perp_{\rho} (X) \rightarrow\mathbb R^N$, $\rho>0$, be the {\it normal exponential} (tautological) map from the $\rho$-spherical normal bundle of $X$ to $\mathbb R^N$, where
 this "$\rho$-spherical normal bundle $T^\perp_{\rho} (X)$"  is the set of vectors normal to $X$ of length $\rho$.
  
   For instance if $X \hookrightarrow\mathbb R^N$ is an embedding and $\rho>0$ is small then the image of this map  
   is equal the  {\it the boundary of the $\rho$-neighbourhood of $X$}, denoted 
 {\color {blue}$$X_{+ \rho}= \partial U_{\rho}(X) =\{y \in \mathbb R^N\}_{dist (y, X)}=\rho.$$}

 In general, if $X \overset {f} \hookrightarrow\mathbb R^N$ is an immersion and 
 if $\rho<\big (curv^\perp X\hookrightarrow \mathbb R^n\big )^{-1}$
 then the exponential map is also an immersion and we abbreviate this by writing 
 $$X_{+\rho} \overset {f_{+\rho}} \hookrightarrow\mathbb R^N$$
  and observe the following.
  
 (a) If  $X \hookrightarrow\mathbb R^N$ is contained in $R$-ball, then  $X_{+\rho} \hookrightarrow\mathbb R^N$ is contained in the 
 $(R+\rho)$- ball, where the relation 
  $  \rho={r\over c}-\rho$, implies that
 {\color {blue}$$curv^\perp\big (X_{+\rho}\hookrightarrow  B^N(1)\big )=1/\rho \leq 1+2c=1=1+2\cdot curv^\perp (X).$$}

(b) The focal radius and the curvature of  $X_{+\rho}$ satisfy the following mutually equivalent relations 
 {\color {blue}$$rad^\perp(X_{+\rho}) =\min(\rho, rad^\perp(X)-\rho)\leqno{\color {blue} \bf[\rho]}$$} 
  and  
 {\color {blue}$$curv^\perp (X_{+\rho}\overset {f_{+\rho}} \hookrightarrow\mathbb R^N)
 =\max \Big(\rho^{-1},  \big (curv^\perp( X\hookrightarrow \mathbb R^n)\big )^{-1}-\rho\big)^{-1}\Big).\leqno{\color {blue} \bf[\rho^{-1}]}$$}
 
 { \bf 3.D.} {\textbf {$\bf [1+2c]$-Example.}} Let 
 $curv^\perp\big (X\hookrightarrow  B^N(1)\big ) \leq c$ and  move $X$ to the smaller   ball $B^N(r)$ by scaling $X\mapsto X'= r X$ for 
 $r=1-\rho$,  for some $0<\rho< 1/c$. 
 
 Then   $X'_{+\rho}$ is contained in the unit ball and 
 $$\Big (curv^\perp\big (X'_{+\rho}\hookrightarrow  B^N(1)\big )\Big)^{-1} \geq \min \left( {1\over \rho},\left({r\over c}-\rho\right)^{-1}\right ).  $$

 {\it Riemannian Remark.} The definition of $X_{+\rho}$ makes sense for an  immersed submanifold  $ X$ in  a Riemannian  manifold $Y$ 
 if the normal exponential map  
  $$\exp^\perp_ \rho:T^\perp_{\rho} (X)\hookrightarrow Y$$
is an immersion, e.g. if $Y$ is complete with non-positive sectional curvature $\kappa$ and also for $\kappa (Y)\leq 1$
and $\rho<\pi$.
Here the  above {\color {blue} $\bf[\rho]$} remains true but  {\color {blue} $\bf[\rho^{-1} ]$}  doesn't.

This suggests  that the reciprocal of the focal radius of $X$, may serve a replacement for the normal curvature for Riemannian submanifolds   
$$curv_x^{foc}(X\hookrightarrow Y)={1\over  rad_x^\perp(X\hookrightarrow Y))}.$$

$Y$ be a complete Riemannian manifold, $X\subset Y$ a smooth immersed submanifold and let us {\it define} the focal curvature of $X$ in $Y$ as the reciprocal of the focal radius of $X$,
  $$curv_x^{foc}(X\hookrightarrow Y)={1\over  rad_x^\perp(X\hookrightarrow Y))}.$$

\subsection { Maximum Principle by Exercises }


{\bf  Definition of  $\bf maxrad^\perp$.} 
 Let $Y$ be a metric space, let  $X\subset Y$ be a subset  and  let $x_0\in X$.
 
Define  $maxrad_{x_0}^\perp (X)$ as the {\it  infimum of the numbers $R$}, such that   there exists a point $y_0\in Y$ such that $dist(x_0,y_0)\leq R$ and 
the distance function $x\mapsto dist_Y(x,y_0)$ assumes  {\it local maximum} at $x_0$.

{\bf 3.1.A.} Consult (i) from the previous section and show that R$$maxrad_{x_0}^\perp (X)\geq  rad_{x_0}^\perp(X), $$ 
for smooth submanifolds $X$ in Riemannian manifolds $Y$.

{\bf 3.1.B.} Show that $$maxrad_{x_0}^\perp (X)= rad_{x_0}^\perp(X)\mbox { for  } dim(X)=1,$$ 
for  smooth 1-submanifolds (curves) $X$ in Riemannian manifolds $Y$, provided the normal exponential map $\exp:T^\perp{x_0}(X) \to Y$ is immersion on the $R$-ball $B^{N-n}_{0=x_0}(T^\perp_{x_0}(X)$.
Show    that the condition $ dim(X)=1$ is necessary.

{\bf 3.1.C.} Show that  if  a compact subset $X\subset Y$ 
 is {\it contained in an $R$-ball}  $B_{y_0}(R)\subset Y$, then  
$$\inf_{x\in X} maxrad_{x}^\perp (X)\leq R$$
 
 {\bf 3.1.D.} Show that the inequality $\inf_{x\in X} maxrad_{x}^\perp (X)\leq R$ remains valid for smooth  immersed {\it complete,    possibly non-compact,} submanifolds  $X\hookrightarrow Y$,  provided $curv^\perp(X)<\infty$.

(The condition $curv^\perp(X)=\sup_x curv^\perp(X)<\infty $  is necessary: there are  examples [Roz 1961]  of complete surfaces $X$   in the unit  3-ball  with negative Gauss curvatures, hence with $maxrad_{x}(X)=\infty$ for all $x\in X$.)

{\bf 3.1.E.}  Let $D(\rho)\subset B^N(1)\subset  \mathbb R^N$, $N\geq 3$ be the {\it boundary of the convex hull of  a  truncated unit ball}, where 
 $D(\rho)$ is equal to  the union of a spherical cap  $C^{N-1}(\rho)\subset S^{N-1}(1)=\partial B^N(1) $, $0< \rho<\pi $ and a flat $(n-1)$-ball $B^{N-1}(r=\sin \rho)\subset B^N(1)$,  
$$D(\rho)= C^{N-1}(\rho)\\cup B^{N-1}(r),$$
where $\rho$ is the  radius of $C^{N-1}(\rho)$ regarded as a ball in the spherical geometry in $S^{N-1}(R)$,
and where the (edge-like) intersection $E$ of the two parts of $D(\rho)$,
	$$E(r)=C^{N-1}(\rho)\subset S^{N-1}(R)\cap B^{N-1}(r)=(\partial C^{N-1}=\partial  B^{N-1}$$
is an $(N-1)$-sphere contained in $S^{N-1}(R)$ of (Euclidean) radius $r$.

Let $x_0\in E_r$ and
show that 

$\bullet_{conv}$ if $\rho \leq {\pi \over 2} $ then $maxrad^\perp_{x_0}(D(\rho))=r=\sin \rho$,

$\bullet_{concv}$  if $\rho \geq {\pi \over 2} $ then $maxrad^\perp_{x_0}(D(\rho))=R$.

{\bf 3.1.F.}  {\color {blue} Non-Smooth Maximum Principle.} Let $X\subset B^N(R)\subset  \mathbb R^N$ be a closed connected  subset in an $R$ ball,  such that 
$$maxrad^\perp_x(X)\geq r \mbox  { for some $r\leq R$ qnd all $x\in X$}.$$

Prove  that if $r=R$, then the intersection 
$X\cap \partial B^N(R)\subset S^N(R)=\partial B^N(R)$ is {\color {blue}non-empty} and 
show that   {\it  no connected component} of this intersection 
$X\cap \partial B^N(R)\subset S^N(R)=\partial B^N(R)$ is {\it contained in a spherical cap} 
$$C^{N-1}\left (\rho<{\pi\over 2}r\right )\subset S^{N-1}(R).$$
Consequently,  this intersection has no isolated points. 

Moreover, \vspace {1mm}

 {\it the topological dimensions of all  connected components  of $X\cap \partial B^N(R)$ satisfy
 {\color {blue}$$dim(comp(X\cap \partial B^N(R))\leq 1.$$}}

 {\bf 3.1.G.} Show that there  exists  a smooth convex (topologically spherical)  rotationally symmetric surface in the unit 3-ball,  
  $X\subset B^N(1)\subset  \mathbb R^3$,   
 which is {\it not equal to the boundary sphere} $S^2(1)=  B^3(1)$  and such that $maxrad^\perp_x(X)\geq 1$ for all $x\in X$. 
  
 {\bf 3.1.H.}  Generalize the above  to subsets $X$ (e.g. smooth submanifolds)  in balls $B(R)$   in Riemannian manifolds $Y$, 
  where the boundary of $B$, as well as the boundaries of concentric balls of radii $0<r\leq R$ are smooth and 
  where the inequality $\rho \geq {\pi \over 2}$ should be  be replaced by  $\rho \geq \delta=\delta(B)>0.$

 Thus show that if a compact connected subset   $X\subset B(R) $ satisfies $$ maxrad^\perp_x(X)\geq R$$ for all $x\in X$, then 
 
$\bullet$  the intersection   $X\cap \partial B^N(R)$ is non-empty,

$\bullet$   the   connected components  of this intersection satisfy

 $$dim(compX\cap \partial B^N(R))\leq 1,$$

  $\bullet$ no connected component of $X\cap \partial B^N(R))$ can be diffeomorphic to  segments $[0,1]$ and/or  $(0,1].$
 
  Consequently, \vspace {1mm}
 
 {\color{blue!44!black}  if {\it  all geodesics 
 $\gamma $} in a smoothly immersed  closed submanifold in a ball $B^N(R))\subset Y$ satisfy 
  $rad^\perp (\gamma\hookrightarrow Y)\geq R$, then   {\it $X$ is contained in the boundary $\partial B(R)$.}}
 \vspace {1mm}

{\bf 3.1.I.} Let $X$ be a  smoothly immersed {\it complete  connected} submanifold in a ball $B(R))\subset Y$, such that 
 the intersection $X$ with the boundary sphere $S(R)=\partial B(R)$ is{\it  nonempty} and such that each point $x_0\in X \cap S(R)$ admits a neighbourhood $X_0\subset X$ such that 
  radial component of the second fundamental form of $X$ at all $x\in X_0$  is non greater than that of the concentric sphere $S(r)$, which contains $x$,  
$$\langle {\rm II}_X(\tau, \tau),\nu\rangle \leq \langle {\rm II}_{S(r)}(\bar\tau, \bar\tau),\nu\rangle, $$
 where $\tau$ and $\nu$ and $\bar\tau\in T_x(S(r)$ is the normal projection of $\tau\in T_x(X)\subset T_x(Y)\supset T_x(S(r))$ to $T_x(S(r))$.
  
  Then $X$ is contained in the boundary  of the ball, $X\subset \partial B(R)$. 
  
  {\it Hint.} Prove convexity of a  $\phi(dist(x, S(R))$ for a suitable function $\phi(d)$.
   (Compare with  \SunshineOpenCircled \hspace {1mm} in section 1.)
  
  {\it\color {magenta} Question} What should be a comprehensive "non-smooth 'maximum principle", which would incorporate all we know in the smooth case? 
  



 



\subsection  {Topologically Defined  Focal Radius}


 Let  $X^n\subset \mathbb R^{n+1}$ be a   smooth hypersurface, let $x\in X$ and observe that 
 the normal curvature $curv_x^\perp(X)$
 is equal to the infimum of the  curvatures $c$ of the spheres  $S_\pm^n(1/c)$,
which  are:\vspace {1mm}

(i) tangent  to $X$ at $x$,  
 
(ii) the balls bounded by these spheres do not intersect (small)  neighbourhoods of $f(x)$ in $f(X)$ minus $f(x)$ itself, 
   
(iii) do not mutually  intersect away from $x$.\vspace {1mm}


Generalise this  to  smooth submanifolds $X^n\subset \mathbb R^N $ for all $N\geq n+1$ as follows.

 Let $\mathcal B(c)$
be a family of  balls $B^N_y(1/c)\mathbb  \mathbb R^N $ with centers $y\in \mathbb R^N $ such that

(i$'$) all balls from $\mathcal B(c)$ contain $x\in X$, 
 
 (ii$'$) 
 the balls do not intersect (small)  neighbourhoods of $x_0$ in $X$ minus $x_0$ itself, 

 (iii$'$) for all $\varepsilon >0$, there exists a family of points in $\mathbb R^N$ continuously parametrized by $\mathcal B(c)$, say
 $$ \phi_\varepsilon :\mathcal B(c)\ni B\to \mathbb R^N, $$
 such that

$\bullet$ $\phi_\varepsilon (B)\in B$  for all  $B\in \mathcal B(c),$
 
$\bullet$  $dist (\phi_\varepsilon (B), x_0)\leq\varepsilon$ for all  $B\in \mathcal B(c),$

 $\bullet$ the set $\mathcal B(c)$ contains an $(N-n-1)$-cycle the $\phi$-imagef this cycle. is non-trivilally linked with $X$  for all sufficiently  small $\varepsilon$.\footnote {Think of $X$ as a relative  $n$-cycle 
 in the pair $(B_{x}^N(2\varepsilon),\partial (B_{x}^N)(2\varepsilon))$. } 
 
Then show that  
 $curv_{x}^\perp (X)$ is equal to

$\CIRCLE$ {\it the infimum of $c>0$, such that a family $\mathcal B(c)$ with all these properties exists.}
 

\vspace{1mm}

  The $\CIRCLE$-definitions of $c$  and $\inf c$ apply to {\it non-smooth   topological} submanifolds   $X\mathbb R^N$. Yet   if this $\CIRCLE$-"curvature"   $curv_x^{\small\CIRCLE}(X)=1/\inf c$
 is uniformly bounded for  $x\in X$,  then $X$ is $C^1$-smooth, moreover, it is $C^{1,1}$-smooth--the partial derivatives are  {\it Lipschitz}.

 {



\section { Products of Spheres in  $B^{n+1}$ with Small Curvatures}


 { \bf 4.A. Products of Spheres Represented by Hypersurfaces} \vspace {1mm}
{\sl}  Let  $X$ be a product of $m$ spheres and $k\geq m-1$.  Then $X_m\times S^k$ admits a codimension one embedding  to the unit ball with normal 
   curvature $1+2\sqrt m.$

{\it Proof.} Imbed   $X$ to $B^{N+m}(1)\subset  B^{N+k+1}(1)$ for $N=dim(X)$  with curvature $c=\sqrt m$ (see  1.A),   let $\rho= 1+2\sqrt m$ and observe that  
$X'_{+\rho}\subset B^{N+k+1}(1)$, (this  is the boundary of the 
$\rho$-neighbourhood of  $X'\subset B^{N+k+1}(1)$ in the present case) 
is diffeomorphic to $X\times S^k$.
Since $curv^\perp(X'_{+\rho})\leq 1+2c$ (see 3.D), the proof follows.\vspace {1mm}

{\sc Two  Examples and one Theorem.}\vspace {1mm}

$(\bullet_1)$ {\sl Products of two  spheres 
   admit codimension one embeddings  to the unit balls with normal 
   curvatures 3:} 
 {\color {blue}$$curv^\perp \big ( S^{n_1}\times S^{n_2}=S_{+1/3}^{n_1}(2/3) \subset B^{n_1+1+n_2}(1)\big)=3,\leqno{\left[2/3]\times [1/3\right]}.$$}
 
 $(\bullet _2)$  
{\sl  Products of three spheres $S^{n_1}\times S^{n_2}\times S^{n_3}$, e.g. 3-tori $\mathbb T^3$,  admit  codimension one embeddings  to the unit balls with curvatures $1+2\sqrt 2 <4$.}\vspace {1mm}
 
 We {\color {magenta} don't know} answers to the following questions:\vspace {1mm}

 {\sl are there  immersions 
 $ S^{n_1}\times S^{n_2}\hookrightarrow B^{n_1+n_2}(1) $ with $curv^\perp  <3$?

 are there  immersions  immersions $S^{n_1}\times S^{n_2}\times S^{n_3}\hookrightarrow B^{n_1+n_2+n_3}(1) $ with $curv^\perp < 1+2\sqrt 2$.}\vspace {1mm}

But   the situation changes starting from $m=4$ and $C=1+3\sqrt 2=5.24264....$ with the following.

  {{ \bf 4.B. Codimension one Immersion  Theorem.} {\sl Let $X$ be a compact orientable $n$-manifold, which admits an immersion to $\mathbb R^{n+1}$, e.g.   $X$ is (diffeomorphic to) a product of spheres $S^{n_i}$ of dimensions $n_i$, 
  $\sum_in_i=n$.

Then, for all $\varepsilon>0$,  the  product $S^{20n^2}\times X$  admits an   immersion $f_\varepsilon $ to the $(20n^2+n+1)$-ball, such that  
   {\color {blue}$$curv^\perp\big (( S^{N}\times X)\overset {f_\varepsilon}\hookrightarrow B^{20n^2+n+1}(1)\big)\leq1+2
  \sqrt \frac {3(n+1)}{n+3}+\varepsilon <4.5.\leqno {(\bf <4.5)}$$}}

{\it Proof.}  The $\sqrt 3$-immersion corollary 1.C with $m=n+1$  delivers an immersion 
$X\to B^{20n^2}(1) $  with $curv^\perp \leq \sqrt \frac {3(n+1)}{n+3}+\varepsilon$  and  the manifold $X'_\rho$ as  in  $[1+2c]$-example (1.G) does the job  since it is diffeomorphic to $X\times S^{20n^2}$ in the present case.

{\bf $[X=\mathbb T^n]$-Case.} If  $N>>n$, then    the $\sqrt 3$-Clifford sub-torus  theorem  {\bf  1.C}  implies that
$S^N\times \mathbb T^n$  admits an   immersion to the $(N+n+1)$-ball, such that  
   {\color {blue}$$curv^\perp\big (( S^{N}\times X)\overset {f_\varepsilon}\hookrightarrow B^{N+n+1}(1)\big)\leq1+2
  \sqrt \frac {3n}{n+2}.$$}

{\it Embedding Remark.} Unlike how it is in   $(\bullet_1)$ and  $(\bullet_2)$, the 
construction of $f_\varepsilon$ in    1.K creates self-intersection 
of $S^k\times X$ in the ball.

{\it\color {magenta} Sharpness Conjectures.} The constant $1+2
  \sqrt \frac {3n}{n+2}$, {\color {magenta}probbaly},  is optimal for tori $\mathbb T^n$ 
of dimension $n\geq 3$

We also {\color {magenta} conjecture }  that there are  
{\it no embeddings} $\mathbb T^n\times S^k\to B^{n+k+1}(1)$
with  $curv^\perp \leq1+2 \sqrt \frac {3n}{n+2}+\varepsilon $
for all $n\geq 3$ and $\varepsilon < 1/n^2$.\vspace {1mm}
   
But it is {\color {magenta}hard to say } if the constant  $\sqrt \frac {3(n+1)}{n+3}$  for general  orientable $X^n\mathbb R^{n+1}$ can be improved, even to $\sqrt \frac {3n}{n+2}$. 

Also it is {\color {magenta} unclear} what to expect in this regard from {\it non-orientable immersed} hypersurface $X^n\hookrightarrow \mathbb R^{n+1}$ 

\vspace {1mm}
    
    {\it   Products of Equidimensional Manifolds. } The codimension one immersion  theorem
  doesn't deliver  immersions of products  of equidimensional manifolds with "interesting" curvature bounds, while by arguing as in $(\bullet _1)$ and $(\bullet _2)$ we   show the following.

    $(\bullet _3)$  {\sl The product of $(m+2)$ copies of $S^m$   admits an embedding to the ball  $ B^{m(m+2)+1}(1)$ with $curv^\perp\leq 1+2\sqrt {m+1}$. }
  
  For instance, (as in $\bullet_2$) the $3$-torus  embeds to the unit 4-ball, such that 
    {\color {blue}    $$curv^\perp\big (\mathbb T^3\subset B^4)\big)\leq 1+2\sqrt {2}<4.$$ }

   {\color {magenta} Conjectrally,} {\sl the constant $1+2\sqrt {m+1}$ is  optimal  
   for all $m=1,2,...4$, possibly, not only for embedding  but also for immersions  $$(S^m)^{m+2}\hookrightarrow  B^{m(m+2)+1}(1).$$  }

  \section {Extremality,  Rigidity, Stability:      Spheres and  Veronese Varieties }


 The natural candidates for {\it extremal immersions} $X\hookrightarrow B^N$, which  implement maximal topological complexity with minimal curvatures are the most symmetric ones that are immersions,  which are equivariant under large isometry groups  $G$ acting on $X$ and $B^N$
    
 For instance the standard ($O(n)$-equivariant)  embedding $S^n \hookrightarrow B^N$ is extremal by  \SunshineOpenCircled \hspace {1mm}Example. 1.A:
{\sl  
 the $n$-dimensional spheres of radius one, are the  only closed submanifolds with $curv^\perp\leq 1$ in $B^N(1)$}.
  (If $n=1$ these may be  multiple  coverings of the circle).

    {\it Rigidity and Stability}.   Most (all?) sharp geometric inequalities are accompanied by  
    the rigidity/stability of the extremal objects\footnote {See stability Gr. for a general discussion}.
   
    To  establish stability for $S^n\subset B^N$ we  observe that the maximum principle argument  used in 1.A. equally applies to {\it complete} (for the induced Riemannian metrics) manifolds $C^{1,1}$-immersed to $B^N(1)$ and that   
     the space of $C^{1,1}$-|immersions $curv^\perp\leq const$ of immersed  complete manifolds to the ball  is compact. 
     
     Thus we conclude that there exists $\varepsilon>0$, such that   if a closed immerses  submanifold satisfies  
     $$curv^\perp \big(X^n\hookrightarrow B^N(1)\big)\leq 1+\varepsilon \mbox {  and  } n\geq 2,$$  then  $X$ can be obtained  by a $\delta$-small  $C^1$-perturbation of a unit $n$-sphere  $S^n\subset B^N$,   where $\delta\to 0 $ for $\varepsilon \to 0$. 
 
 A priori, this $\varepsilon$ could  depend on $n$ and $N$,  but when this argument is applied to geodesics in $S^n$,
 it  provides an effective, albeit rough, bound on $\varepsilon$,  e.g.  $\varepsilon =0.01$  (See  below and section 12 for  Petrunin's sharp result.) 
    
  {\it Immersions to Tubes}. The maximum principle applied to closed immersed $n$-submanifolds  in "unit tubes" $B^N(1)\times \mathbb R^k\subset \mathbb R^{N+k}$    shows 
  that 
     $$curv^\perp \big(X^n\hookrightarrow B^N(1)\times \mathbb R^k \big)\geq 1 
     \mbox {  for $k\leq n+1$},$$
   where  extremal $X$, i.e. where $curv^\perp \big(X^n\hookrightarrow B^N(1)\times \mathbb R^k =1$  for $k\geq 1$ are by no means   unique. (see section 11)

    {\it About  Mean Curvature.} The maximum principle argument also applies to immersed $n$-submanifolds $X$ in $B^N$ with $mean.curv\leq n-1$  and shows that 
    these $X$ lie in $S^{N-1}$, where they  are {\it minimal}, i.e.  have zero mean curvatures.
    
    An abundance of minimal surfaces in $S^{N-1}$ makes it  {\color {magenta}plausible} that all   
 $n$-manifolds  admit 
   and $curv^\perp(X)\leq const =const(n)$, say $const(n)=100^n$,
for given   $\varepsilon,\delta>0$. \footnote {Possibly the case $n=2$   can be approached with the techniques  from  [Nadirashvili 1996] and its generalizations.}
      
    \vspace {1mm}

    {\bf 5.A. Veronese Manifolds.} Besides  $n$-spheres, there are  other  $O(n+1)$-equivariant immersion $S^n\hookrightarrow B^N(1)$, where the most interesting ones are the (quadratic)  {\it Veronese maps}.
    
     These are (minimal) isometric immersions of the $n$-spheres of radii $R_n=\sqrt\frac {2(n+1)}{n}$      to the unit balls, which  factors through embeddings of the projective spaces   $\mathbb RP^n=S^n(R_n)/\{\pm1\}$ to the balls 
    $B^\frac {m(m+3)}{2}$, where these embedding have {\it amazingly small} curvatures:

 $$curv(Ver_n)=curv^\perp\left (\mathbb RP_{Ver}^n \hookrightarrow  B^\frac {n(n+3)}{2}\right)=\sqrt \frac {2n}{n+1}, \mbox{  e.g.} $$
   $$curv(Ver_2)=curv\left (\mathbb RP_{Ver}^2 \hookrightarrow  B^5\right)=2\sqrt \frac {1}{3}<1.155,$$

Observe that  the radii  $R_n$ of the Veronese   $n$-spheres, which covers $\mathbb RP_{Ver}^n$,    
satisfy
$$R_n={2\over curv(Ver_n)}.\leqno {[2/curv^\perp]}$$

{ \color {magenta} Conjecture.} $$curv^\perp(X^n,B^N)< \sqrt \frac {2n}{n+1}\implies X=_{diffeo} S^n.$$ 

The "homeo-version" of this  {\color {blue} proven by Petrunin for   $n=2$}. (See [Petrunin 2023] and  section 12, where we also  explain the above  and say more about Veronese maps and their generalizations.)
  


  \section{ Exercises: Hypersurfaces Inscribed in Convex Sets}


  Given a subset $V\subset \mathbb R^{ n+1}$, let $ext_{+r}(V)$
denote the {\it $r$-neighbourhood of $V$}, that is the subset of points in $\mathbb R^{n+1}$ within distance $\leq r$ from $V$.
  $$ext_{+r}(V)=\{y\in \mathbb R^n\}_{dist (y, V)\leq r}\subset \mathbb R^{n+1},$$
   and let 
  $$int_{-r}(V)\subset  V$$ 
  be the complement of the interior of  the $r$-exterior of the   complement of $\mathbb R^{n+1}\setminus V$, that is equal to 
  the  set of points in $V$ with   distance $\geq r$  from  the boundary of  $V$,
  $$int_{-r}(V)=\{v\in V \}_{dist (v,\partial V)\geq r}\subset  V.$$
  Clearly, 
 $$ext_{+r}(int_{-r}(V))\subset  V\mbox { and } int_{-r}(ext_{+r}(V))=V.$$
 
  Let $R=R(V)$ denote    the {\it  in-radius  of} $V$, that is   the maximal  distance  from  the boundary of  $V$ in $V$,
   $$R=inrad (V)=\sup _{v\in V} dist(v,\partial V)$$
  and let 
  $$cntr(V)= int_{-r}(V)$$
 be the set of the centers of the $R$-balls in $V$, 
that is the subsets of $v\in V$ with $dist(v, \partial V)=R=inrad (V).$

  Let $V\subset \mathbb R^{ n+1}$ be a  {\it compact  convex domain},  e.g. the  $(n+1)$-cube $\square ^{n+1}=[-1,1]^{n+1}$
   or an $(n+1)$-simplex $\triangle^{n+1}$.

     Then, clearly,  the  {\it $r$-interior }of $V$ is convex and if $r=R=inrad(V)$
then $int_{R}(V)$ called  {\it the central locus} in $V$,
$$int_{R}(V)=cntr(V)$$ 
is a {\it non-empty  compact {\it convex} subset in $V$ of dimension $\leq n=dim(V)-1$}.

For instance, if $V$ is a cube or a simplex, then $cntr(V)$ consists of a single point 
and $ext_{+r}(int_{-r}(V))$ is equal to the (unique maximal)  ball inscribed into $V$.

 (If $V$ is a general $(n+1)$-dimensional  rectangular solid then  
  $int_{-r}(V)$ is a subsolid of   certain  dimension $0, 1,..., n$.)
 
  {\bf 6.A.  $inrad$-Convex Exercises.}    (a) Let the boundary of  $V$ be {\it $C^{1,1}$-smooth}\footnote{ Locally, the hypersurface 
    $\partial V\subset \mathbb R^{n+1}$ is representable by the graph of a $C^1$-function with bounded measurable second derivatives.} (e.g. piecewise $C^2$-smooth) with curvature bounded by a constant $c$,
    $$curv^\perp(\partial V\subset \mathbb R^{n+1})\leq c.$$
  
   Show that if $r\leq {1\over c}$, then, 
    the $r$-balls  $B\subset \mathbb R^{n+1}$ tangent  to $\partial V$ either are  fully contained 
    in $V$ or lie outside $V$, meeting $W$ at a single contact point between  the boundaries of  $B$ and  $ V$; consequently:
   $$ext_{+R} (cntr (V))=V\mbox {  for  }  R=inrad (W).$$

 (b) {\sl Let $X\hookrightarrow \mathbb R^{n+1}$ be a $C^2$-smooth compact  
    immersed  hypersurface} in $\mathbb R^{n+1}$ and let 
    $$W=conv(X)$$
  be the {\it convex hull}   of (the image of)  $X\hookrightarrow \mathbb R^{n+1}$. 
  
  Show that the boundary of  $W$ is {\it $C^{1,1}$-smooth}\footnote{ Locally, the hypersurface 
    $\partial W\subset \mathbb R^{n+1}$ is representable by the graph of a $C^1$-function with bounded measurable second derivatives.} with curvature bounded by that of $X$, 
    $$curv^\perp(\partial W\subset \mathbb R^{n+1})\leq curv^\perp(X\hookrightarrow \mathbb R^{n+1}).$$

  (c) {\it Sphericity}. Let $V\subset \mathbb R^{n+1}$ be a convex bounded domain, e.g. a polytope, such as   $(n+1)$-cube $\square ^{n+1}=[-1,1]^{n+1}$
   or an $(n+1)$-simplex $\triangle^{n+1}$, and let $X\overset {f}\hookrightarrow V$ be a $C^2$-smooth immersion, where $X$ is a closed $n$-manifold. 
   
   Apply (a) and (b) to the convex hull  $W=conv(X)\subset V$ of $X$ and 
   show that if 
   $$inrad(V)=R\leq {1\over curv^\perp (X\hookrightarrow V)},$$
  then, in fact, 
   $$inrad(V)={1\over curv^\perp (X\hookrightarrow V)}.$$
   
  Fithermore, if $cntr(V)$ consists of a single point $o\in V$, (e.g. $V =\square ^{n+1}$ or $V =\triangle ^{n+1}$), show that 
 the image of  the immersion $f$ is contained  the $R$-ball centred at $o$  for $R=inrad V$.
  
   Consequently, (see  1.3.A)
   \vspace {1mm}
  
  \hspace {9mm} {\it the image of $X\overset {f}\hookrightarrow V$  is equal to the $R$-sphere centered at $o\in V$.}
    \vspace {1mm}

 (d) {\it  Stability.}  Argue as  in  section 1.3 and,  assuming as above that  $cntr(V)$ consists of a {\it single point} $o\in V$,  show that the (only)  $R$-sphere in $V$ is stable: 
 
 {\sl there exists  an $\varepsilon =\varepsilon (V)>0$, such that all immersed closed 
  hypersurfaces $X\hookrightarrow V$ with $curv^\perp(X)\leq R+\varepsilon$ are $\delta$-close in the $C^1$-topology to to the $R$-sphere $S_o^n(R)$, where $\delta\underset{\varepsilon \to 0} \to 0.$}
 \vspace {1mm}

 {\it More on Stability}. 
Unlike  1.3.B, this  $\varepsilon$ is sensitive to dimension.

For instance if $V$ is the regular  unit simplex then   $\varepsilon (\triangle^{n+1})\sim\varepsilon_0/ n$ and if it is the cube  $\square^{n+1}=[-1,1]^{n+1}$, then 
$\varepsilon (\square^{n+1})\sim\varepsilon_0/\sqrt n.$

 \vspace {1mm}

{\it On $dim(cntr(V)>0$}. If {\it $int_{-r}(V)$ has positive dimension}, then there are many 
non-spherical $C^2$-immersed (and even more   $C^{1,1}$)   hypersurfaces in $V$ with curvatures 
$\leq {1\over inrad V}$, see  section 11.


{\bf On $
\bf dim(cntr(V))>0$}. Let    $dim(cntr(V))=k>0$, let $Z\subset \mathbb R^{n+1}$
be an affine  $k$-dimensional subspace which contains the (convex!) subset $cntr(V)
\subset \mathbb R^{n+1}$ and let $B^{n+1}_Z(R)= ext_{+R}(Z)\subset \mathbb R^{n+1}$ be the $R$-neighbourhood of $Z$.

   \section { Bowl Inequalities}


     Let  $X\hookrightarrow \mathbb R^N$ be  an immersed complete (e.g. closed) connected  $n$-dimensional submanifold in the Euclidean $N$-space,
     let $x_0\in X$, let $T=T_{x_0}\subset \mathbb R^N$ be the tangent space to $X$ at $x_0$ (represented 
     by an affine subspace in $\mathbb R^N$) and let $P_{x_0}:X\to T_0^n$ be the normal projection map.

  Let $U_{x_0}\subset X$ be
       the   maximal connected neighbourhood of $x_0$,   such that the normal projection $P=P_{x_0}$, from
       $ U_{x_0}$ to $T_{x_0}$
    is a one-to-one diffeomorphism onto a domain $V_+{x_0}\subset T_{x_0}=\mathbb R^n$, which is {\it star 
       convex} with respect to $x_0.$
       
       Clearly such a $U_{x_0}$ exists and unique.  
       
     Let $\underline S=S^n(R)$  an  $n$-sphere  of radius $R$, which is tangent to $X$ at the point $x_0$,  (such spheres $\underline S=\underline S_\nu$ are parametrised by the  unit  normal vectors $\nu\in T_{x_0}^\perp(X)$) let $\underline P:\underline S\to T_{x_0}$, be the normal projection map and observe that the corresponding neighbourhood 
     $\underline U_{x_0}\subset \underline S$ is  the hemisphere $\underline S_+$  that is  the ball $B_{x_0}({\pi\over 2} R)\subset \underline S= S^n$ around $x_0$. 

Let $d(x)=dist_{T} (P(x),x_0)$ and let $\underline d(s)=dist_T(\underline P(s),x_0)=R\sin \frac {1}{R}dist_{\underline S}(s,x_0)$  be the corresponding  function for the sphere $\underline  S$.

 Let $h(x)=dist (x, y=P(x)) $, $x\in X$,  and let $\underline h(s)=R\cos \frac {1}{R}dist_{\underline S}(s,x_0)$ be the corresponding 
  function for the sphere $\underline  S$. 
    
{\it Remark.} Both $d$--functions and both $h$-functions have their gradients   {\it bounded by one}, 
    in fact, 
    $$\|grad_{\underline S}(d(s))\|^2+\|grad_{\underline S} (\underline  h(s))\|^2=1\mbox { and } \|grad_{X}(d(x))\|^2+ \|grad_{X}( h(x))\|^2\leq1, $$
    
    The gradients of  both $d$-functions have unit norms at $x_0$,\footnote{These  functions are non-differentiable at $x_0$ but the norms of their gradients continously extend to one at $x_0$.}, they don't vanish   in the interiors of the domains $U_{x_0}$ and $\underline U_{x_0}$ correspondingly;  $grad (\underline h)$ vanishes on the boundary. of $\underline U_{x_0}$ and $\underline U_{x_0}$ vanishes at at least 2 ponts at the boundary   of $U_{x_0}$.

    The gradients of  the $h$-functions have $norms < 1$ in (the interiors of)   domains $U_{x_0}$ and $\underline U_{x_0}$  correspondingly, and these norms. are equal to one  the boundaries of these domains.

       In fact, $\underline U_{x_0}$ is the same as   the   maximal connected neighbourhood of $x_0$,   where $\|grad_X (h)\|< 1$
   and 
  the $P$-image   of which is star convex.

       The following proposition, says that $U_0$ lies at least  as close in the $C^1$-metric  to $T_{x_0}$ as 
   $\underline S_+$.

     {\bf 7.A.  Hemisphere Comparison Inequalities.} Let $$curv^\perp(X)\leq curv^\perp(\underline S)=1/R. $$}
      Then:

      The gradient of the  $h$-function on $X$, 
    $$h:x\mapsto dist (x, P(x)) $$ for  $x\in U_{x_0}$ is
      bounded by that for the $\underline h$-function on $\underline S_+$
     { \color {blue}$$  \|grad (h)\|\leq \|grad( \underline h)\|\mbox {  for  }dist_X(x,x_0) \leq  dist_{\underline S} (s,x_0) \|\mbox 
      { and } s\in \underline S_+$$}
     Consequently, the domain  $U_{x_0}\subset X$ {\it contains an open $R$-ball} centered at $x_0$.
     
      $\bullet_d$ 
      The gradient  of the  $d$-function on $X$  in the radial direction is  bounded from below by  that for $\underline d$:
      
     if $s\in \underline S_+$ and a unit vector $\tau \in T_x(X) $  which is tangent to a geodesic segment  $\gamma$ in $X$
      issuing from $x_0$ and termnating at $x$ satisfy
     $$length (\gamma) \leq dist (s,x_0)$$ 
     then 
    $$\langle grad(d),\tau\rangle \geq \|grad (\underline d(s))\|,$$
     Consequently, the $P$-images in $T$ of the $r$-balls from  $ U_{x_0}\subset X$ centered at $x_0$,  contain the $\underline P$-images    of the corresponding  spherical balls from $\underline S_+$,
     $$P(B_{x_0}(r ))\supset B_{x_0}\left (R\cdot \sin {1\over R} r\right )\subset T\mbox  { for all 
     $r\leq {\pi\over 2}R$,}$$
 
     $\bullet^{-1}_h$ The inverse function $h^{-1}(y)$, $y\in B_{x_0}(R)\subset T$, and the norm of its gradient are bounded by  $\underline h^{-1}(y)$, and  $\|grad (\underline h^{-1}(y))\|$ correspondingly.

    {\bf 7.B.   Out of Ball Corollary.} Let $\underline B^N(R)\subset \mathbb R^N$ be a ball, such  that the  boundary sphere 
    $\underline S^{N-1}(R)=\partial \underline B^N(R)$ is tangent to $X$ at $x_0$, i.e. 
    $$ T_{x_0} (\underline S^{N-1}(R))\supset T_{x_0}(X).$$
    
    If $curv^\perp(X)\leq 1/R$,  then the  subset $U_{x_0}\subset X$ doesn't intersect the interior of this ball. Thus,  $U_{x_0}$ lies in the closure of the complement of the union of the  $R$-balls tangent to $X$ at $x_0$.
     \vspace {1mm}
    
    {\bf 7.C. Spherical Bowl Theorem.} Let  $U_{x_0}(+r)\subset X$ be the $r$-neighbourhood of $U_{x_0}$ in $X$. Then 
    the gradient of the function $d(x) = dist_T (P(x), x_0)$  doesn't vanish  in the interior of the complement 
        $U_{x_0}(+R)\setminus  U_{x_0}$ 
    and 
    the $P$-mage  of the complement  $ U_{x_0}(+r)\setminus  U_{x_0}$, $r\leq R$,  doesn't intersect the interior of the ball $B_{x_0}(R-r)\subset T$. 
    
    {\it Proof} The bounds on the gradients of the functions $h$ in the hemisphere comparison inequalities
    follow from  the angular arc inequality. 1.1.E,
    while the bowl theorem follows from these inequalities applied to $X$ at $x_0$ and at all  ponts $x\in \partial U_{x_0}$.

 If $dim(X)=1$,  then the  bowl theorem, where the   proof\footnote {Hopf Schimd}
  becomes especially  transparent  \footnote {he abive Hopf Schimdt, oter proofs}  implies the following.
   
    {\bf 7.D. Circular Bow Inequality {\color {blue}$[dist_X\geq dist_A]$}.}\footnote{According to  [Sch 1921] this goes back H. A. Schwarz, 1884. }  Let a planar circular arc  $A\subset \mathbb R^2$  (a segment 
   of a circle)  and a smooth spatial curve $X\hookrightarrow \mathbb R^N$ satisfy:
 $$length  (X)=length(A)=l\mbox { an ]d } curv^\perp(X)\leq curv^\perp(A).$$
Then the  distance between the endpoints  of $X$ is {\it greater than or equal to that}  in $ A$,
where the equality holds if and only if $X$ is {\it congruent to $A$.}\vspace {1mm}

{\it Remark}   This inequality applied to geodesic segments from $X$ yields  most essential geometrical  properties $X$
with no use of the full ($n$-dimensional)  spherical bowl theorem.

  {\it Example.} Let $X\subset\mathbb R^N$ is a closed curve of length $2\pi$. If $curv^\perp(X)\leq 1$,   then the  Euclidean distances between opposite points $x,x_{opp}\in X$ are $\geq 2$, where an equality  $dist(x_0,(x_0){opp}=2$ implies that 
 $X$ is circlular. 
 
 {\it Exercises.} (a) Show that all closed curves of length $2\pi$ in the Euclidean space contain  pairs of opposite points $x,x_{opp}\in X$,n(i.e. with the $X$-distance $\pi$ between them),    such that
 $dist(x,x_{opp})\leq 2$.


         .

(b) Let  $X\subset\mathbb R^N$  be a smooth  simple\footnote {"Simple" means  {\it homeomorphic to  $[0,1]$}.}   curve, such that  the   two ends of $X$, say  $x_1,x_2\in X$, are  positioned in two parallel hyperplanes,
$$x_1\in H_1\subset\mathbb R^N\mbox { and }  x_2\in H_2\subset\mathbb R^N,$$
where $X$  is {\it tangent} to these hyperplanes at {\it the points $x_1,x_2$}.

 If  {\it no tangent line} to  $X$, except for  these  at the points $x_1 x_2$  is parallel to $H_i$, $i=1,2$, 
if there exists a point $x_\circ\in X$, such the tangent to $X$ at $x_\circ$ is {\it normal to} $H_i$
and if
$$curv^\perp(X)\leq c,$$ 
then the distance between the hyperplanes is at least $\pi/c$,
$$dist(H_1,H_2)\geq \pi/c,$$
where the equality implies that $X$ is a  circular arcs with curvature $c$ and  length $\pi/c$.

(c){\it Remark}.   If 
we allow {\it $C{1,1}$-curves},  then the extremal $X$  are {\it unions of pairs} of circular arcs of length $\pi/2c$,  where the (normal) curvature (vector function) of $X$, may be {\it discontinuous at the middle} of $X$.

Moreover, if  also allow $X$ to be tangent   $ H_i$  at some intervals  $X_i\subset X$ around the ends $x_i$ and  nowhere else (the subset $X\parallel subset X$, where the  of tangent lines to $X$  are parallel to $H_i$, consist of {\it  two}
 components), then gthe inequality still $dist(H_1,H_2)\geq \pi/c,$ holds valid  and  the extremal $X$, where the equality holds, contains,  
 besides two  circular $\pi/2c$-arcs, two straight segments that  are $X_i\subset H_i$, $i=1,2.$


 \subsection {High Dimensional  Applications of the  Circular   Bowl inequality}


   Basic  geometry properties  of immersed $n$-submanifolds $X$ in  Euclidean spaces  with  
     $$curv^\perp (X\overset{f}\hookrightarrow \mathbb R^N)\leq c,$$
   can be reduced to the case $n=1$  applied  
 to the geodesic segments from $X$, which are, by the definition of the normal curvature $curv^\perp (X)$, are  {\it curves in  $\mathbb R^N$ with $curv^\perp \leq c$.}



The circular bow inequality applied to geodesics in immersed $n$-submanifolds 
$$X\overset {f}\hookrightarrow \mathbb R^N, \mbox { } dim(X)=1,2,...n,...$$ yields the following.

{\bf 7.1.A. $\bf [2\sin]_{bow}$ and  $2\sin]_{dist}$ Inequalities.}  Let $\gamma\hookrightarrow X$ 
be an (not necessary minimising) {\it geodesic segment}\footnote{Recall "Geodesic" refers to the  induced (inner) Riemannian metric in X,}  between two points   $ x_0,x_1\in X$.
If the normal curvature of $X$ is bounded by $1/R$ and if  $length(\gamma)=l\leq 2\pi R$, then 

Then  the Euclidean  distance between these points is bounded from below: 
 $$dist_{\mathbb R^N}\big(f(x_0),f(x_1)\big)\geq 2R\sin {l\over 2R}\leqno {[2\sin]_{bow}}$$
and,  
the equality implies that the $f$-image of $\gamma$ is a circular ark in a plane in $\mathbb R^N$.

  If $X$ is connected and the induced metric in $X$ is {\it complete} (e.g, $X$ is compact without boundary), then
    {\sl $$dist_{\mathbb R^N}{\mathbb R^N}\big(f(x_0),f(x_1)\big)\geq 2\sin \left({dist_X(x_0,x_1)\over 2}\right) \leqno {[2\sin]_{dist}}$$
   for all $x_0,x_1\in X$, such that $dist(x_0,x_1) \leq 2\pi$.}
  
   {\it Remarks.}  The ${[2\sin]_{dist}}$-inequality for infinitesimally close points $x_0,x_1$ is equivalent to the inequality $curv^\perp\leq R$.

     The ${[2\sin]_{bow}}$-inequality holds for immersions to (complete simply connected) manifolds $Y$  with {\it non-positive} sectional curvatures. ( see ..)

     \vspace{2 mm}
.
 
   {\bf 7.1.B  $2\pi$-Injectivity } Let $TB_x(r)\subset T_x(X)$ gent space  be the $r$-ball in the tangent space at a point $x\in X$  and let  $\exp_x:TB_x(r)\to X$ be the exponential map.  If 
  $r<\pi$, then the composition of this map with our immersion $f:X\hookrightarrow
  \mathbb R^N$ is one-to one.

Here are two obvious sub-corollaries.

 {\bf  $2\pi$-Geodesic Loop Inequality.} Geodesic loops $\gamma$ in $ X$ have 
 $length(\gamma)\leq 2\pi$.
 
 {\bf $2\pi$-Diameter Inequality.} If the {\it intrinsic diameter}, i.e. the diameter   with respect to to the induced Riemannian metric, of $X\overset {f}\hookrightarrow \mathbb R^N$, satisfies 
$$diam_{int}(X)< 2\pi,$$
then $X$  is {\it  embedded} to $\mathbb R^N$:  the map $f$ is one-to-one.

 {\it This inequality is sharp}:  \ the equality holds for 
 $S^n_{Ver}(R_n)\to \mathbb RP_{Ver}^n  \hookrightarrow  B^{n(n+3\over 2}(1)$ by the  above $[{2\over curv^\perp}]$
  $$diam_{int}(S^n_{Ver})=\pi R_n=
 {2\pi\over curv^\perp(S^n_{Ver})}.$$
 
   {\it \color {magenta} Question.} Are Veronese the only ones with this property? (Compare with [Petr 2024] and also with  section 12)

{\it  Exercise.} Let   $X^n\overset {f}\hookrightarrow B^N(R)$ be an immersion, such that all geodesics in $X$ issuing from a point $x_0\in X$ have their Eucldiean  curvatures $curv^\perp$ bounded by one.

 Show that the image $f(X)\subset B^N(R)$ is  equal to an equatorial $n$-sub-sphere in $S^{N_1}(R)=\partial B^N(R)$ and if $X$ is connected and $dim(X)\geq 2$ then  the immersion $f$ is an embedding.


\section {Optimal Control, 
 Bow Theorem and Arm Lemma}


The circular bow inequality  {\color {blue}$[dist_X\geq dist_A]$} from section 7 represents the solution  of the following variational problem  for 
 curves, now denoted $y(s)$, $ s\in [0,l]$,    of length $l$
 issuing from the origin in  $\mathbb R^N$, i.e. $  y(0)=0$, and such that the   curvatures of these curves are bounded by a positive constant $c$, 
 $$curv^\perp(y(s))\leq c.$$
 
 If $l\leq 2\pi c$, then,  according to this inequality, 
 
 {\it the minimum of the distance from the second end of $y(s)$ to the origin, that is $\|y(l)\|$, is achieved  by planar circular  arcs of   curvature $c$.}

More generally Let $h(y)$ be a piece-wise smooth  function on a Riemannian manifold $Y$, let $\tau_0\in T_{y_0}(Y)$ be a unit tangent vector
and let $\sigma\geq 0$ be  a positive 
  Borel measure on the segment  $[0,l]$, e.g. $\sigma=c(s)ds$, $s\in [0,l]$, for a bounded measurable  function,
c(s) or $\sigma=\sum_ic_i\delta(s_i)$ as in 1.1.I.

{\sl \color {blue!44!black} Find an  isometric, i.e length preserving, immersion $[0,l]\to Y$, written as  
$s\mapsto y(s)$, such that     
$$y( 0)=y_0, \mbox { }  y'(0)=\tau_0\mbox { and } curv^\perp(y(s)\leq c(s),$$
which {\it minimises $h(y(l))$.}}\vspace {1mm}

This  is an instance of {\it an  optimal control problem}\footnote {Think of piloting a  jet plane, where  acceleration must be limited by a couple of $G$ for your comfort.\url {https://en.wikipedia.org/wiki/Optimal_control}.} where
 solutions are often piecewise smooth rather than smooth
according to 
  Fel'dbaum's {\it n-interval theorem} from the optimal systems control theory.\footnote {see  Optimal Control Systems
{\url https://www.scribd.com/document/390018919/Optimal-Control-Systems-Feldbaum-pdf}
and 
\url{https://encyclopediaofmath.org/index.php?title=Pontryagin_maximum_principle}}

The variational principle   behind  this theorem suggests  an effortless proof  of the classical combinatorial antecedent of the bow theorem.



 {\bf 8.A.  Cauchy(1813)-Legendre(1794) Arm Lemma.}\footnote{See [Sab 2004] and references to  the contributions   by Legendre, Cauchy and  Steinitz.} 
   Let  $Q\subset \mathbb R^N$ be a polygonal curve  with vertices $q_1,q_2,...q_k\in \mathbb R^N$ and segments $s_i=[q_i,q_{i+1}]\subset \mathbb R^N$,
   where the external angles between consequitive edges are bounded by positive  numbers $0\leq c_i\leq \pi$, that is 
   $$\angle _{q_i}=\angle (s_{i-1},s_i)\geq \pi-c_i, \mbox { } i=2,...,k-1.$$
  Let
  $P\subset \mathbb R^2$ be a locally  convex polygonal curve  in the plane     with vertices $p_1,p_2,...p_k$, where
  the segments have the same lengths $l_i$ as those in $P$,
   $$\|p_i-p_{i+1}\|=l_i=\|q_i-q_{i+1}\|$$
   and the external angles at the vertices $p_i$ are equal $c_i$.
   $$\angle _{p_i}=\pi-c_i,\mbox { } i=2,...,k-1.$$
  
  If the curve $P$ is {\it convex}, that is if the union of this curve with the segment $[p_,p_]$ makes a {\it closed convex}
   curve, then
  then the end points $q_1$ and $q_k$  of $Q$ lies {\it further apart} than these of $P$,
 $$ dist (q_1,q_k)\geq dist (p_1,p_k),\leqno{ [dist_Q\geq dist_P]}$$
  where {\it the equality implies that $Q$ is congruent to} $P$.

  {\it  Proof.} Let  $Q=Q_{extr}\subset \mathbb R^N$ be a polygonal curve which {\it minimises} the distance  
  $dist (q_1,q_k)$ among all curves with 
  $$ length(s_i)=l_i \mbox {  and }\angle_{p_i} \geq \pi-c_i.$$

  {\bf1}. {\it Extremality$\implies$Rigidity}. \footnote {Compare with Connelly} 
  If the angle of an  extremal curve $Q$ at some vertex point lies strictly between $0$ and $\pi-c_i$, 
  $$0<\angle_{q_{i_0}}<\pi-c_{i_0},$$ 
  and if $dist (q_1,q_k)>0$ 
  then the   points $q_{i_0},q_1,q_k$ are {\it collinear}. 
 
  Otherwise, the spheres around $q_{i_0}$ of radius $R_0=\|q_{i_0}-q_k\|$ and around $q_{1}$ of radius $R_1=\|q_1-q_k\|$ 
  would be transversal, 
 $$S_{q_{i_0}}^{N-1}(R_0) \pitchfork S_{q_{1}}^{N-1}(R_1),$$ 
 and a small variation of this angle along with small  rotation of the part of $Q$ following $q_{i_0}$  around the edge $[q_{i_0},q_{i_0+1}]$ would decrease the distance between $q_1$ and $q_k$. 
 
 {\bf 2.} {\it Vertex cut off:} [$k-1\implies k$]\footnote {compare p. 28  in [Cauchy] and with section 3 in Zaremba} Let the  angle of a (not necessarily extremal) $Q$ at   some point be minimal possible, 
 $$\angle_{q_{i_0}}=\pi-c_{i_0},$$ 
  let $Q'$ be the curve obtained  by removing $p_{i_0}$ and by joining  $p_{i_0-1}$ and $p_{i_0+1}$  by an edge 
 and let $P'$ be a  similarly truncated  curve $P$.
 Then the curve $P'$  is convex and the (old as well as  new) external angles of $Q'$ remain bounded  by those of $P'$.

 If $dist (q_1,q_k)>0$, then  
  induction on $k$   and  Lemma 1 reduce the arm theorem to where either all vertices  are collinear 
  or $k=3$ and where the validity of the lemma is obvious.
 
  {\bf 3}.{\it Terminal Edge Cut off:} [$(j-1){l_k\over m}  \implies j{l_k\over m}  $].
 Let $m$  be a positive   integer, such that  $\delta={l_k\over m}>0$, where $l_k$ is the length of the terminal edge in  $P$,  is  smaller than the distance from $p_1$ to this edge,
 $$ \delta={l_k\over m}< dist [p_{k-1}, p_{k}], \mbox { }  l_k=\|p_k-p_{k-1}\|$$
   Let $P_j\subset P$, $j=0,1,...,m$, be obtained from $P$   by cutting away   the terminal  $l_k-(m-j)\delta$  segment of length $l_k-(m-j)\delta$ from the edge  $[p_{k-1}, p_{k}]\subset P$; thus 
    $P_m=P$ and $P_0$ is a polygon with $k-1$ vertices.  
   
   Let $Q_j\subset Q$ be the corresponding parts of $Q$ and observe that (validity of) the inequality
   $[dist_Q\geq dist_P]$ for $Q_j$ implies that the distance between the end vertices in  $Q_{j+1}$ doesn't vanish, 
   hence, by {\bf 1} and {\bf 2},  $Q_{j+1}$ satisfy $[dist_Q\geq dist_P]$ for $Q_j$ as well, and  
     since $P_0$ has  only $k-1$ vertices, the proof of  $[dist_Q\geq dist_P]$ for  $P=P_m$ is conclude by 
  induction in $j$ and $k$.
   
   Finally notice that 
     the  congruence of  $Q$ and  $P$ in the case where  $dist (q_1,q_k)=dist (q_1,q_k)$
   follows by  tracking  strictness of angle inequalities in {\bf 2} of  this a argument.
  
 { \bf 8.B.  Euclidean  $C^2$-Bow Corollary} (A.  Schur 1920, E. Schmidt 1925)
  Let $y(s)$, $0\leq s\leq l$, be a smooth  curve in $\mathbb R^N$ parametrized by the arc length
  and let
    $$curv^\perp(y(s))=\|y'(s)\|\leq c(s), \mbox { } c(s)\geq 0.$$
Let $a(s)=a_c(s)$ be a  planar curve  with 
     $$curv^\perp(a(s))=c(s).$$
    (Such a curve, 
    which is locally convex, is unique up to congruence.) 
    
  If $a(s)$  is convex, that is if the union of the image $a[0,l]\subset \mathbb R^N$ and the straight segment 
    $[a(0), a(l)]\subset \mathbb R^2$ constitutes a (planar) closed convex curve (of length $l+\|a(l)-a(0)\|$),
then $y(s)$ satisfies the {\color {blue}$[dist_Y\geq dist_A]$} inequality:
$$dist(y(0),y(l))\geq dist(a(0), a(l)).$$

{\it In fact}, this follows from the above  $[dist_Q\geq dist_P]$ by approximation (see section 1.1.I)  of smooth curves by  the  polygonal ones.

 {\it Remarks.}  (a) There is  a   calculus proof  of this inequality, commonly   attributed  to E.Schmidt [Hopf 1946],
 which we reproduce in the next section.

 (b) The approximation argument works both ways:   "$C^2$-bow" implies   "polygonal arm". 
 
 (c) The approximation argument   also delivers the bow inequality for arbitrary convex curves $a(s)$, a direct proof of
  which is presented in [Sull 2007]
 \vspace {1mm} 
 
  {\bf 8.C. \it \textbf {Spherical and Hyperbolic Arms and  Bows.}} (i) The above proofs of arm and bow inequalities  extends verbatim to curves   
 in Riemannian manifolds $Y$ with constant sectional curvatures, spheres and hyperbolic spaces,\footnote {According to [Ni 2023] a proof of the hyperbolic bow inequality is presented in a   1985 preprint by C. L. Epstein, which I was unable to locate on the web. }where it also applies to    
 minimization of distance functions $h(y)$ for curves  with arbitrary initial conditions $(y(0),y'(0)))\in T)(Y)$
 
Also the spherical  bow inequality reduces to the Euclidean  one [Connely 1982, p. 31] as follows.
  
  Given positive numbers $r_i$, $i=1,...,k$, let us  denote the  corresponding "radial lift" of points $q_i\in S^{N-1}= S^{N-1}(1)$ to $\mathbb R^N\supset S^{N-1}(1)$  by
  $$q_i\mapsto \tilde q_i=r_i\cdot q_i\in \mathbb R^N,$$
  write 
  the distances  $\tilde d$ between  $\tilde q_i$ as  $\tilde d= \Phi(d)=\Phi_{\{r_i\}}(d)$, i.e. 
    $$dist_{\mathbb R^N}(\tilde q_i,\tilde q_j)=\Phi(dist_{S^{N-1}}( q_i,q_j)),$$ 
 similarly  express  the   {\it angles $\tilde \alpha$ between the  segments} $[\tilde q_{i-1},\tilde q_{i}],[\tilde q_{i},\tilde q_{i+1}]\subset \mathbb R^N $  as the function  in  the angles $\alpha$  between the  corresponding {\it geodesic segments}  in the sphere  $S^{N-1}$
  $$\angle_{\tilde q_i}=\Phi_\angle (\angle_{q_i}),$$
  and observe that  the functions $\Phi(d)$ and $\Phi_\angle(\alpha)$,  (determined by $\{r_i\}$) are {\it monotone increasing.}
  
If points $p_i\in S^2\subset S^{N -1}$ (consequently joined by geodesic segments)  form a {\it convex}  polygonal curve $P\subset S^{N-1}$,  
 then there exist numbers $r_i$, such that the polygonal curve   $\tilde P \subset \mathbb R^N$ with vertices  
$\tilde p_i=r_i\cdot p_i\in  \mathbb R^N$ is  {\it planar convex}.

  Indeed, since $P$, being convex, is contained in a hemisphere, say $S^{N-1}_+\subset S^{N-1}$,  the {\it radial} projection  $\rho$ from $P$ to the hyperplane $T\subset \mathbb R^N $, which is  tangent to $S^{N-1}_+$ at the centre, sends $P$ to a planar convex 
 polygonal curve  $\tilde P=\rho (P)\subset T$ and the  required $r_i$ are taken from the relations 
 $\rho(p_i)=r_i \cdot p_i.$

   Then, due to monotonicity of the functions $\Phi$ and $ \Phi_\angle$,
 the Euclidean arm lemma 
   for $\tilde p_i =\rho(p_i), \tilde q_i r_i\cdot q_i \in\mathbb R^N $
   with these very $r_i$  yields  the spherical   lemma for $ p_i,q_i \in S^{N-1}$.


  



\subsection {Semi-Circle Lemma and Calculus Proof of the Bow Inequality}


Let us look closer at the minimization problem for $h(y(l))$ for smooth curves $y(s)$, $s\in [0,l]$ in $Y=\mathbb R^N$ 
(see the beginning of the pevious section),
 such that     
$$y( 0)=y_0 \mbox { and }  y'(0)=\tau_0\mbox { and } curv^\perp(y(s))\leq c(s).$$
and where $h(y)$ is a linear function on $\mathbb R^N$.

Since  $Y=\mathbb R^N$, this can be  reformulated in terms of the derivative $z(s)=y'(s)\in S^{N-1}(1)$ as minimization of
problem for the integral $$\int_0^l  \langle z(s), grad(h)\rangle ds$$
under constraint $$\|z'(s)\|\leq c(s).$$

Let $h(y)$ be a {\it unit}  linear function, i.e. $\|{\sf grad} (h)\|=1$ and let 
    $y(s)$ start at the origin, $y(0)=0$ and satisfy 
     $$curv^\perp(y(s))\leq 1.$$

 Represent ${\sf grad}(h)$ by a point in the unit sphere, say $g\in S^{N-1}$, 
  and let $z(0)\in S^{N-1}$ represent the initial 
  derivative $z(0)=y'(0)$. 
  
{\bf  8.1.A.  Semi-Circle Lemma}.   The derivative $z(s)$  of the extremal $y(s)$, which minimises $h\circ y(l)$,
  follows the shortest geodesic arc in $S^{N-1}$ from $z(0)$ to $-g$ with constant speed $c$ for 
 $$s\leq   c^{-1} dist_{S^{N-1}}(z(0),-g)$$ and $z(s)$ is constant equal $-g$ for 
  $$s\geq   c^{-1} dist_{S^{N-1}}(z(0),-g)$$ 
 and  the extremal  curves $y(s)$ are planar circular arcs,  of length $\leq \pi /c$, which may be followed by straight segments for large $l$.

{\bf  8.1.B.  Semi-Circle Example/Corollary.} Let  the ends of a smooth immersed  curve $y(s) $, $0\leq s\leq l$, such that $curv^\perp (y(s))\leq 1$,  lie in two parallel hyperplanes $H_0,H_l\subset \mathbb R^N$, which are tangent to $y(s)$  at the corresponding ends $y(0)$ and $y(l)$. 

If   no tangent to the curve $y(s)$. at an  interior point $s\in (0.l)$ is parallel to $H_0$
and if there exists  a point $s_\perp\in [0,l] $,  where the tangent to the curve  is normal $H_0$, then
$$dist(H_0,H_l)\geq 2,$$
where the equality $dist(H_0,H_l)\geq 2$ implies that the   $[s_1,s_2]$-part of the curve for  an interval 
$[s_1,s_2]\subset  [0,l]$ is congruent to a planar semicircle of unit radius.
\vspace {1mm}

 Lemma 8.1.  If $curv^\perp(y(s)\leq c(s)$, where  $c(s)\geq 0$ is a non-constant  (bounded measurable), then arguing as in 8.1.A we obtain the folloing.
 
 {\bf  8.1.B. Minimal displacement Lemma}.   The  derivative $z_{ext}(s)$ of the  extremal $y_{ext}(s)$, which minimises $ h\circ y(l)$ follows the shortest geodesic segments from $z(0)$ to $- g$
    with (now  variable) speed  $c(s)$ for $$\int_0^s c(s)ds\leq dist_{S^{N-1}}(z(0), -g),$$  
    where  (this part  is irrelevant  for the  proof of the bow  inequality below) continues  with constant  $z(s)=-g$ for $s$ for
    $$\int_0^s c(s)ds\geq dist_{S^{N-1}}(z(0), -g).$$ 

It follows that

{ $\bf[a_{ext}]$}{ \sl if $y'(0) =-grad (h)$ and 
$$\int_0^l c(s)ds \leq \pi=dist_{S^{N-1}}(g, -g),$$ 
then the  extremal curves  $y_{ext}(s)=\int_0^s z_{ext}(s)ds$ are a planar\footnote{"Planar means contained in a 2-plane in $\mathbb RF^N.$} locally convex   curves, call them $a(s)=a_c(s)$, such that  
     $$curv^\perp(a(s))=c(s),$$}
  where, as we know, all these $a$-curves  are mutually congruent.\vspace{1mm}
    
\vspace {1mm}
  
  Now we are ready to prove   "bow inequality", that is,  recall,  the bound $$dist(y(0),y(l))\geq dist(a(0), a(l))$$
  for smooth curves $y(s)$, $s\in[0,l]$  in $\mathbb R^N$, such that 
 $$curv^\perp(y(s))\|y'(s)\|\leq c(s)$$ 
  and such that the corresponding (locally convex) curves $a(s)$ with $curv^\perp(a(s))=c(s)$
  are convex.
  
  This is done by dividing both curves $a(s) $ and  $y(s)$ into halves  by a point $s_\circ\in [0,l]$, such that 
 
 (1) The total  curvatures of the halves  of both curves are bounded by $\pi$, 
$$\int_0^{s_\circ} c(s)ds,  \int^l_{s_\circ} c(s)ds\leq \pi. $$
We apply  $\bf[a_{ext}]$ to these halves and to the linear function 
$$h(y)=h_\circ(y)=\langle  y, y' (s_\circ)\rangle,\mbox { } y\in \mathbb R^N.$$ 
which is characterised by  the equality  
$grad (h)(y(s_\circ))=y'(s_\circ)$   
   and obtain two inequalities
    $$ h_\circ y(s_\circ)-h_\circ y(0)\geq  h_\circ a(s_\circ)-h_\circ a(0)$$
    and
    $$ h_\circ y(l)-h_\circ y(s_\circ)\geq  h_\circ a(l)-h_\circ a(s_\circ),$$
  where the abbreviation $h_\circ y(s)$ stands for $h_\circ (y(s))$, etc.  and where we use same notation $h_\circ$ 
for the planar linear function $x\mapsto \langle  x, a'(s_\circ)\rangle$.
    
    These two add up to
    $$ h_\circ y(l)-h_\circ y(0)\geq  h_\circ a(l)-h_\circ a(0),
$$
where the first term is bounded by 
 $$ h_\circ y(l)-h_\circ y(0)\leq  dist(y(0),y(l)),$$
since $\|grad(h_\circ)\|\leq 1$. Therefore,
 {\color {blue} $$h_\circ a(l)-h_\circ\circ a(0)\leq dist(y(0),y(l)).$$}

 (2) Finally, let 
 $s_\circ\in [0,l]$ be a point,  where the tangent  to the (image of the)  curve $a(s)$ at $s=s_\circ$   is {\it parallel to the straight  segment} between the ends of this curve, that is 
  $[a(0),a(1)]\subset \mathbb R^2$.\footnote {Such a point, e.g. the one which  maximises the distance between  pairs of  points  $a(s)\in \mathbb R^2$, $s\in [0,l]$,  and $x\in  [a(0),a(1)]$, exists on all smooth immersed curves.}
  Then, by convexity of $a(s)$,
 {\color {blue} $$h_\circ a(l)-h_\circ a(0)=  dist(a(0),a(l))$$}
 and  since, also by convexity of $a(s)$, the total curvatures of the both halves of  $a(s)$, hence of $y(s),$ are bounded by $\pi$, (1) applies. Then confronting  the above  "blue" relations yields
     the bow inequality     
      {\color {blue}  $$dist(y(0),y(l))\geq dist(a(0),a(l)).$$ } \vspace {1mm} 
     
    {\it On the Error Term in the Bow Inequality.}
    One can sharpen the bow inequality (at least the circular one  for $\int curv^\perp\leq \pi$) with a bound on the total curvature, where the extremal ones, i.e. with minimal $dist(y(0),y(l))$ are circular arcs extended by straight segments at one of the ends. } Also one can evaluate    non-planarity of curves by   non-additivity of angles, e.g. in  decompositions of polygonal curves into triangles.

 
\vspace {2mm}

\hspace {46mm}{\sc Examples.}

\section {Angular    Bow InequaIity and non-Angular Corollaries }


Let  $Y$  and    $A$ be Riemannian manifolds and $y(s)$   and $a(s)$,
  be smooth  curves of length $l$ in these manifolds   parametrized by the arc length   $s\in [ 0,l]$, where
 $a(s)$ is a {\it convex ark} (as it is explained below)  in $A$ and where 

  the normal   curvature of $y(s)$ is everywhere bounded by that of $a(s)$: 
 $$curv^\perp(y(s))\leq curv^\perp(a(s)\mbox {  for all } s\in [0,l].$$

 {\bf 9.A.  Angular Bow Inequality for Negative Curvature.}

  
  Let both
  manifolds be complete simply connected with nonpositive sectional curvatures, 
let  $dim(A)=2$ and 
$$sect.curv(Y)\leq \inf_{a\in A} sect.curv(A,a).$$

Let the  ark $a(s)$  be {\it locally convex and  "acute":} the angles $\alpha _1$ and $\alpha_2$  between geodesic  segments $[a(s_1),a(s_2)]\subset A$ and the curve 
 $a(s)$ at both  ends are $\leq \pi/2$, where the  curve and the segments are oriented in same direction , where short segments are acute. (Arks in the unit circle  are "acute" for $ length\leq \pi$.)

Then  {\it angles $\beta_0$ and $\beta_l$  between  geodesic  segment $[y(0),y(l)]\subset Y$ and the curve $y(s)$ are bounded} by 
the corresponding angles between the  segment $[a(0),a(l)]\subset A$ and $a(s)$,
{\color {blue} $$\beta_0\leq \alpha_0\mbox  { and }\beta_l\leq \alpha_l.\leqno [\beta\leq \alpha]$$}



{\it Proof}. Let us enumerate the relevant properties of  
  the distance functions  $d=d(a_1,a_2)=dist_A (a_1, a_2)$ on  $A$ and $D=dist_Y (y_1), y(2))$ on$Y$ restricted to our curves. 
 
(1) Since both  manifolds  
complete simply connected with nonpositive sectional curvatures, 
the functions $d$ and $D$ are smooth away the diagonals.
   
  What is relevant for our purpose, is that they are

 {\it smooth  at the 
 pairs $(a(s_1),a(s_2))\in A\times A$ and 
the  pairs $(y(s_1),y(s_2))\in Y\times Y$ respectively for $s_2>s_1$.}

(2) Let
$$\alpha_{s_1} [\overrightarrow {s_1s_2}]=\angle( a'(s_1), -grad_{s_1} (d)),$$
where $ a'(s)={da(s)\over ds} \in T_{a(s)}(A)$ stands for  the derivative of $a$,  be the angle between the $a$-curve  at the point $a(s_1)$ and the directed  geodesic segment \footnote{Smoothnss of $d$ imy the uniqueness of the minimal segment.}
 $\overrightarrow{[a(s_1),a(s_2)]}\subset A$
and let 
 $$\alpha_{s_2}[\overrightarrow {s_1s_2}] =\angle (a'(s_2), grad_{s_2}(d)) $$
  be the angle between this curve  with same  segment   at the second end $a(s_2)$.

 Similarly, define such angles for the $y$-curve and denote them 
 $$\beta_{s_1} [\overrightarrow {s_1s_2}]\mbox {  and }\beta_{s_2}[\overrightarrow {s_1s_2}].$$
 
 Thus our objective is the inequalities 
 $$\beta_{s_i} [\overrightarrow {s_1s_2}]\leq \alpha_{s_i} [\overrightarrow {s_1s_2}]\mbox { for } i=1,2. $$

 Recall that the inequality {\sf [$sect.curv(Y)\leq 0$] and [smoothness of the distance function $D(y_1,y_2$]} imply that 
 \vspace {1mm}

({\Large$\ast$}) the $R$-spheres in $Y$ are {\it  greater than 
 $R$-spheres in  $\mathbb R^N,n=dim(Y)$}; in fact the exponential map $\exp_{y_1}:T_{y_1}(Y)\to Y$
 is {\it distance increasing}. \vspace {1mm}

This,  or rather the corresponding contraction property of the differential of   the (rescaled) inverse exponential map, 
can be expressed in terms of  the distance function as follows.

Let $$y_2\overset  {\chi}\mapsto grad_{y_2} (D(y_1,y_2))(y_1)\in T_{y_1}(Y), \mbox { } y_1, y_2\in Y$$
and let 
 $$\chi'=\chi'_{y_1,y_2}:  T_{y_2}(Y)\to T_{y_1}(Y)$$
 be the differential of ${\chi}$

 Then ({\Large$\ast$}) is equivalent to the bound of the norm $\chi'$ by the inverse distance function:
 $$ \|\chi'_{y_1,y_2}\|\leq D^{-1}(y_1,y_2). \leqno{(\mbox {\Large$\ast'$})}$$

 {\it Flat  Example}. If $Y=\mathbb R^N$, then 
 $\chi'(t)=t \cdot {1\over D(y_1,y_2)}$, $ t\in  \mathbb R^N=T_{y_1}(\mathbb R^N)= T_{y_2}(\mathbb R^N)$.

 The conditions [$sect.curv(Y)\leq 0$]\&[smoothness of $D$] imply the following  stronger geometric property than ({\Large$\ast$}):

({\Large$\ast \ast$}) the   $R$-spheres $S(R) \subset Y$ are "more convex"  than $R$-circles in $\mathbb R^2$:  the curvatures
 of these spheres at all unit tangent vectors $\tau\in T(S(R))$ are $\geq 1/R$,
 where the implication $(\ast\ast)\implies (\ast)$
is seen with  the metric definition of curvature (see sect1.1.E).
 
The above was meant to motivate  the following relation between the norms of the 
 the (linear) maps  $\chi'$ in $X$ and in $A$ which follows from the inequality $$sect.curv(Y)\leq \inf_{a\in A} sect.curv(A,a).$$

(3) The norms of the  maps  $\chi'$ in $X$ are   bounded by these in $A$ as follows
 $$\|\chi'_{y(s_1),y(s_2)}\|\leq \|\chi'_{a(s_1)a(s_2)}\|\mbox {  for $D(y(s_1), y(s_2)\geq d(a(s_1), a(s_2)$,} $$
 where this inequality is strict, $\|\chi'_{y(s_1),y(s_2)}\|< \|\chi'_{a(s_1)(s_2)}\|$, for $D(y(s_1), y(s_2)> d(a(s_1), a(s_2)$.

{\it Corollary.} The derivatives of the above defined angles. $\alpha$ and  $\beta$ satisfy
$$ \alpha'_{s_1} [\overrightarrow {s_1s_2}]= \chi'
 \cdot \sin\alpha_{s_2} [\overrightarrow {s_1s_2}]$$
and 
{\color {blue}$$\| \beta_{s_1} [\overrightarrow {s_1s_2}]\|\leq \chi'(s_1,s_2)
 \cdot \sin\beta_{s_2} [\overrightarrow {s_1s_2}]\mbox { for $D(y(s_1,y(s_2)\geq d(a(s_1,a(s_2)$},$$}
where $\chi'$ in both cases is for $A$, i.e.  $\chi'=\chi'_{a(s_1)a(s_2)}$,  where 

{\color {blue}this $\chi'$ is {\it strictly monotone
decreasing in the distance  $d(a(s_1)a(s_2))$}} 
 
 which makes  
 
 {\color {blue}the $\beta$-inequality strict for 
    for $D(y(s_1,y(s_2)> d(a(s_1,a(s_2)$.}
 
 {\it Remark}.This "$\sin$" plays no special role in our proof of the inequality $[\beta\leq \alpha]$
 and it can be replaced by another positive monotone strictly increasing function on the segment $[0,\pi/2]$.

 (4) The inequality $sect.curv(Y)\leq \inf_{a\in A} sect.curv(A,a)$ impies that 
the curvatures of the spheres in $X$ are greater or equal than  of the  spheres of same radii in $A$, nd $a(s_2)$.
 
 I fact we shall need this bound only for 
 the  spheres  $S_{y(s_1)}\hspace {-1.6mm}\uparrow^{y(s_2)} $  centred at  $y(s_1)\in Y$  which contains $y(s_2),$ 
 and will use  it
 at the points ${y(s_2)}$, where "greater or equal" applies to the curvatures of these spheres at all their unit tangent vectors,
 $\tau \in T_{s_2}(S_{y(s_1)}\hspace {-1.6mm}\uparrow^{y(s_2)})y(s_2)$:
$$curv_\tau^\perp(S_{y(s_1)}\hspace {-1.6mm}\uparrow^{y(s_2)})y(s_2)\geq curv^\perp(S_{a(s_1)}\hspace {-1.6mm}\uparrow^{a(s_2)})a(s_2). $$
(The sphere $S_{a(s_1)}\hspace {-1.6mm}\uparrow^{a(s_2)}$ is 1-dimensional with a $\pm$single
tangent vector at $a(s_2)$.)

(5) The {\it derivatives} of the distance functions $ D(s_2)=D_{s_1}(s_2)=D(y (s_1),y(s_2))$ and $d(s_2)=d(a (s_1),a(s_2))$  in the $s_2$ variable are {\it monotone decreasing} in the angle   $\beta_{s_2}=\beta_{s_2} [\overrightarrow {s_1s_2}]$. \vspace {1mm}
 In fact, 
 $$D'(s_2)=\cos \beta_{s_2} \mbox { and } d'(s_2)=\cos \alpha_{s_2}$$ 
 and by symmetry of the distance function,
 $$D'(s_1)=-\cos \beta_{s_2} \mbox { and } d'(s_2)=-\cos \alpha_{s_2}$$ 
for  $ D(s_1)=-D_{s_2}(s_1)=D(y (s_1),y(s_2)$ and $d(s_1)=d_{s_2}(s_1)=d(a (s_1),a(s_2)$.

(6) Finally we turn to the contribution of the $curv^\perp$-curvatures of the curves  $a(s)$ and $y(s)$ 
to the derivatives of the angles $\alpha(s_2)=\alpha_{s_2} [\overrightarrow {s_1s_2}]$ 
and $\beta(s_2)=\beta_{s_2} [\overrightarrow {s_1s_2}]$ and thus to the second derivatives of the distance functions 
$d(a(s_1),a(s_2))$  and $D(y(s_1),y(s_2))$ in the second variable.

Let $c_a(s_2)= curv^\perp(a(s_2))$ be the normal curvature of the curve $a(s)$ in $A$  at the point $a(s_2)\in A$
and $c_y(s_2)=curv^\perp(y(s_2))$ be this for $y(s)$ in $Y$.

Let $c_d(s_2)$ be the normal curvature of the circle $S^1={a(s_1)}\hspace {-1.6mm}\uparrow^{a(s_2)})$  (with centre $a(s_1)$ and  radius 
$dist_A(a(s_1), a(s_2))$))
at the point $a(s_2)$
and  let $c_D(s_2)$ be the normal curvature of the sphere $S^{N-1}=S_{y(s_1)}\hspace {-1.6mm}\uparrow^{y(s_2)})$  (with centre $y(s_1)$ and  radius 
$dist_Y(y(s_1), y(s_2))$)
 at the unit tangent vector $\tau\in T_{y(s_2}(S^{N-1})$ that is the normalized normal projection of the unit tangent vector to the curve $y(s)$ at $s=s_2$, (the derivative $y'(s_2)\in T_{y(s_2)}(Y)$)  to the tangent space 
 $T_{y(s_2)} S^{N-1}\subset  T_{y(s_2)}(Y).$

Then $$\alpha'(s_2)=c_a(s_2)-c_d(s_2), \mbox { and  }\|\beta'(s_2)\|\leq c_y(s_2)-c_D(s_2),  $$

Therefore, our assumptions $curv^\perp(y(s_2)\leq curv^\perp(y(s_2)$ and

$sect.curv(Y)\leq \inf_{a\in A} sect.curv(A,a)$ imply that 

{\it the $s_2$-derivatives of the 
$\beta$-angles} of the curve  $y(s)$ in the manifold $Y$ {\it are bounded by these  of the $\alpha$-angles} of  (the convex "acute"  arc) $a(s)$ in the surface  $A$,
{\color {blue}$$\|\beta'(s_2)[\overrightarrow {s_1s_2}] \| \leq \alpha'(s_2)[\overrightarrow {s_1s_2}]\leqno {[\beta'\leq \alpha']}$$}
for all pairs of points $0\leq s_1\leq s_2\leq l$, such that {\color {blue}$ D=d$}, i.e. where 
{\color {blue}$$dist_Y (y(s_1), y(s_2)=dist_Y (y(s_1), y(s_2)).$$}

The proof of the 
{\color {blue} $ [\beta\leq \alpha]$} inequality is concluded with the following.

{\bf 9.B. \it \textbf {Elementary Calculus  Lemma}}. Let $d(s_1,s_2)$ and $D(s_1,s_2)$, $0\leq s_1, s_2\leq l$, be  positive symmetric functions,  which are strictly positive and   smooth for $s_1\neq s_2$ and which 
are approximately equal to $|s_1-s_2|$ at  the diagonal $\{s_1=s_2\}$,
$$d(s_1,s_1+\varepsilon)=|\varepsilon|+O(\varepsilon^3)\mbox { and also } D(s_1,s_1+\varepsilon)=|\varepsilon|+O(\varepsilon^3). $$

 If these functions  satisfy the above (1)-(6) relations, where the $\alpha$- and $\beta$-angles are defined as
the derivatives of these functions,
then these angles  satisfy the $ [\beta\leq \alpha]$ inequality, where either of the equalities   
$\beta_0=\alpha _0$ or $\beta_l=\alpha _l$ implies that the two functions are equal,
$d(s_1,s_2)=D(s_1,s_2)$.
ely curved metrics 
 
\hspace {40mm} {\sc  Corollaries} \vspace {1mm}

Since  the angles $\alpha_s$   and $\beta_s$ are  represent the derivatives of the distance functions 
 $d(a(0),a(s)$ and $D'(y(0),y(s)$, the   $\beta \leq \alpha$ implies that 
$dist_Y(y(0),y(l) \geq dist_A(a(0),a(l)$ for convex "acute" arks .

More interestingly, the same holds  for some non-acute arks which can be divided two acute segments by a point 
$s_\circ$.

Thus we arrive at the following

{\bf  9.C. (De)composed  Bow Inequality}.    Let $A$ and $Y$ be complete simply connected manifolds  with non-positive  sectional curvatures 
 let
$$sect.curv(Y)\leq \inf_{a\in A} sect.curv(A,a)$$
and let 
 $$curv^\perp(y(s))\leq curv^\perp(a(s)\mbox {  for all } s\in [0,l].$$ 

Let the curve  $a(s)$ be locally convex and let  $s_\circ\in [0,l]$ be a point, such that   
the two  parts  $a[0,s_\circ]$ and  $a[s_\circ,l]$ of the curve  are "acute":

 if either  $s_1, s_2\leq s_\circ$ or
$s_1, s_2\geq s_\circ$, then 
the  angles between   the  geodesic segment $[a(s_1), a(s_2)]$ at its ends with 
the curve $a(s)$ are acute.

Then 
$$dist_Y(a(0),a(l))\geq dist_A(a(0), a(l)),\leqno {\mbox {\color {blue}$[dist_Y\geq dist_A]$}}$$
provided  one of the following two conditions is satisfied.
 $$dist_A(a(0),a(s_\circ))=dist_A(a(s_\circ), a(l)),\leqno {\bf Condition \hspace {1mm} 1 }$$
$$dist_A(a(0),a(s_\circ))\geq l-s_\circ.\leqno {\bf Condition\hspace {1mm} 2 }$$

{\it In fact}, the angular inequality bow inequality "$ [\beta\leq \alpha]$" shows that the angle $\angle_y(s_\circ)$  between the geodesic segments   
$[y(s_\circ),y(0)]$ and $[y(s_\circ),y(l)]$ in $Y$ is greater or equal than the the angle  $\angle_a(s_\circ)$ between $[a(s_\circ),a(0)]$ and $[a(s_\circ),a(l)]$ in $A$ and the proof follows by the standard  comparison theorems  for
 geodesic triangles.\vspace{1mm}

 {\bf 9.D. \it\textbf { Subcorollary: Riemannian Circular Bow Inequality.}}
    Let $A$ and $Y$ be complete simply connected manifolds, where $A$ has constant  sectional curvature $\kappa$ and  
     $sect.curv(Y)\leq \kappa$, where for $\kappa>0$ we require that all geodesic segments of length $<\pi/\sqrt \kappa$ are distance minimising.

Let the curve  $a(s)$ be (planar) convex with constant  curvature $c\geq 0$ and 
 $$curv^\perp(y(s))\leq c.$$ 
Then
$$dist_Y(a(0),a(l))\geq dist_A(a(0), a(l)),\leqno {\mbox {\color {blue}$[dist_Y\geq dist_A]$}}$$

{\it Proof}. If $A$ and $a(s)$ have constant curvatures (each of its own kind), then the  above two conditions are satisfied  and the case of $\kappa\leq 0$ follows. \vspace{1mm}

 {\bf 9.E. Non-Simply Connected  Generalization and $\kappa>0$}.  The above argument applies to arc-length parametrized  curves  $[0,l]\ni s\mapsto y(s)\in Y$  in manifolds $Y$, such that 
 
$ \bullet$ every subsegment $[0,l] \supset [s_1,s_2]\to Y$ admits a  homotopy by curves with fixed ends and 
with  length$\leq s_2-s_1$  to a  geodesic segment $[y(s_1),y(s_2)]$ in $Y$ and  there is only one segment achievable  by such a homotopy, where   the length  of this segment takes the role 
 of the distance $D$ between $s_1$ and $s_2$;

 $ \bullet$ the exponential map $\exp_{y(s)}$ is  defined   on the balls $B(R)\subset T_{y(s)}(Y)$ of radii
 $R\leq \max(s, l-s)$ for all $ s\in [0,l]$, where this map is an immersion which is    strictly locally convex on the $R$-spheres.  

 For instance, these conditions are satisfied by complete manifolds  $Y$ with $sect.curv (Y)\leq \varepsilon^2 $
 and curves $y(s)$ of length $l\leq \pi/2\varepsilon$ and the $\beta$-angles of these curves are bounded by those of 
 convex  acute arcs  $a(s)$ in the 2-sphere $S^2(1/\varepsilon ^2)$, if $curv^\perp (a(s))\geq curv^\perp (y(s)).$
 
 Now, if the distance between the ensds of a curve $C\in Y$ is understood as the length of the geodecsic obtained by shortening homotopy of $C$ the above argument delivers the so modified  Riemannian Circular Bow Inequality remain valid for all complete manifolds $Y$ with 
 $sect.curv (Y)\leq \kappa$. and all $-\infty <\kappa< \infty$

 \vspace {1mm}

 Alternatively 
if one insists on manifolds being simply connected and on keeping true metric  distances, then
the condition $sect.curv(X)\leq \kappa$ for $\kappa>0$  must be augmented by $inj.rad (X)\geq \pi/\sqrt \kappa$. 

Then  the problem reduces to the case of $\kappa=0$ by taking the 
$1/\sqrt \kappa$-cone $Cone(X)$, which albeit is singular,  has $sect.curv\leq 0$ in the sense of  Alexandrov's $CAT(\kappa)$-theory, where the   above  analytic proof can be carried over with
minor modifications.

\vspace {1mm}

{\it Remark.} The circular bow inequality in the $(N-1)$-sphere $S^{N-1}\subset\mathbb R^N$ effortlesly follows from that for $\mathbb R^N$, 
since curves of constant (geodesic) curvature in the sphere are {\it planar}, i.e. contained in planes in $4\mathbb R^N$, 

{\it Exercise.} Show that (domains in) spheres are the only smooth submanifolds in  $\mathbb R^N$, where  
curves of constant geodesic)] curvature  are {\it planar}.\footnote{I didn't solve this exercise.}



\section { Hypersurfaces in   Balls and Spheres.}


 {\bf   Hypersurfaces  with Unit Curvature in  $B(2)$ and in $B(2+\delta)$. }  Let 
 the image of an  immersion  $f: X\hookrightarrow\mathbb R^{n+1}$, $n=dim(X)$, such that  $curv^\perp(X)\leq 1$,  is contained in the ball $B^{n+1}(2)$.

{\bf10.A. \it \textbf{ {\small \bf $\Circle$\hspace {-0.35mm}$\Circle$}-Extremality.}} If $n=1$ and the {\it degree of the Gauss map $S^1=X\to S^1(1)\subset \mathbb R^2$ equals zero}, then (this was stated in 1.C)  the image $f(X)\subset B^2(1)$ equals the union
 of {\it two unit  circles  } which tangentially  meet at the center of the  disk $B^2(2)$. 

{\bf10.B. \it \textbf{Star Convex Rigidity.}} If $n\geq 2$, then either $f(X)$ is star convex and  the radial projection $X\to S^n(2)$ is a  diffeomorphism, or $f(X)$ is equal to a unit $n$-sphere $S^n_{y_o}(1)$,  where the  centre of this sphere  is positioned  half way from the boundary of the ball $B^{n+1}(2)$,
i.e. $\|y_o\|=1$.

{\bf10.C. \it\textbf{Star Convex Stability Switch.}} There exists   $\delta > 0.01$  {\color {magenta} (probably $\delta >0.2$}),  such that if $n\geq 2$ and the image $f(X)$ is contained in the ball
 $B^{n+1}(2\varepsilon))$ 
then $f(X)\subset B^{n+1}(2+\delta)$ is star convex with respect to some point in $B^{n+1}(2+\varepsilon)$.

{\it Proof.} If $f(X)$ is not star convex with respect to the centre of the ball $B^{n+1}(2)$ then
 some radial ray is tangent to $ f(X)$ at some point $y_0=f(x_0)\in f(X)$ and  the semi-circle lemma 8.1.A  implies that $y_0$ is equal to the centre of $B^{n+1}(2)$ and 
 that $f(x)$ equals a unit sphere passing through   $y_0$.

This proves (b) while the stability of the bow proof  argument (see section 8.1) yields  an approximate unit sphere in  $B^{n+1}(2+\varepsilon)$ and (c) follows as well.

{\it Remark/Example}.  The boundary $X_{+1}$ of the $\rho$-neighbourhood  for $\rho=1$ of a  circular ark $S$ with radius $2$  has curvature bounded by $1$. If such an 
 $S$ is  slightly shorter than half circle, then, because of "shorter",  $X_{+1}$  can be fit to the  ball of radius $3-\epsilon$  and  $X_{+1}$ and it is 
 non-star convex because of "slghtly".

{\it  \color {magenta}Question} Do $\delta$ and $\epsilon$ ever meet or there is a definite gap between their possible values?

  {\bf10.D. Hypersurfaces in  $S^{n+1}$.}
  It is {\color {magenta} unknown} if there are closed connected $n$-manifolds $X$ non- diffeomorphic  to spheres  $S^n$, which admit immersions to the unit balls $B^{n+1}$ with 
  $curv^\perp(X\hookrightarrow B^{n+1}(1))<3$. 
  
  {\color {magenta} Nor do we know}  if there exist codimension two  immersions of these $X$  with 
   $curv^\perp(X\hookrightarrow B^{n+2}(1))<\sqrt 2$. 
    
      {\color {magenta} Conjecturally}, the only non-spherical immersions with {\it "critical curvatures"}, that are  {\color {magenta} 3 for codimension 1} and {\color {magenta} $\sqrt 2$ for codimension 2}, are the standard embeddings\footnote { If either $n_1$ or  $n_2$ is equal to 1, one may have covers of these embedded $X$.} of $S^{n_1}\times S^{n_2}$ to $B^{n+1}(1))$ and to $B^{n+2}(1))$, where the latter are Clifford's and the former are encircling of round spheres. 
   
But the (first) critical curvature is unquestionably is {\it equal to one}  for hypersurfaces {\it in the unit spheres}. Here there are  non-spherical submanifold with unit curvatures, namely 
Cliffords product of spheres 
  $$S^{n_1}(1/\sqrt 2)\times S^{n_2}(1/\sqrt 2)\subset S^{n_1+n_2+1}(1),$$
  which have curvature $curv^\perp=1$ and 
there is nothing of the kind for smaller curvature:

\SnowflakeChevronBold {\it Closed connected orientable\footnote {This {\color {magenta}must be} redundant. Anyway this is relevant only for those even $n$,  where  $\mathbb RP^N  $ admits an immersion to $\mathbb R^{n+1}$.} immersed hypersurfaces $X\hookrightarrow S^{n+1}(1)$ with 
  $curv^\perp<1$ are diffeomorphic to the  $n$-sphere $S^n$}.(Compare [Ge 2021].)

  \vspace {1mm}

   {\it\color {blue} More generally,}
    let  $X$ be a smooth closed connected $n$-manifold, let
$Y$ be a  closed connected simply connected Riemannian $(n+1)$-manifold,     let $ f:X\to Y$ be a cooriented (two-sided) immersion and 
   $f^\perp_ {\pm t}:X^n\hookrightarrow Y(1)$, $t\geq 0$, be (the composition 
   of  $f$ with) the normal exponential map $\exp^\perp: X\times \mathbb R_{\pm}\to Y$ at $t\in \mathbb R_+$. 
  
  Define two  {\it one-sided focal radii}: $rad{_\pm}^\perp(X)$ that are  the suprema of $t>0$ for which the maps f$^\perp_ {\pm t}$ are immersions. 
  
  Notice that $rad^\perp(X)=\min(rad_+^\perp(X),rad_-^\perp(X)$.

   {\sl If $sect.curv(Y) \geq 1$ and 
   $$rad_+^\perp(X)+rad_-^\perp(X) > \pi/2,$$
   e.g. $rad^\perp(X\hookrightarrow Y)\geq \pi/4,$
    $X$ is diffeomorphic to the sphere $S^n$. }

   {\it Proof.}       Observe following [Ge 2021]  that
  if $rad{_+}^\perp(X)>t_1$ and 
  $rad{_-}^\perp(X)\geq t_2$, then    the  (kind of composed) map
  $$(f^\perp_ {+t_1})_{-(t_1+t_2)}:X\to Y$$ is defned and it  satsfies 
  $$(f^\perp_ {+t_1})_{-(t_1+t_2)}=f^\perp_ {-t_2}.$$
  
  It follows that if  $rad_+^\perp(X)+ rad_-^\perp(X)>\pi/2$, 
 then  
  the map $f^\perp_{-t_2}: X\to Y$ is an immersion 
   (by the definition of $rad^\perp$) and if 
   $$\mbox{$sect.curv(Y)\geq 1$ and $t_1+t_2>\pi/2$}$$
   this immersion is {\it locally strictly concave}.\footnote { A smooth immersion $f:X\to Y$ is locally strictly  convex/concave if the second fundamental form of $f(X)$ is positive/negative definite.   
 "Convex" is distinguished from"concave" 
    for families of immersions $f_t$: {\it convexity}   indicates  i{\it increase} 
    of 
    the induced metric in $X$ and concavity corresponds to  decrease of this metric.
   Thus the boundaries  $X$ of convex sets are convex for outward deformations of $X$ and concave for inward deformations. Similarly, we  attribute  convexity and concavity to families of non-smooth hypersurfaces.}

   In fact, if $sect.curv(Y)\geq 1$ and
   $f_\ast :X\to Y$ is an immersion, such that $rad_-^\perp(f_\ast(X))>\pi/2$,
   then the map $(f_\ast)_{-t} X\to Y$ is  a  locally strictly concave immersion in the range   
   $\pi/2<t< rad_-^\perp(X)$
   by the Hermann Weyl tube curvature  formula.\footnote {Here and in the contraction argument below we follow [Esch 1886] and [Gro 1990].}

   
  Then   Gromoll-Meyer's  {\it concave  contraction }
   adjusted to immersions\footnote {This argument needs $dim(X)=n\geq 2$, which can be assumed in the present case.}
 and followed by  smoothing 
   delivers   a regular homotopy of the (locally concave) immersion 
   $f_{t_2}=f^\perp_{-t_2}:X\to Y$, say   
  $\bar f_t:X^n\to Y$, where  $t_2\leq t <t_\bullet$,  for 
   some $t_\bullet>t_2$, where $\bar f_{t_2}=f_{t_2}$ and   such that the locally concave immersions  $\bar f_t$ become {\it  concave  embeddings}\footnote {These are boundaries of (geodesically  convex subsets in $Y$ with inward directed  normal fields on them.} for $t$ close to $t_\bullet $,  which eventually  converge to a constant map for $t\to t_\bullet$. 
   
   Since  "small" closed convex hypersurfaces in  (all  compact Riemannian manifolds) $Y$  are diffeomorphic to spheres, the proof     follows. \footnote {Ge assumes that  $f$ is an embedding and he  proves the existence of a diffeomorphism $X \cong_{diff}S^n$ for $n\neq 4$  by the $h$-cobordism theorem.}
   
   \vspace{1mm}
   {\it Remark.} The above applies to  $X$ which satisfy a {\it point-wise version} of the inequality $rad_+^\perp(X)+rad_-^\perp(X) > \pi/2$,that is 
$rad_+^\perp(X,x)+rad_-^\perp(X,x) > \pi/2$ for all $x\in X$,
$rad_\pm^\perp(X,x)$ is the maximal $r$,where the immersion condition on the exponential map is required only  at thepoint $x\in X.$

{\color {magenta}Problem.} Show that 
{\sl the only topologically non-spherical  smoothely  immersed closed  connected hypersurfaces  $X$ with $curv^\perp(X)\leq 1$  
in unit spheres  are Cliffords $S^{n_1}\times S^{n_2}\subset S^{n_1+n_2+1(1)}$.}
  
  More generally,  if a closed $n$-manifold, which is  non-diffeomorphic to $S^n$,  is  immersed  
to  a simply connected $(n+1)$-manifold $Y$ with  $sect.curv(Y) \geq 1,$  
such that $rad_+^\perp(X)+rad_-^\perp(X) \geq  \pi/2$, or  if 
$rad_+^\perp(X,x)+rad_-^\perp(X,x) \geq  \pi/2$ $\forall x\in X$,
then $Y=S^{n+1}(1)$, {\color {magenta} conjecturally} 
and $X$ is Clifford's product of spheres.\footnote {Consulting    [Ge 2021] [Luis Guijarro and Frederick Wilhelm, Focal Radius, Rigidity, and Lower Curvature Bounds, 2017],[Grisha Perelman. Proof of the soul conjecture of cheeger and gromoll. Journal of Differential Geometry, 
1994]   [D. Gromoll and K. Grove, A generalization of Berger’s rigidity theorem for positively curved
manifolds], is instructive .}

\vspace {1mm}


  {\it Remarks.} (a) Since  the sectional curvature  of the induced metric in  $X$ 
  is non-negative by Gauss' formula, $\bullet_{\leq }$ the only non-spherical surfaces 
  in $(1)$ with $curv^\perp\leq 1$ are flat tori,  which are, by a simple argument, are
  coverings of Clifford's tori.

  \vspace {1mm}

{\it Exercise}\footnote {Compare with [Ge 2021] and section 3.7.3(F) in [Gro 2021]}. Let $W$ be a compact Riemannian $(n+1)$-manifold with two boundary components  $\partial_\pm(W)$ with the distance $dist(\partial_-,\partial_+)\geq  {\pi \over 2}+\varepsilon$, e.g. $V$ equals the 
$\rho$-neighbourhood, $\rho=\frac {\pi}{4} +\varepsilon/2$, of a hypersurface $X\subset S^n$ with $curv^\perp (X)<1$.


 Let $sect.curv (W)\geq 1 $ and let $W$ admit a locally isometric immersion to a complete ${n+1}$-manifold $Y$ with $sect.curv(Y)>0$.
 
 Show that there exists an embedded  $n$-sphere   $S^n\subset W$, which separates 
the two  boundary components of  $W$.

\section{Immersed Submanifolds in  Bands and
in Tubes}


  Let an $n$-dimensional manifold $X$ be immersed 
 to the $k$-{\it tube $B_{\mathbb R^k}^N(R)$ of radius} $R$,
 $$X \overset {f}\hookrightarrow B_{\mathbb R^k}^N(R)=B^N(R)\times \mathbb R^k\subset \mathbb R^{N+k},$$ 
  (where  $B^N(R)=B_0^N(R)\subset \mathbb R^N$ is the $R$-ball). 
 
  Let $p:X\to \mathbb R_{ax}^k=\{0\}\times \mathbb R^k$ be the projection of 
  $X\hookrightarrow B_{\mathbb R^k}^N(R)$ to the  central axes of the tube,  let 
 $$\mathcal K=\mathcal K(p)\subset T(X)\hookrightarrow T(B_{\mathbb R^k}^N(R))$$
 be the {\it kernel of the differential} $dp:T(X)\to T(B_{\mathbb R^k}^N(R)).$

Let
$$\Sigma=\Sigma(p) =\{x\in X\}_{rank(\mathcal K_x)> 0}\subset X$$ 
be the support of $\cal K$.\footnote {If $k<n$, then $\Sigma=X$ and   $rank (\mathcal K_x(p))=n-k$, for generic maps $p:X^n\to \mathbb R^k$ and generic points $x\in X $.  If $k\geq n$, then ether  $p$ is an immersion, i.e. $\Sigma=\emptyset$,  
or $dim(\Sigma(p))=2n-k-1$ for generic $p$ and 
 and $rank (\mathcal K_x(p))=1$ at generic $x\in \Sigma$.}

 Let the induced Riemannian metric in $X$ be {\it geodesically complete}, e.g. $X$ is compact without boundary, and let 
$$\gamma_\tau(l)\hookrightarrow X \hookrightarrow  B_{\mathbb R^k}^N(R),\mbox { } \tau\in \mathcal K_x$$ 
be the geodesic segment of length $l$  issuing from $x\in \Sigma$ in the $\tau$-direction, where $\tau$ is a {\it non-zero vector } in  
the vector  (sub)space $\mathcal K_x\subset T_x(X), x\in \Sigma.$

If
$$curv^\perp(X\hookrightarrow B_{\mathbb R^k}^N(R)\subset  \mathbb R^{N+k})\leq 1/R,$$ 
then the half circle lemma  applied to the curves $\gamma_{\pm \tau}({1\over 2}\pi R)$ in the $R$-tube  $ B_{\mathbb R^k}^N(R)$ and to the hyperplane $H=H_{\perp \tau}\subset \mathbb R^{N+k}\supset B_{\mathbb R^k}^N(R)$, which contains $f(x)\in B_{\mathbb R^k}^N(R)$ and is   normal to $\tau$    imply the following.

  [${\pi\over 4}\ast {\pi\over 4}$] Either $\Sigma=\emptyset$, i.e. $p:X\to\mathbb R^k$ is an immersion, (in this case one may have $curv^\perp(X)<1/R$) the  curves  $f(\gamma_{\pm\tau}({1\over 2}\pi R))$ and $f(\gamma_{\pm\tau}({1\over 2}\pi R))$ are composed of  quarter's of  planar  
 circlers, both of which reach the boundary of the tube and where  the ( normal) curvature vectors of these curves are parallel to 
the  central axes  $\mathbb R_{ax}^k$ of the tube, i.e. the planes containing these curves are perpendicular to the $N$-(sub)space, which  contains the ball $B^N(R)$
 (which is normal to $\mathbb R_{ax}^k$).
Thus, $f(X)\subset B_{\mathbb R^k}^N(R)$ intersect the boundary of $B_{\mathbb R^k}^N(R)$ at at least two points.

 [${\pi\over 2}$]   If an immersion  $f$ is $C^2$-smooth, so is the $\pi R$-curve in the tube made of $f(\gamma_{\tau}({1\over 2}\pi R))$ and $f(\gamma_{- \tau}({1\over 2}\pi R))$. This necessarily makes this   curve  a planar half circle.

 In general, $C^{1,1}$-curves composed  of circular arc of same  curvature $1/R$ are  not always planar arks themselves.

   However  if   the  geodesic segments  $\gamma_\tau$ issuing from a point $x$ in  all  directions $\tau\in T_x(X)$ in an    $n$-dimensional   $C^{1,1}$-immersed $n$-manifolds $X\hookrightarrow \mathbb R^N$,      are  planar circular arcs of same  curvatures $c>0$,  an if $n\geq 2$, then
    $x\in  \Sigma \subset X^n$, serves as the centre of a  geodesic $R$-hemisphere $(S^m_+)_x\subset X$ of dimension 
$m=rank (\mathcal K)$, such that the map $f$ isometrically sends $(S^m_+)_x$ to an equatorial $m$-hemisphere in the $(N-1)$-sphere $S^{-1}_{p(x)}$, where  the boundary of this hemisphere is contained in the boundary of the tube $ B_{\mathbb R^k}^N(R)$. 
\vspace {1mm}

{\it Exercises.}  (a). Recall that the real projective spaces of dimension $n=2^l$,   admit no immersions to  $\mathbb R^k$ for $k\leq 2n-2$, and show that they admit  no immersions to the tubes $B_{\mathbb R^k}^N(R)$ with $curv^\perp(f)<1/R$.

  (b)  Let a closed connected $n$-manifold $X$, $n\geq 2$, be immersed  to a (cylindrical)  $(1,R)$-tube 
  $$X\overset {f}\to B^N_{\mathbb R^1}(R)\subset \mathbb R^{N+1}.$$

 (b$_1$) Show that {\it the only  critical points} $x\in X$ of the function $p:X\to \mathbb R=\mathbb R^1_{ax}$, i.e. where 
 $rank(\mathcal  K_x)=n$,  are the   maximum  $x_{max}$ and the   minimum  $x_{min}$ ponts of $p$, and that the $f$-images of  both of them in  the tube 
 are positioned on the axial line 
$\mathbb R_{ax}^1=\{0\}\times \mathbb R^1$, where they  serve as the centers of $n$-hemispheres $(S^n_+)_{max}(R)$ and $(S^n_+)_{min}(R)$, both  of radius $R$ and where both  are contained in the 
$f(X)\subset  B^N_{\mathbb R^1}(R)$ and where the spherical ($S^{n-1}(R)$) boundaries of them  are contained in the boundary of the tube.
 
 (b$_2$)  Show that the  $(n-1)$-hemispheres  $(S^{n-1}_+)_{f(x)}(R)$ in the tube tangent at their centres $y=f(x)$ to the (topologically $(n-1)$-spherical) fibers of the map $p$  for all non-critical points $x\in X$ {\it continuously} depend on $x$ .
 
 Observe that the resulting continuous map from the $(n-2)$-sphere bundle 
 $U$ over $X\setminus \{x_{max},x_{min}\}$ of unit vectors tangent to the fibres to $X$, 
 say $\Phi:U\to X$,
 sends $U$ to the intersection of   $X$ with the boundary of the tube.

(b$_3$) Show that the image of the immersion $f:X\to B^N_{\mathbb R^1}(R)$ equals the union of the two hemi-spherical cups  $(S^n_+)_{max}(R)$ and $(S^n_+)_{min}(R)$ and a region between them {\it contained in the boundary of the tube}.
 that is equal the $\Phi$-image of $X\setminus ((S^n_+)_{max}(R)\cup (S^n_+)_{min}(R)$
 
{\it Hint } Start with case $n=2$

 (i) Let  $n\geq 2$ and $N=n$and  show that $f$ is an embedding, the image of which is equal the  $+R$-encircling of a segment  in the  central line  
 $\mathbb R^1$ in  $B_{\mathbb R^1}(R)^N$, that is 
 a (convex) region between to half-$R$-spheres normal to this line, which is  equal in the present case to the boundary of the convex hull of $f(X)\subset B_{\mathbb R^1}(R)^N$. 
  (Unless  the half-$R$-spheres have a common boundary, this region is only {\it piecewise} $  C^2$.) 
   
 (ii)  Let  $n\geq 2$ and   $N>n$.  Show that  $f$ is an embedding into the  $+R$-encircling of a central segment  $[a,b]\subset
 \mathbb R^1$ in  $B_{\mathbb R^1}(R)^N$, where this image contains two $n$-hemispheres of  radius $R$ and a cylindrical region between them which is fibered over $[a,b]$, where the fibers are equatorial   $(n-1)$-subspheres in the 
$(N-1)$-spheres $S_y^{N-1}(R)\subset \partial B_{\mathbb R^1}(R)^N$, $y\in [a,b]$.
     
     \vspace{2 mm}

   {\bf $\bf \infty$-Circles In Bands}. \footnote{The proof of this presented below was pointed out to me by Anton Petrunin. (compare  circles in 1.A  and with "crcles in discs" from  the previous  section).}  Let $f:S=S^1\to \mathbb R^2$ be a $C^{1,1}$-immersion 
     with curvature $curv^\perp(f) \leq 1$. If  the Gauss map $G_f:S\to S^1(1)$ has zero degree, then {\it width} of the image  \footnote{The width of a subset $X$ in the Euclidean space  is the supremum of  the widths of the bands between parallel hyperplanes  in the space, which contain $X$.}$f(S)\subset \mathbb R^2$ is at least 4, where the equality $width(f(S))=4$ implies that 
     
     {\sl the image $f(S)$, contains 
    the left and the right halves    of the figure {\sf 8}, where these halves are composed of pairs  of semicircles of unit radii.   
    }

      {\it Proof}.  Let $\tilde G_f :S\to \mathbb R$  be a lift of $G_f$  to the universal cover $\mathbb R\to S^1(1)$
and let the function  $\tilde G_f(s)$ assumes its {\it minimum} at  $ s_0\in S$. 

     Then there exist (exactly)  
     two {\it disjoint minimal} segments $S_1,  S_2\subset S$ in the circle, say  $S_1=[s_0,s_1]$ and $S_2=[s_0,s_2],$
   such that $G_f(s_i)=-G_f(s_0)$, $i=1,2$, where "minimal" signifies  that 
   they contain no points  $s$ where $G_f(s)=-G_f(s_0)$, except for their second ends $s_i$ and "disjoint" means $S_1\cap S_2=\{s_0\}$.
  
   Thus, the Gauss map sends $S_1$ and $S_2$ to two disjoint semicircle in $S^1$
   "disjoint" means hat the two meet only at their end points.
   
   Let $H_0, H_1, H_2 \subset\mathbb R^2$ be {\it three lines 
   tangent to} $  f(S^1)$ at the points $s_0, s_1$ and $s_2$ correspondingly, which are,
by their construction,  {\it mutually parallel} and where,
    by the minimality of $\tilde G (s_0)$ and minimality of the segments $S_i$, the only tangent lines to the curves $f(S_i)$,   which, are parallel 
    to $H_0, H_1,H_2$,  are these lines themselves.     
   
   Also. observe that  the {\it tangents} to these curves  at the points  $(s_\circ)_i$, where the Gauss map 
   values $\tilde G_f(s)_i\in S^1(1) $ lie  in the centre of the $\pi$-segment $[G_f(s_0),-[G_f(s_0)]$, are {\it normal} to $H_i$. \vspace {1mm}
   
   Now semi-circle example 8.1.B  applies and the proof completed by applying the above 
  to the {\it maximum} as well the minimum point of the function $\tilde G_f$.

    \vspace {2mm}

\hspace {46mm}{\sc Examples.}

(a) Let $X_o\hookrightarrow \mathbb R^k\subset \mathbb R^{N+k}$ be a smooth closed immersed  submanifold with 
$curv^\perp\leq 1/2$. 

Then the $1$-encircling\footnote  {"$R$-Encircling" is a generalization of "boundary of the $R$-neighbourhood" 
 for    embeddings, see section 3.}
 $X_{o+1}=(X_o)_{+1}\hookrightarrow \mathbb R^{N+k}$ of  $X_o\hookrightarrow \mathbb R^k\subset $ is an smooth immersed hypersurface in the tube
$B_{\mathbb R^k}^N(R)=B^N(R)\times \mathbb R^k\subset \mathbb R^{N+k}$,
 
  (b) If $N=1$ and $B_{\mathbb R^k}^1(R)$ is a band of width 2R rather than a "tube"
and if $dim(X)=k$ then $X_{o+1}$ contain a flat domain inside the boundary of the band. the 
There is much freedom in  deforming with curvature $\leq1/R$, while keeping $X_{o+1}$ within  the band,  for all $k=dim(X)$, but especially for $k=1$.

This is quite different from what happens submanifolds of dimensions$\geq 2$ in tubes 
$B_{\mathbb R^1}^N(R)$ and a similar rigidity is {\color {magenta} probably}  present 
whenever $dim(X)> k$.

(c) Let $X_\circ \hookrightarrow \mathbb R^2$ be the figure $\infty$ curve made of two unit circles (as in 1.C) and let $S^1\times S^{n-1}=X_{\circ \circ} \hookrightarrow\mathbb R^{n+1}$  be obtained by rotating $X_\circ$ around an axes 
$A\subset \mathbb R^2\subset \mathbb R^{n+1}= \mathbb R^{2+(n-1)}$.

If this axes is normal to the line between the centers of the circles, then the image of the immersion  $X \hookrightarrow\mathbb R^{n+1}$ is contained in the unit tube  $B_{\mathbb R^2}^{n+1}(1)$
and if $dist(A, X_\circ)\geq 1$ then $curv^\perp(X_{\circ \circ})\leq 1$.
This $X_{\circ \circ} \hookrightarrow\mathbb R^{n+1}$ is $C^1$-smooth and piecewise $C^2$ smooth 
as in $X_{o+1}$  (b) but the geometry of  $X_{\circ \circ}$ is significantly different from that of $X_{o+1}$. 

These (a)(b)(c) reasonably well represent  immersed hypersurfaces  with curvatures one in the unit "tubes".  $B_{\mathbb R^k}^{n+1}(1)$,especially 
for $k=1$, where  all immersions of closed $n$-manifolds to $B_{\mathbb R^k}^{n+1}(1)$ for $n\geq 2$ are embedding, which are  1-encirclings (boundaries of 1-neighbourhoods) of segments in the line $\mathbb R^1$.



    \section {Veronese Revisited}


Besides   invariant  tori, there are other submanifolds in the unit sphere   
$S^{N-1}$,  which have  small curvatures  and which are transitively acted upon by  subgroups in
the orthogonal   group $ O(N)$.
\vspace {1mm} 

{\it The  generalized  Veronese maps}   are a {\it minimal equivariant isometric}  immersions of spheres to spheres, with respect to certain homomorphisms ( representations)  between the orthogonal
groups $ O(m+1)\to O(m+1)$,
 $$ver=ver_s=ver_s^m:
 S^m(R_s)\to S^m=S^{m_s}= S^{m_s}(1),  $$  
 where  
 $$m_s= (2s+m-1)\frac{s+m-2)!}{s!(m-1!}<2^{s+m}\mbox { and } R_s=R_s(m)=\sqrt \frac {s(s+m-1)}{m},$$
for example, 
$$\mbox {$m_2= \frac{m(m+3)}{2}-1$,    $R_2(m)= \sqrt\frac {2(m+1)}{m}$   and  $R_2(1)=2$,}$$\vspace {2mm}
 (see [DW1971]
If $s=2$ these, called  {\it classical Veronese maps}, are defined  by taking squares of linear
functions (forms) $l=l(x)= \sum_i l_ix_i$  om $\mathbb R^{m+1}$,
$$Ver: \mathbb R^{m+1}\to \mathbb R^{M_m},\mbox { } M_m=\frac {(m+1)(m+2)}{2},$$
where tis   $\mathbb R^{M_m}$ is represented by the space $ \mathcal Q=\mathcal Q(\mathbb R^{m+1})$  of quadratic functions (forms) om $\mathbb  R^{m+1}$,
$$Q =\sum^{m+1,m+1}_{i=1,j=1}q_{ij}x_i x_j.$$

  The Veronese map, which is (obviously)    equivariant  for
the natural action of the orthogonal group  group $O(n+1)$ on  $ \mathcal Q$, where, observe, 
 this action fixes the line 
$\mathcal Q_\circ$ spanned  by the form $Q_\circ=\sum_ix^2$ as well as the complementary subspace $\mathcal Q_\diamond$ of the {\it traceless forms $Q$}, where the action of $O(n+1)$ is irreducible and, thus, it has  a {\it unique, up to scaling} 
Euclidean/Hilbertian structure. 

Then  the normal projection\footnote {The splitting $\mathcal Q=\mathcal Q_\circ \oplus \mathcal Q_\diamond$ is necessarily normal for all $O(m+1)$-invariant  Euclidean  metrics in $\mathcal Q$.}
defines an equivariant map to the  sphere in $\mathcal Q_\diamond$  
$$ver: S^m\to S^{M_m-2}(r)\subset \mathcal Q_\diamond, $$
where the radius of this sphere, a priori, depends on the normalization of the $O(m+1)$-invariant  metric in $ \mathcal Q_{\diamond}$. 

Since we want the map to be   isometric,  we  either take $r=\frac {1} {R_2(m)}=\sqrt\frac {m}{2(m+1)}$  and keep $S^m=S^m(1)$ or 
 if we let $r=1$ and  $S^m=S^m(R_2(m))$ for 
$R_2(m)=\sqrt\frac{2(m+1)}{m}$.

Also observe that the
Veronese maps,   which are not  embeddings themselves,  factor via embeddings of  projective spaces 
to spheres 
$$S^m\to \mathbb RP^m \subset S^{M_m-2}\subset
 \mathbb R^{M_m-1}=\mathcal Q_{\diamond},\mbox { } M_m=\frac {(m+1)(m+2)}{2}.$$

{\it \textbf {Curvature of Veronese.} } Let is show that 

CURvature of veronese  by Petrunin formula

$$curv^\perp_{ver} \left(S^m(R_2(m))\hookrightarrow S^{M_m-2}(1)\right)   =\sqrt { \frac {R_2(1)} {R_2(m)}-1}=\sqrt\frac  {m-1}{m+1}.$$

Indeed, the Veronese map sends equatorial circles from $S^m(R_2(m))$ to planar circles of radii 
$R_2(m)/R_2(1)$, the curvatures of which in the ball $B^{M_m-1}$ is 
$R_2(1)/R_2(m)=2\sqrt\frac {m}{m+1}$  and the curvatures of these  in the sphere,  
$$curv^\perp(S^1\subset S^{M_m-2}(1))=\sqrt {curv (S^1\subset B^{M_m-1}(1))^2-1}=\sqrt {\frac {4m}{m+1}-1}=\sqrt \frac {3m-1}{m+1}$$
is equal to the curvature of  the Veronese $S^m(R_2(m))\hookrightarrow S^{M_m-2}(1)$ itself

$\sqrt{R_2(1)/R_2(m)}=\sqrt\frac {2m}{m+1}$, and the curvatures of these  in the sphere,  
$$curv^\perp(S^1\subset S^{M_m-2}(1))=\sqrt {curv (S^1\subset B^{M_m-1}(1))^2-1},$$
is equal to the curvature of  the Veronese $S^m(R_2(m))\hookrightarrow S^{M_m-2}(1)$itself.
QED.

\vspace {1mm}

It may be hard to prove (conjecture in section 1)   that  {\sf  Veronese manifolds  have  the smallest possible curvatures among  non-spherical $m$-manifold in the unit ball:}
{\sl if a smooth compact $m$-manifold $X$ admits a smooth  immersion to the unit ball  $B^N=B^N(1)$ 
 with curvature $curv^\perp(X\hookrightarrow B^N)<\sqrt \frac {2m}{m+1}$, then $X$ is diffeomorphic to $S^m$.}

It is  more realistic   to show that the Veronese  have smallest curvatures among  submanifolds $X\subset B^N$ {\it invariant under  subgroups in $ O(N)$, which  transitively act on $X$. }

\vspace{1mm} 

{\it Remark.} {\sf  Manifolds $X^m$ immersed to $S^{m+1}$ with curvatures $<1$ are diffeomorphic to $S^n$}, see 5.5,
   but, apart from   Veronese's, we {\color {red!60!black} can't rule out } such $X$  
 in  $S^N$ for $N\geq m+2$ \footnote { Hermitian Veronese maps from the complex projective spaces $\mathbb CP^m$ to the spaces $\mathcal H_n$ of  Hermitian forms on $\mathbb C^{m+1}$  are among the prime suspects in this regard.
 } and,   even less so,  non-spherical $X$ immersible with curvatures $<\sqrt 2$ to $B^{N}(1)$, even for $N=m+1$.

It seems hard to decide this way or another, but it may  be realistic to try to prove  {\it sphericity  of simply connected} manifolds immersed with curvatures $<1$  to $S^N(1)$   for all $N$.

\vspace {1mm}

The curvatures of Veronese maps can be  also evaluated with the {\it Gauss formula,} (teorema egregium),
  which also gives the following formula for curvatures of all $ver_s$:

$m=2$  $1-2c^2=1/3$, $2c^2=2/3$ $c\sqrt{1/3}$

$C=\sqrt{1+1/3}=2/\sqrt 3$

{\it \textbf {From Veronese to Tori.}}  The restriction of the  map 
$ver_s:S^{2m-1}(R_s)\to S^{N_s}$ to the Clifford torus $\mathbb T^m\subset S^{2m-1}(R_s)$  obviously satisfies 
$$curv^\perp_{ver_s}(\mathbb T^m)\leq A_{2m-1,s}+\frac {\sqrt m}{R_s}= 
\sqrt {3-\frac {5}{2}m+\varepsilon(m, s)}$$
 for 
 $$\varepsilon(m, s)= \frac {2}{4m^2}-\frac {4m-2}{s(s+2m-2)}+\frac{5(2m-1)}{2ms(s+2m-2)}-
  \frac {2m-1}{(ms(s+2m-2))^2}.$$

 This, for $s>>m^2$, makes $\varepsilon(m, s)=O\frac{1}{m^2}$

 Since $N_s< 2^{s+2m}$, 
 
 {\it starting from $N=2^{10m^3}$
 $$curv^\perp_{ver_s}(\mathbb T^m)<\sqrt{ 3- \frac {5}{2}m}.$$}
 where it should  be noted that

 {\it  the Veronese maps restricted to the Clifford tori are $\mathbb T^m$-equivariant}\vspace {1mm}

and that

{\sf this bound is {\it  weaker than  the optimal one}  $ \frac{||y||^2_{l_4}}{||y||^2}\geq\sqrt {3-\frac {3}{m+2}}+\varepsilon$ from 
 the previous section}.\vspace {1mm}

{\it Remarks}. (a) It is not hard to go to the (ultra)limit for $s\to \infty$  and thus obtain an 

{\sf equivariant  isometric immersion $ver_\infty$ of the 
Euclidean space $\mathbb R^m$ to the unit sphere in the Hilbert space, such that
$$curv^\perp_{ver_\infty}(\mathbb R^m\hookrightarrow S^\infty)=\sqrt \frac {(m-1)(2m+1)}{(m+1)^2}= 
\sqrt {2-\frac {5}{m+1}+\frac{2}{(m+1)^2}},$$}
where  equivariance is understood with respect to a certain unitary representation of the isometry group of $ \mathbb R^m$. 

{\color {red!50!black}Probably,} one can show that this $ver_\infty$ realizes the {\it  minimum} of the curvatures   among all equivariant  maps $\mathbb R^m\to S^\infty$.\vspace {1mm}

(b) Instead of   $ver_s$,  one could achieve (essentially) the same result with a use of compositions of the classical Veronese maps, $ver: S^{m_i}\to S^{m_{i+1}}$,  $_{i+1} = \frac {(m_i+1)(m_i+2)}{2}-2$,
$$S^{m_1} \hookrightarrow  S^{m_2}\hookrightarrow  ...  \hookrightarrow S^{m_i},$$
starting with   $m_1=2m-1$ and going up to $i=m$. (Actually, $i\sim \log m$ will do.) 

      \subsection {Petrunins Veronese Rigidity Theorem}


{ \bf  Large Simplex  Property}.(Compare with section 5 in pet.) Let the curvature of   a  complete \footnote{"Complete" refers to the induced Riemannian metric .} connected $n$-submanifold in an $n$-ball of radius $r$ be bounded by one, 
$$curv^\perp (X\hookrightarrow B^N(r))\leq 1,$$
and let $x_0,...,x_m \in X$ be $m+1$ points (e.g. m=n), such that 
$$dist_X(x_i,x_j)=\pi, 0\leq i<j\leq m.$$
Then 
$$r\geq\sqrt {2m\over m+1}.$$

{\it In fact,}  the  Euclidean distances between $x_i$ are $\geq 2$ by  
${[2\sin]_{bow}}$inequality, the minimal ball which contains these  point cant be smaller than 
the ball 
 circumscribed about regular $m$-simplex  with the  edge  length 2 by the   Kirszbraun theorem.

{\bf Petrunin's two Balls Covering and the Sphere Theorem.} Let  the $f$-mage of $X$ be contained in the ball of radius 
$r<2/\sqrt 3$ and   let
$x_-,x_+\in X$ be two points joint by a geodesic segment of length $\pi$.
Then the two geodesic balls $B_{x_\pm}(\pi)\subset X$ cover $X$. 

It follows that  $X$ is homeomorphic to  the sphere and,  except for $n= 1$, the map $f: X\hookrightarrow \mathbb R^N$ is an embedding.

{\it Proof}. The above for $m=2$ shows that the boundaries of these balls don't intersect
and since these boundaries are connected  for $n\geq 2$ the balls  do cover $X$.

{\bf Petrunin's   Veronese Planes Rigidity Theorem.} If  
the image $f(X)\mathbb R^N$ is contained  the ball $ B^N(2/\sqrt 3)$ and is not homeomorphic to the sphere
 then $f$ is an embedding and  all geodesic segment in $f(X)$ are planar (contained in planes). 

Consequently, $X$ is either (congruent  to) a  {\it Veronese plane or its complex, quaternionic or  Cayley numbers counterpart.}

{\it Proof .} Track the  two balls covering argument in the extremal case with the bow rigidity at you hand or consult [Petr 2024].

{\it Embedding Remark}.  Petrunin  requires that $f$ is embedding, but this seems {?} unneeded for his argument.

  \section {Hilbert’s Rational Spherical  Designs and Optimal Tori }

   Let $$E:\mathbb R^N\to  B^{2N}(1)\subset \mathbb R^{2N}$$ be the composition of the Clifford embedding $\mathbb T^n\subset  B^{2N}$ 
  and the exponential (locally isometric covering) map
    $$\mathbb R^N=T_0(\mathbb T^N)\overset {\exp} \to \mathbb T^N.$$
 A simple computation  shows [Gro 2022] that the Euclidean  curvature of $E$ on the line $\bar x \subset  \mathbb R^N$ generated by a non-zero vector $x\subset \mathbb R^n$ is
{\color{blue} {$$curv^\perp(\bar x\overset {E}\hookrightarrow  \mathbb R^N)= 
\left ( \frac {\|x\|_{L_4}}{\| x\|_{L_2}}\right)^2,\leqno \mbox {({\huge $\star$})}$$}
where $x=(x_1,...,x_N)$ for the standard Euclidean (corresponding to the cyclic torical) coordinates $x_i$ and 
 $$||x||_{L_p}=\sqrt [p]{\frac  {\sum^N_1 |x_i|^p}{N}}.$$}

 Let $P(n,4)$ be the linear space of homogeneous polynomials of degree
  4 on $\mathbb R^n$, this has dimension  ${n+4 \choose n}=\frac {n(n-1)(n-2)(n-3)}{24}$, and let 
  $$V_+4: \mathbb R^n\to P(n,4), \mbox  { }   V_+4:(c_1,...,c_n)\mapsto (c_1x_1+...c_nx_n)^4$$ be the {\it 4th degree Veronese map}.

 Then A $(n-1)$-spherical $N$-multi-set, that is   map from a set  $\Sigma$  of cardinality  $N$ to the  unit  sphere 
  $S=S^{n-1}\subset \mathbb R^n$ written as $\sigma\overset {\mathcal D} \mapsto s(\sigma)$, is called is a called {\it a design of degree $4$ and cardinality $N$ in $S=S^{n-1}$}  if
 
 {\sl the center of mass of the $N$-multi-set $V_+4D$ in the image $V_+4(S^{n-1})\subset P_(n,4)$ is equal to the center of mass of $V_+4(S^{n-1})$ itself with respect to 
 the usual spherical measure} or, equivalently, if
 $$\frac {1}{N}\sum_{\sigma\in \Sigma}  l^4(\mathcal  D(\sigma))=\int_{S} l^4(s)ds  $$   
  for all linear functions $l$ on $S=S^{n-1}$, where $ds$ is the normalised (i.e, of the full mass one) spherical measure.

  Yet another way to characterise the design property of a muti-set $D$ on $S^{n-1}$
of cardinality $N$ is via  the tautological map
  $$\mathbb R^n=\mathbb R_D \hookrightarrow \mathbb R^N$$ 
 from  the  Euclidean $n$-space of linear functions $l(s)$ on $S^{n-1}$ 
  to the space $\mathbb R^N$ of (all) functions on $\Sigma$.   
 
In these term 
   $\cal D$ is a {\it a design (of degree $4$ and cardinality $N$} in $S=S^{n-1}$)  if and only if --
  this follows  by the standard $\Gamma$-formulas for the  
$\int_{S} l^p(s)ds$-integrals,  
the $L_2$ and the $L_4$ norms on the non-zero  vectors $x\in \mathbb R^N $
which are contained in 
$\mathbb R_{\mathcal D}$
satisfy:
{\color {blue}$$\frac {\|x\|_{L_4}}{ \|x\|_{L_2}}=\sqrt[4]{\frac {3n}{n+2}}$$} 

Thus, in view of {\color{blue} \huge $\star$, }

{\sl every Design $\cal D$ of degree 4 and  cardinality $N$ on $S^{n-1}$ defines a homomorphism (which is a locally isometric immersion), call it $ E_{\mathcal D}$,   from 
$\mathbb R^n=\mathbb R_{\mathcal D}$ to the Clifford $N$-torus, such that 
the curvature of $ E_{\mathcal D}$  in the ball  $B^{2N}(1)\supset \mathbb T^N$ satisfies:
\color {blue}$$curv ^\perp \big (\mathbb R^n\overset {E_{\mathcal D}}\hookrightarrow  B^{2N}(1)\big )=\sqrt\frac {3n}{n+2}.$$}
 
 A design $D$ is {\it rational} if all points in $D$ are rational.

 {\bf Hilbert's  Lemma.}\footnote{In his solution of the Waring problem, Hilbert uses  this lemma  
(for all even degrees)  in the form of  an identity $\sum_{i=1}^N l(x_j)^{2d}=(\sum _{j=1}^N(x_j^2))^d$ for some linear form $l_i$ with rational coefficients.}
  {\sl If $N>>n$, then $S^{n-1}$ contains a rational design of cardinality $N$.}
 
 {\it Proof}. Use three simple  facts.
 
(i)  the center of mass ${\bf c_o}\in P(n,4)= \mathbb R^{{n+4 \choose n}}$ lies in  the {\it interior} of the convex hull of the image  $V_+4(S^{n-1})\subset P(n,4)$
 
 (ii) $\bf c_o$ is a {\it rational} point in  $P(n,4)$,
 
 (iii) rational points in $S^{n-1}$ are dense
 
 and proceed in four steps;
 
(1) Because of (i) and (iii) there exit finitely many  {\it rational} points $s_i\in S^{n-1}$, $i=1,....,M$, such that the convex hull of these ponts contains $\bf c_o' $. 
 
(2) Because of rationality of $\bf c_o$, there exist {\it rational} numbers $p_i\geq 0$, $p_1+...+p_M =1$,  such that $p_1 V_+4(s_1)+...p_M V_+4(s_M)1=\bf c_o$.
 
(3)  Let $Q$ be the common denominator of these numbers and write them as 
 $P_i\over Q$ for integer $P_i$, $i=1,...,M$, where $P_1+...+P_M=Q$.
 
 (4) Let $D$ be the multi-set in $S^{n-1}$, which consists of the points $s_i$, each taken with multiplicity $P_i$.
 
 Then  the center of mass of $V_+4D$ is  
 $$ {1\over Q} \sum_i P_1V_+4(s_i)=\sum_i p_i V_+4(s_i)=\bf c_o.$$
QED.


\textbf {A. $2n^2$-Designs}. The number $N$ delivered. by the above proof is very big, a rough estimate is $N\leq $ but 
non-rational designs are known to exist for much smaller $N$.

For instance {\it If  $n$ is a power of 2, then  there exists a   design of cardinality  $N=2n^2+4n$.}
 \footnote {This was stated and  proved in a  written  message by Bo'az Klartag to me.  Also, Bo'az pointed out to me that 
   the  Kerdock code used in   [Kon 1995] yields  designs for $N=4^k$ and $N=\frac {n(n+2)}{2}$.}
  
[K1995] H. Konig, Isometric imbeddings of Euclidean spaces into finite di-
mensional lp -spaces, Banach Center Publications (1995) Volume: 34, Issue: 1,
page 79-87.

  {\sl homomorphism,   (which is a locally isometric immersion) from
 the Euclidean $n$-space to the Clifford $N$-torus in the ball $B^{2N}$ for 
 $N=8(n^2+n),$
  such that the normal Euclidean curvature  of this immersion is   
 {\color {blue} $$curv^\perp\big (\mathbb R^n\hookrightarrow B^{16(n^2+n)}(1)\big)=\sqrt\frac {3n}{n+2}\leqno \mbox {({\huge $\star \star$})}$$}}
  
      Since rational points are dense in the sphere, we conclude to the extence of 
   subtori $\mathbb T^n_\varepsilon \subset \mathbb T^{ 8(n^2+n)}$, such that  
   {\color {blue}$$curv^\perp\big (\mathbb T_{\varepsilon}^n\hookrightarrow B^{16(n^2+n)}(1)\big)\leq\sqrt\frac {3n}{n+2}+\varepsilon \mbox { for all $n$ and all $\varepsilon>0$}.\leqno \mbox {({\huge $\star \star\star$})}.$$}

[K1995] H. Konig, Isometric imbeddings of Euclidean spaces into finite di-
mensional lp -spaces, Banach Center Publications (1995) Volume: 34, Issue: 1,
page 79-87

    {\sl if $N>>n$ as in Hilbert's Lemma, then there exist $n$-subtori  $\mathbb T^n \subset B^{2N}$,  fsuch that} 
 {\color {blue} $$curv^\perp \big(\mathbb T^n\hookrightarrow  B^{2N}\big )=\sqrt\frac {3n}{n+2}.$$}

 {\it Example/Non-Example.} Regular pentagons serve as   designs of cardinality five and   degree four on the circle; these are irrational and  there is no
  apparent {\it simple} rational design on $S^1$.


  \section{Link with the Scalar Curvature via the  Gauss Formula}


  The  $curv^\perp$  problem came up in the context  of Riemannian geometry of manifolds $X$ with {\it positive scalar curvatures} [Gro 2017], where 
   
   {\sl the scalar curvature of an $X$ at $x\in X$,  denoted $Sc(X,x)$, is 
     the sum of the values of the {\it sectional curvatures} $\kappa$ at the $n(n-1)$  (ordered) orthonormal  bivectors in  $T_x(X)$, for $n=dim(X)$}.\footnote{One knows 
 that  $Sc(X,x)>0$ if and only if the volume of the ball $B_x(\varepsilon) \subset X$ is smaller than the volume of the 
     $\varepsilon$-ball in  $\mathbb R^n$, {\it provided $\varepsilon>0$ is  sufficiently small}:   $\varepsilon\leq \varepsilon(X,x)>0$.   Albeit  looking explanatory,   this is   only an  illusion of understanding the geometric  meaning of the inequality $Sc(X)>0$.}      
     \vspace {1mm}

 For instance, scalar curvatures of  surfaces are equal to twice their sectional (Gauss) curvatures.

\vspace {1mm}

{\bf Spheres Example.}  The $n$-spheres of radii $R$  in
  the Euclidean space  $\mathbb R^{n+1}$  (which have constant sectional curvatures $1/R^2$), satisfy:
      {\color {blue}$$Sc(S^n(R))=n(n-1)/R^2 \mbox {  for all } n.$$ }

 {\bf Additivity.} It follows from the definition that the scalar curvature is additive under Riemannian products,
 {\color {blue}$$Sc(X_1\times \underline X)=Sc(X_1)+Sc(\underline X).$$}
  
 For instance, the scalar curvature of the  $n$-th power of the unit 2-sphere is 
 $$Sc(\underset {n}{\underbrace {S^2\times S^2\times ...\times S^2}})=2n=Sc(S^{2n}(R=\sqrt{ {2n-1}})$$
 
This also shows  that    the topology of manifolds with positive scalar curvatures of dimensions $n\geq 4$, can be arbitrary complicated\footnote {Three manifolds with $Sc>0$ are not too  simple either :  {\it connected sums  lens spaces} and copies of $S^1\times S^2$ admit metrics with $Sc>0$ by a theorem by Schoen and Yao.} for
 {\color {blue}$$Sc(X\times S^2(\varepsilon ))\underset{\varepsilon \to 0}\to +\infty \mbox {  for {\it all compact} Riemannian  manifolds $X$.} $$}

 
 
 Yet, there are limits to this complexity: there are compact manifolds of all dimensions, which admit no metrics with $Sc>0$,  called {\color {red!50!black}$\mathbf\nexists$}PSC,
 where the three basic examples are as follows.\vspace {2mm}

  \hspace {40mm}{\sc Basic  {\color {red!50!black}$\mathbf\nexists$}PSC Manifolds}\footnote{See [Gro 2021] for a survey on  topological properties and examples of  $\nexists$PSC Manifolds.}
  \vspace {1mm}

{\bf 14. A. Lichnerowicz Theorem.}   The (Kummer) surface defined by the equation  $z_1^4+ z_2^4+ z_3^4+z_4^4=0$  in the complex projective space $\mathbb CP^3$ and, more generally orientable spin manifolds with non vanishing $\hat A$ genus 
(dimensions of these are   multiples of $4$) admit no Riemannian metrics with $Sc>0$.
  
 {\it Proved} in 1963 with the first (1963) {\it  Atiyah–Singer index theorem for the Dirac operator}.

 {\bf 14. B. Hitchin theorem}: there exist manifolds $\Sigma$ homeomorphic (but non-diffeomorphic!) to the spheres $S^n$ for all $n = 8k+1,8k+2$,  $k = 1,2,3...,$ which admit no metrics with $Sc > 0.$

{\it Proved} in 1974 with the second (1971)  {\it Atiyah–Singer index theorem}. 

{\bf 14.C. Geroch Conjecture.} $n$-Tori admit no metrics with $Sc>0$.

 {\it Proposed} in 1975,  {\it proved} in 1979 by  Schoen-Yau  for $n\leq 7$ with via {\it minimal hipersurfaces} by induction on $n$ and  by Gromov-Lawson in 1980 for all $n$ with the {\it index theorem for the Dirac operators  twisted with almost
 flat bundles.}
 \vspace{1mm}
 {\bf 14.  D. Product Manifolds.} Products of the above manifolds, e.g. of tori by Hitchins spheres are also   {\color {red!50!black}$\mathbf\nexists$}PSC.
 
 This is proven with the {\it index theorem} for the (generalized) Dirac operators.

 {\it Sectional Curvature Remarks.} Although  the   inequality   $Sc>0$ is much weaker then $sect.curv>0$ (which is equivalent to  geodesic triangles having the sums of the angles $>\pi$)
  {\it no alternative proofs} of non-existence of metrics with $sect.curv>0$ on manifolds from 
  {\bf A} and {\bf B} are available,
 while  the  $sect.curv>0$ (and  $Ricci>0$) version of {\bf C}     follows by an  elementary argument
  relying on the geometry of geodesics in $X$.
  
  (The ancient Bonnet-Myers theorem says that  
  $Ricci (X)\geq \kappa>0 \implies  diam(X)\leq \sqrt {1/\kappa}$, which rules out closed   manifolds with infinite universal coverings, such as tori.) 
  
  \vspace {1mm}

{\it Turning to  Constant Sectional Curvature}. If one requires  the strongest possible condition of this kind, namely the   sectional curvature  to be constant as well as  positive, then   everything about $X$ appears   100\% transparent.

Indeed, one knows. that     these  metric are {\it locally spherical}; hence all simply connected $n$-manifold $X$ with  $sect.curv(X)=\kappa>0$ admit {\it locally isometric immersions to $S^n(R)$ for $R=\sqrt{1/\kappa}$}.

Consequently, 

{\it the universal coverings of {\it closed} (compact without boundaries)   manifold $X$ with $sect.curv(X)=\kappa$  are isometric to $S^n(R)$.}
 This is the end of the story. \vspace {1mm}

 Yet, this may be hard to believe, there are {\it \color {blue} non-trivial links between geometry and topology} of manifolds  $X$ with {\it constant sectional curvatures}  if these $X$ have {\it non-empty boundaries}, where the available proofs of  these properties 
 rely on   the scalar curvature  inequality  $Sc (X) \geq n(n-1)/R^2$ and  where  one doesn't  know how to  exploit   to full power of the condition 
 $sect.curv=const =1/R$ (see section 15)


 \subsection {Gauss Formula and Petrunin's Curvature}


  Let $X\subset Y$ be a smooth n-dimensional submanifold in a Riemannian $N$-manifold, e.g. in    $Y=\mathbb R^N$ and let II=  II(X,x)II=II$(\tau_1,\tau_2)$ 
  be the {\it second fundamental form} (corresponding  to the {\it shape operator}) of $X$ at $x\in X$, 
   where $\tau_1,\tau_2\in T_x(X)$ are tangent vectors to $X$ and the form II  takes  values in the normal space $T_x^\perp(X)$ and where  
  II$(\tau,\tau)$ is equal to the second derivative of the geodesic in $X$ issuing from $x$ with the velocity  $\tau$.
   
 The normal curvature of  $X\subset Y$ at $x\in X $,  in these terms  is
   $$curv^\perp_x=\sup_{ \|\tau\|=1}\|\rm II(\tau,\tau)\|.$$

  {\it The  $l_2$-norm of II at $x$} is 
 
    $$\|\rm II\|^2_{l_2}=
    \sum_{i_1, i_2=1,...m}\|{\rm II}(\tau_{i_1},\tau_{i_2})\|^2,$$
  where  $\{\tau_i\}$, $i=1,...,n=dim(X),$   is a frame of  orthonormal vectors  in the tangent space $ T_x(X)$.

  We shall need    the  simple inequality 
   $$\|{\rm II}\|^2_{l_2}
  \leq kn\cdot curv^\perp(X)^2,  $$
  which is useful for $k<n$.
    One can also show  that     
      $\|{\rm II}\|^2_{l_2}\leq  n^2\cdot curv^\perp(X)^2,$ for all $k$ but
       the following  inequality. will serve us better.

  { \it Petrunin curvature}  ${\sf \Pi=\Pi_x}(X\subset Y)$ is the average of 
 $$\|{\sf II}(\tau,\tau)\|^2$$
 over the unit vectors $\tau\in  S^{m-1}_x\subset T_x(X)$, where clearly,
  $$
  \frac {(curv^\perp_x)^2}{n-1}\leq  {\sf \Pi}_x\leq (curv^\perp_x)^2
  $$ 
  and where the equality $\frac {(curv^\perp_x))^2}{n-1}=  {\sf \Pi}_x$
  holds if the form II has rank one  and   
  ${\sf \Pi}_x= (curv^\perp_x)^2$
   if $||$II$||^2_{l_2}=||mean.curv(X,x)||^2$.

For instance, if  $codim(X)=1$, the latter means that all principal curvatures  $X$ at $x$ are mutually equal.

 More interestingly   [Petr 2023]) 
  $$ \mbox {${\sf \Pi}={2\over n(n+2)} \big(||$II$||^2_{l_2}+{1\over 2}|\hspace{-0.2mm}|mean.curv^\perp|\hspace{-0.2mm}|^2\big)$}$$
   or 
  $$\|mean,curv\|^2-\|{\rm II}\|^2_{l_2}={3\over2} mean,curv^2-\frac {n(n+2)}{2}{\sf \Pi},$$
   which is proven with the {\it same}  $\Gamma$-function formula for the integrals of  polynomials of degree four on $S^{n-1} $,   which goes along with spherical designs 
 and   used for construction of immersions $\mathbb T^n\to \mathbb R^N$ with 
   $curv^\perp =\sqrt {3n/(n+2)}+\varepsilon$.
  
   (One wanders if there is  a geometric reason for this,  e.g.  a "Riemannian curvature  averaging formula" of some kind.)
   

  For instance, if  $n=dim(X)=2$, $N=dim(Y)=3$ and  $\alpha_1$ and   $\alpha_2$ dnote the  principal curvatures of $X$ at $x$,
  then 
 $$curv^\perp(X,x)=\max( |\alpha_1|,|\alpha_2|),$$
  $$\mbox{ $||$II$||^2_{l_2}=\alpha_1^2+\alpha_2^2$},$$ 
   $$|\hspace{-0.2mm}|mean.curv^\perp|\hspace{-0.2mm}|=|\alpha_1+\alpha_2|$$
   and $${\sf \Pi}={1\over 4} ( \alpha_1^2+\alpha_2^2)+{1\over 8} (\alpha_1+\alpha_2)^2={3\over 8}( \alpha_1^2+\alpha_2^2)+\frac{1}{4} \alpha_1\alpha_2;   $$
   if
   $X=S^2\subset Y=\mathbb R^3$, where  $\alpha_1=\alpha_2=1$, this  makes 
${  \sf \Pi}=1$ as well.

   \vspace {1mm}
   
  {\bf Gauss Formula.} Let $Y$ have constant sectional curvature $\kappa$ and let 
  $Sc_{|n}=Sc_{|n}(Y) =nk(k-1)$. Then the scalar curvature of $X$ satisfies:
  $$\mbox {$Sc(X,x) = Sc_{|n}+ ||mean.curv^\perp(X,x)||^2-\|\rm I I\|^2_{l_2}$}, $$ 
 where by Petrunin's formula 
 $$Sc(X,x) = Sc_{|n}+ {3\over 2}\|mean.curv(X,x)\|^2-\frac{n(n+2)}{2}\cdot{\sf \Pi}, $$ 
 
Hence, the inequality  $Sc_{|m}(Y)\geq \sigma_n$ implies that 
$$\mbox {$Sc(X) \geq \sigma_n-||$ II$(X,x)||^2$}.$$
Therefore 
$$Sc(X) \geq \sigma_n- kn\cdot curv^\perp(X)^2\leqno [kn]$$ 
for $k\leq n$
and 
$$Sc(X) \geq \sigma_n-n^2curv^\perp(X)^2.$$
for all $k$, 
where  Petrunin's formula yields  better, in fact optimal,  inequality    for $k>>n$
$$Sc(X) \geq \sigma_n-{n(n+2)\over 2}{\sf \Pi} \geq \sigma_n-{n(n+2)\over 2}   curv^\perp(X)^2.$$

It follows that 
{\it  if  the manifold $X$ is {\color {red!50!black}$\mathbf\nexists$}PSC, i.e. it admits no metric with $Sc>0$, then 
\color {blue}$$ curv^\perp(X)\geq \sqrt {\sf \Pi} \geq  \sqrt \frac {2\sigma_n} {n(n+2)}\mbox {   for all   $k$  and   }  N=n+k=dim (Y),\mbox { }  Y\hookleftarrow X,$$
and   
 $$ curv^\perp(X)\geq \sqrt \frac {\sigma_n}{kn}\mbox {  for $k<n/2$}.$$}        \vspace {1mm}
 
                    \hspace {31mm} {\sc Examples and  Corollaries.}
                    \vspace {1mm}

 {\sl Let $X$ be an $n$-dimensional {\color {red!50!black}$\mathbf\nexists$}PSC manifold, 
e.g. the $n$-torus $\mathbb T^n$,  Hitchin's exotic $n$-sphere $\Sigma^n$ or a product 
$\Sigma^m\times  \mathbb T^{n-m}$.}

 ({\Large$\ast$}$_{S^{{n+k}}}$)  {\sf Then immersions from $X$ to the unit sphere satisfy
\color {blue} $$ curv^\perp\big (X\hookrightarrow S^{n+k}(1)\big)\geq
 \sqrt \frac {n-1}{k}\leqno{\bf [A]} $$
and 
$$ curv^\perp\big (X\hookrightarrow S^{n+k}(1)\big)\geq\sqrt{\sf \Pi}  \geq \sqrt  \frac {2 n-2}{n+2}.\leqno{\bf [B]}$$}
 
Inequality  {\bf [A]} is better than {\bf [B]}  roughly for   $ k \leq n/2$, while Petrunin's  {\bf [B]} takes over for larger $N$, where it is, as we known (see sections 2 and 14) , optimal for $k>>n^2$. 


\subsection { Petrunin's $\sqrt 3$ Extremality Theorem }


The above doesn't directly apply to immersions to 
the Euclidean balls, since these have $Sc_{|n}=0$, where  the Gauss and Petrunin formulas 
for the induced metric $g$, reduce to 
 $$Sc(g) = ||mean.curv^\perp||^2-\|\rm I I\|^2_{l_2}\leqno{\bf[a]}$$
 and
 $$Sc(g) = {3\over 2}\|mean.curv\|^2-\frac{n(n+2)}{2}{\sf \Pi}. \leqno{\bf[b]}$$ 

Yet, inequality  {\bf [A],} applied to the image   of $X\hookrightarrow B^N(1)$ in $S^N$ under the radial projection of 
of the unit ball in  tangent hyperplane  $B^N\subset \mathbb R^N=T_s(S^N)\subset \mathbb R^{n+1}\supset S^N$ to $S^N$
shows that
{\color {blue} $$ curv^\perp\big (X\hookrightarrow B^{n+k}(1)\big)\geq
 \sqrt \frac {n-1}{(8-\varepsilon_{n,k})k}.\leqno{\bf [{1\over 8-\varepsilon}]}
 $$}
 for some  (moderately small) $\varepsilon_{n,k}>0.$
 
 This is  crude, but 
  in the $\sf \Pi$-case Petrunin proves the sharp $curv^\perp$-inequality 
{\color {blue} $$ curv^\perp\big (X\hookrightarrow B^N(1)\big)\geq 
\sqrt { \frac {3 n}{n+2}}. \leqno{\mbox {\bf [B{\Large$\star$}]}}$$}
{\it  for all  $n$-dimensional {\color {red!50!black}$\mathbf\nexists$}PSC manifolds $X$, all
$n$ and $N$.}

This is done by  showing that 
{\it if
$$curv^\perp\big (X\overset {f}\hookrightarrow B^N(1)\big)< \sqrt  \frac {3 n}{n+2},$$
then a conformal change of the induced metric $g$ on $X$  has positive scalar
curvature. Namely, if $n\geq 3$\footnote{If $n=2$ then the average value of $\sf \Pi$ is $\geq\sqrt {3\over 2}$, see [Petr 2023]}, then 
$$Sc(u^\frac {4}{n-2}g)>0\mbox { for $u(x)=\exp -l{1\over 2}\|f(x)\|^2$ and $l={3\over 4}\cdot\frac {n-2}{n-1} \cdot n$}.$$} 
 
{\it Remark.} One might think, that Petrunin's argument  with the Gauss formula 
$$Sc(g) = \|mean.curv\|^2-\|\rm I I\|^2_{l_2}\geq  \|mean.curv\|^2 -k(curv^\perp)^2$$
rather than Petrunin's $$Sc(g) = {3\over 2}\|mean.curv\|^2-\frac{n(n+2)}{2}{\sf \Pi}\geq {3\over 2} \|mean.curv\|^2 -\frac{n(n+2)}{2}(curv^\perp)^2$$ 
would improve the above inequality $[{1\over 8-\varepsilon}]$.

In fact, if one uses 
 Petrunin's formula for the Laplace operator  $\Delta=\Delta_g$ applied to the above function $u(x)$
on $X$:
$$ -\frac {\Delta u}{u}= lrc \cdot |H|  +(ln -l^2r^2 s^2) ,$$
where 
$H=mean.curv \big(X\overset {f}\hookrightarrow B^{n+k}(1)\big)$,   $r=r(x)=\|f(x)\|$,  and  $c=c(x)$, $s=s(x) $  are  function ($\cos$ and $\sin$ of certain angles), which are  bounded in the absolute values by one,
$|c|, |s|\leq 1$,
 one arrives at the following  version of $[{1\over 8-\varepsilon}]$.
:
 $$curv^\perp\big(X\hookrightarrow B^{n+k}(1)\big)\geq\sqrt \frac {n}{k(8+ (4/(n-2)))}   $$
This is no better $[{1\over 8-\varepsilon}]$. but   can be slightly improved  with the inequalities 
$c^2+s^2\leq 1$ and $r^2+s^2\leq 1$ proved in [Petr 2023r
under the assumption $curv^\perp\leq 2$.\vspace {2mm}


\subsection {   Lower Bounds on $curv^\perp (X\hookrightarrow Y)$  
for  Manifolds $Y$ with $Sc_n\geq \sigma_n$. } \vspace {1mm}


 Let us define the $n$-dimensional scalar curvature
  $Sc_n(Y)$ for general
Riemannian manifolds  $Y$ of dimension $N\geq n $, that is  a function on the tangent $n$-planes $T_y^n\subset T(Y) $ in  $Y$, which is eual to   the sum of the sectional curvatures $\kappa$ of $Y$ on the bivectors in  $T_y^n$  at $y$.

   Equivalently, $Sc_n(Y, T_y)$ is   the scalar curvature of the  submanifold $\exp(T_y) \subset Y$ at $y$,  that is is the germ of the image of the exponential map from $T_y$ to $Y$.

Then  the Gauss' and  Petrunin's formulas for the scalar curvature of $X\hookrightarrow Y$ remains as they were for manifolds $Y$ with constant sectionl curvatres 
$$\mbox {$Sc(X,x) = Sc_{|m}(Y,T_x(X))+ ||mean.curv^\perp(X,x)||^2-\|\rm I I\|^2_{l_2}$}, $$ 
and
$$\mbox{$ ||mean.curv^\perp(X,x)||^2-||$ II$(X,x)||_{l_2}^2=||{3\over 2}mean.curv^\perp(X,x)||^2- {n(n+2)\over  2}{\sf \Pi}$}. $$ 
 
 Thus, the above inequalities {\bf [A]} and {\bf [B]} concerning immersions of $n$-manifolds  $X$ to the unit sphere  $S^{n+k}$ generalize to immersions to 
$(n+k)$-dimensional manifolds $Y$, such that $Sc_n(Y)\geq n(n-1)$:

{\color {blue} $$ curv^\perp\big (X\hookrightarrow Y\big)\geq
 \sqrt \frac {n-1}{k}\leqno{\bf [A_Y]} $$
and 
$$ curv^\perp\big (X\hookrightarrow Y\big)\geq\sqrt{\sf \Pi}  \geq \sqrt  \frac {2 n-2}{n+2}.\leqno{\bf [B_Y]}$$}

   {\it Example.} Let  $Y=S^{n+k_0}(R)(1)\times H^l_{-1}$, where the sphere $S^{n+k_0}(R)$
   has constant curvature  $+1/\rho^2$ and  $H^l_{-1}$ is the hyperbolic space with  
   the sectioanal curvature $-1$ and let $n\geq l+2$.
    Then
        $$Sc_n(Y) \geq {1\over \rho^2}(n-l)(n-l-1)-l(l-1)$$
   and the two above inequalities hold with $k=k_0+l$, if 
   $$\rho^2\leq \frac{ (n-l)(n-l-1)}{n(n-1)+l(l-1)}.$$
    
    For instance, if $l=2$, and $n\geq 4$ one needs
    $\rho^2\leq\frac{ 1} {7}$. for this purpose.\vspace {1mm}
    
    Notice in conclusion, that neither 
    
    the above inequalities $ [{1\over 8-\varepsilon}]$ and  Petrunin's ${\bf [B{\Large\star}]}$  for immersion to unit balls  
    
    nor such inequalities from the previous sections  based on the ${2p\over n}$inequalities
    
    admit (not at lest obvious) counterparts for these $Y$.


\section{ Second Link with the scalar Curvature: Width Inequalities for Riemannian Bands}


{\textbf {15.A. Example: Torical  $\frac {2\pi}{n}$-Inequality.}} \footnote {See [Gro 2021]  for an account on known results in   geometry of manifolds with $Sc\geq 0$, which are formulated in sections 16.1, 16.2, 16.3 without further references.}
Let $V$ be a Riemannian manifold  homeomeorphic to the product of the  $n$-torus by the unit interval
$V=\mathbb T^{n}\times[-1,+1]$, such that 
$Sc(V)\geq \sigma>0$.  
Then the distance between the two components of the boundary of $V$ is bounded as follows:
{\color {blue}$$dist( \mathbb T^{n}\times \{-1\},\mathbb T^{n}\times \{+1\})\leq 2\pi\sqrt\frac { n}{\sigma (n+1)}.$$}
(See section 16.1 for a few words about   the proof.)

{\bf 15.B. Corollary: No Wide Torical Bands in the Spheres. } {\it If a Riemannian $(n+1)$-manifold $V$  homeomorphic to $\mathbb T^{n}\times[-1,+1]$ admits a locally  isometric immersion  to the  $(n+1)$-sphere of radius $R$
then 
{\color {blue} $$dist( \mathbb T^{n}\times \{-1\},\mathbb T^{n}\times \{+1\})\leq \frac { 2\pi R}{ n+1}.$$}}

{\bf  15.C.  Large Normal Curvature Sub-corollary.}
 {\it Let 
 $$f:\mathbb T^n\hookrightarrow B^{n+1}(1)$$
be a smooth immersion from the $n$-torus to the unit Euclidean $(n+1)$-ball $B^{n+1}\subset \mathbb R^{n+1}$.
 Then the curvature of $f$ 
 is bounded from below by:}
{\color {blue}$$ curv^\perp\big (\mathbb T^n \overset {f}\hookrightarrow B^{n+1}(1)\big)\geq \frac {n+1}{\pi}-1.$$}
  
 {\it Proof of {\bf 15.B $\implies$  15.C.}}  Let 
 $$E_f: \mathbb 
 T^N\times \mathbb R^1\to \mathbb R^{n+1}\supset B^{n+1}(1)$$ 
 be the {\it normal exponential map}, i.e. such that the restriction  $E_f|\mathbb T^N\times \{0\}=f$ and where $E_f$  isometrically sends   the lines $\{t\}\times \mathbb R^1$, $t\in \mathbb T^n$, to the  straight lines in  $\mathbb R^{n+1}$ normal to the immersed torus $f(\mathbb T^n)\subset \mathbb R^{n+1}$ at the points  $f(t)\in f(\mathbb T^n)$. 
 
If $curv(f)<c$, then, (this is the same as it is for circles of radii $1/c$ in the plane)    the map $E_f$ is an {\it immersion} on 
 $\mathbb 
 T^N\times [-r,r] \subset 
 T^N\times \mathbb R^1$ for $r=1/c$, while the image of $f(\mathbb T^n)$ is contained in the ball $B^{n+1}(1+r)$.

  Let  
  $$\mathbb R^{n+2} \subset  S^{n+1}_+(1+r)\overset {p}\to\mathbb R^{n+1} \supset B^{n+1}(1+r)$$ 
    be the normal projection from the hemisphere, 
     compose $E_f$ on  $\mathbb 
 T^N\times [-r,r]$ with the inverse map to $p$  and let 
    $$\tilde E: p^{-1}\circ E_f: \mathbb T^N\times [-r,r] \to S^{n+1}_+(1+r).$$

    Since the projection $p$ is {\it distance decreasing}, the spherical distance between the two components of the boundary of $T^N\times [-r,r]$ with respect to  the Riemannian metric $\tilde g$ in $\mathbb T^N\times [-r,r]$  induced by $\tilde E$ from  the spherical metric in $S^{n+1}_+(1+r)$
    $V$ is bounded from below by $2r$. Then 
 {\bf D} applied to
 $$(\mathbb T^N \times [-r,r],\tilde g)\overset {\tilde E}\to S^{n+1}_+(1+r)\subset S^{n+1}(1+r)$$
shows that 
$$\tilde d=dist_{\tilde g}( \mathbb T^{n}\times \{-r\},\mathbb T^{n}\times \{+r\})\leq 
\frac { 2\pi (1+r)}{ n+1}$$
and since $\tilde d> 2r=2/c$ the inequality 
$c\geq \frac {n+1}{\pi}-1$ follows. QED.

 \vspace{1mm}
    
    {\it Exercise.} Generalise the large normal curvature sub-corollary to immersions of tori to
   products of balls: 
{\color {blue}$$ curv^\perp\big (\mathbb T^{n+k} \overset {f}\hookrightarrow B^{n+1}(1)\times B^k(R)\big)\geq \frac {n+1}{\pi}-1.$$}
for all $k=0,1,2,...$ and all $ R\geq 0$.

{\it On Low Dimensions. }  The inequality $curv^\perp\big ((\mathbb T^n\hookrightarrow B^{n+1}(1)\big)\geq \frac {n+1}{\pi}-1$ may be asymptotically optimal  for $n\to \infty$ but  its performance 
for small $n$ is poor.

For instance, if $n\leq 5$ then $\frac {n+1}{\pi}-1<1$ and our inequality is weaker than
$curv^\perp(X^n\hookrightarrow B^{n+k}(1)\geq 1$, which follows for all closed $n$-manifolds $X$  and all  $n,k$ by the obvious "maximal principle" argument.

Furthermore,  since
$$curv^\perp\big ((X^n\hookrightarrow B^{n+1}(1)\big)> 2$$ 
for all {\it non-spherical} $X$ (this is elementary, see section ...), our $\left(\geq\frac {n+1}{\pi}-1\right )$-bound  is of any interest only for $n\geq 9.$

{\it $\mathbb T^\rtimes$-Remark.} In section 15.2   we introduce the notion of 
 {\it $\mathbb T^\rtimes$-stabilized  scalar curvature,} $Sc^\rtimes (X)$, improve the inequalities {\bf E} and {\bf F} and will see, for example, that 
$$curv^\perp\big (\mathbb T^n\hookrightarrow B^{n+1}(1)\big)>2.5\mbox { for $n\geq 7$}.$$

\vspace{1mm}

 {\it Codimension two Remark.}  The inequality $\bf E$ applied to the unit tangent bundles of immersed $n$-tori  with codimensions  2,\footnote {If the Euler class of such an immersion is non-zero one needs a mild generalisation of {\bf E}.}    shows  (see  $[1+2c]$-Example in section 3)  
  $$curv^\perp \big (\mathbb T^{n+1} \hookrightarrow  B^{n+2}(1)\big)\leq 1+2curv^\perp\big (\mathbb T^{n} \hookrightarrow  B^{n+2}\big)$$
and 
{ \color {blue}$$curv^\perp\big (\mathbb T^{n} \hookrightarrow  B^{n+2}\big)\geq{1\over 2}curv^\perp \big (\mathbb T^{n+1} \hookrightarrow  B^{n+2}(1)\big)-{1\over2}\geq \frac {n+2}{2\pi}-1. $$}

 This  has any merit   only  for $n\geq 11$, where $\frac {n+2}{2\pi}-1>1$, and it becomes better than  Petrunin's inequality only  for  
 $n\geq 15$, where
$\frac {n+2}{2\pi}-1>\sqrt{3n\over n+2}$.  

(The   improvement   with  the $\mathbb T^\rtimes$-remark doesn't significantly change the picture.)
 
 {\bf  Conjectures}. (a)  Immersed $n$-tori in the unit $ (n+k)$-ball satisfy 
 {\color {magenta}$$curv^\perp \big(\mathbb T^n\hookrightarrow B^{n+k}(1)\big )
 \geq{ n\over k}.$$}

(b)  All immersions of all  $n$-manifolds  $X\overset {f_0}\hookrightarrow  B^N(1)$ a regularly 
homotopic to immersions $f_1:X\to  \hookrightarrow  B^N(1)$ where
$curv^\perp(f_1:X) \leq Cn$  for someuniversal constant $C$,
(probably, $C\leq 100)$.) 

These , by no means (not even conjecturally) optimal,  inequalities  are motivated only
by their simple forms.


{\bf  15.D. Immersions with curvatures  $\sim n^\alpha$.} It {\color {magenta}not impossible (but unlikely)} that all immersion of $n$-tori to unit  balls satisfy
{\color {magenta}$$curv^\perp \big(\mathbb T^n\hookrightarrow B^{n+k}(1)\big )
 \geq  {cn^\alpha\over k}$$}
for some  small $c> 0$, $\alpha >1$, e.g. $c=0.001$ and $\alpha ={3\over 2}$,
where the exponent $\alpha ={3\over 2}$ is maximal possible.

Indeed,
$n$-tori embed to  $B^{n+n}(1)$ with curvatures $ n^{1\over 2}$  and 
also there exit
 codimension one embedding of $n$-tori with curvatures about $n^\frac {3}{2},$ 
{\color {blue}$$curv^\perp\big (\mathbb T^n\subset B^{n+1}(1)\big)<6n^\frac {3}{2}.$$}

In fact, arguing as in  section 4.A one  construct 
 $X_m=S^{n_1}\times...\times  S^{n_m}\subset B^{{n_1+...n_m}+1}(1)$  by induction on $m$ as boundaries of $\rho_m$-neighbourhoods of 
 $$X_{m-1}=S^{n_1}\times...\times  S^{n_{m-1}}\subset B^{{n_1+...n_{m-1}}+1}(1-\rho_m)\subset B^{{n_1+...n_m}+1}(1),$$ 
 where the   curvatures of these embeddings grow {\it exponentially} with $m$, roughly as  $2^{m-1}$.

Thus    one  embeds     $X_m$ to the ball $B^{{n_1+...n_m}+1}(1)$  with   the curvature growing {\it polynomially} in $n=dim(X_m)$ (rather than in $m$):
  {\color {blue}$$curv^\perp\big (X_m\subset  B^{n+1}(1)\big)\leq const_\mu n^\frac {\mu+2}{\mu+1}, \mbox { $n=dim(X_m)=n_1+...n_m$, $\mu=\min_in_i$. }$$} 
   
{\color {magenta} For all we know}, if all $n_i$ are equal to a single $n_o$, then all immersions of $(S^{n_o})^m$
 immersions to the unit $(mn_o+1)$-ball satisfy

{\color {magenta}$$curv^\perp\big ((S^{n_o})^m\hookrightarrow  B^{mn_o+1}(1)\big)\leq const_{n_o} (mn_o)^\frac {\mu+2}{\mu+1}. $$}


\subsection {On Three Proofs of $\frac {2\pi}{n}$-Inequalities}



All three  proofs apply to  manifolds $V$,  where their boundaries are decomposed into  two disjoint parts 
  $\partial V=\partial_-\sqcup\partial_+ $, and show that 
  $$dist(\partial_ -,\partial_+)> 2\pi\sqrt\frac { n}{\sigma (n+1)}\mbox { for  $\sigma=\inf_{x\in X} Sc(X,x).$}$$ 
 under certain topological assumptions on $V$ specific to each proof.

 {\bf 1.} The first proof  applies  to suitably {\it enlargeable} manifolds\footnote  {See [GL 1973], [Gro 2021] and references therein.}  $V$, e.g.  to $V=X\times [-1,1]$, where $X$ admits a metric with 
 $sect.curv\leq 0.$

 This proceeds by induction on $n$ with minimal hypersurfaces with boundaries as in  \S12 from GL, where  the  original Schoen-Yau argument
 was augmented  with Fischer-Colbrie\&Schoen warped product symmetrization 
 idea.
  
If $dim(V)>7$, the proof  encounters  a technical difficulty  where minimal hypersurfaces may have singularities, but this can be  resolved modulo the partial regularity theorem 4.6 from [SY 2017].
 
  {\bf 2.} The second proof 
   whenever applies, delivers 
   a hypersurface ({\it $\mu$-bubble})  $X\subset V$ which separates  $\partial_-$.  from $\partial_+$ and which admits a metric with positive scalar curvature. This shows, in particular that  in the following three cases, 
   
   {\it $V$  can't be diffeomorphic  to $X\times [-1,1]$,  where $X$ admits no metric with $Sc>0$},

   (i) $X$ is a {\it spin manifold} , e.g. as  in the above {\bf A }and {\bf B}.

  (ii) $X$ is  $SYS$ as in [SY 1979]   
  or a manifold  as in  [GH 2024].   
  
  (iii) $X$ is an aspherical manifold of dimension| $\leq 5$ or a closely related manifold (see  [Cho-Li 2020], [Gro 2021])

  (These  (i), (ii) and (iii)  cover all {\it known} classes of manifolds, except  for dimension 4, which admit no metrics with $Sc>0$.)
 
 This second  proof  also encounter the  singularity problem  for $dim(V)>7$,  
 where it is more serious than in the first proof, since the Schoen-Yau  partial 
 regularity theorem is not sufficient in this case.
 
 However if $dim(V)=8$ then a required desingularisation follows by a version of Nathan Smale argument  and if $n=9, 10$, then the desingularisation   from
  [Cho-Ma-Sch 2023] {\color {magenta}most probably} apply in the present case.

{\bf 3.} The third proof, which relies on  the  generalized Callias-Dirac  operators technique (see Cecc-Zeid 2023], [Guo-Xie-Yu 2022]),  needs $V$ to be a spin manifold.

This proof applies, in particular,   to  $V$   diffeomorphic to $X\times [-1,1]$, 
where $X$ admits no metric with $Sc>0$, ad  where  non-existence of such a metric    follows via the index theorem for a generalized Dirac operator,  as for instance, for $X$ from the above     
 {\bf A }and {\bf B}.\space{1mm}
 
 As far as the curvature of immersion is concerned, this is  most useful for the Hitchin's spheres $\Sigma^{n}$ for $n=8l+1, 8l+2$ and  which admit immersions to $\mathbb R^{n+1}$ by Hirsch theorem\footnote  {Lichnerowicz's manifolds, which have non-zero $\hat A$-genus admit no Euclidean immersions with codimenension one and two.}
 and all immersions $\Sigma^{n}$ to the unit $(n+1)$ ball satisfy the same inequality as  tori
 {\color {blue}$$curv^\perp \big(\Sigma^{n}\hookrightarrow  B^{n+1}(1)\big )\geq \frac {n+1}{\pi}-1$$}
 and, by a similar argument, 
 { \color {blue}$$curv^\perp \big(\Sigma^{n}\hookrightarrow  B^{n+2}(1)\big )\geq \frac {n+2}{\pi}-2.$$}
  These inequalities can be improved for small $n$ the same way as in the above (b)  for tori, but unlike conjecture {\bf G}  for tori, 
 there is {\it \color {magenta} no (known)  reason to expect} that immersions of $\Sigma^n$ to the unit  balls
 $B^{n+k}$, 
  $k\geq 3$, satisfy $curv^\perp \geq const_k n$. 
  )

  {\it \bf \color {magenta} Question.} Do {\it all}  Milnor's spheres $\Sigma^n$, including those, which carry metrics with $Sc>0$,  develop   large  normal  curvatures when   immersed
   to the balls $B^{n+1}(1)$?



    \subsection {
  $\mathbb T^\rtimes$-Stabilized  Scalar Curvature.}


 Given a compact Riemannian manifold $X$, let 
$$Sc^\rtimes (X)= 4\lambda_1^\rtimes(X),$$
where  $\lambda^\rtimes_1(X)$ is the lowest eigenvalue of the operator 
$-\Delta +\frac {1}{4}Sc$ on $X$  with the Dirichlet (vanishing on the boundary) condition.\footnote{See [Gr 2024] for  justification of this definition/notation and for the proofs of the properties of this $Sc^\rtimes$-curvature used in this paper.}

It is easy to see that $Sc^\rtimes$
 is additive for Riemannian products
$$Sc^\rtimes (X_1\times \underline X)= Sc^\rtimes (X)+Sc^\rtimes (\underline X).$$
and,  more relevantly,

 $Sc^\rtimes (X)$ is  {\it decreasing} under  equidimensional locally isometric immersions:\vspace {1mm}

 \hspace {22mm}{\sl if $X$ immerses to $Y$ then $Sc^\rtimes (X)\geq Sc^\rtimes (X)$.}
\vspace {1mm}

{\bf About $-\Delta +\beta\cdot Sc$.}
The two above  relations remain valid for the first eigenvalues  of the operators 
$$f(x)\mapsto -\Delta f(x)+\beta\cdot Sc(X,x)\cdot f(x)$$ 
for all $\beta\geq 0$,
but $\beta =1/4$ is essential  for the { $2\pi\over \sqrt{Sc^\rtimes}$-inequality} below.\vspace {05mm}

Besides 
$1/4$, a significant value is $\beta = {1\over 4} \frac {n-2}{n-1}$, where positivity of the operator  $-\Delta_X +\beta\cdot  {1\over 4} \frac {n-2}{n-1}Sc(X)$ for $n\geq 3$ on $X$ implies 
that $X$ admits a  metric with positive scalar curvature (as in the proof of the Petrunin's inequality in section 14.2).

Since   ${1\over 4} \frac {n-2}{n-1}< {1\over 4}$ the inequality  $Sc^\rtimes >0$ also implies   
 the existence of a  metric with positive scalar curvature on $X$.
  
  This shows that\ the conditions {\it {\color {red!50!black}$\mathbf\nexists$}PSC and {\color {red!50!black}$\mathbf\nexists$}PSC$^\rtimes$ are equivalent}.
  
  But unlike how it is with the effects  of the positive  signs of $Sc(X)$ and of $Sc^\rtimes(X)$   on the {\it topology} of 
  $X$,  the   $Sc(X)$ and $S^\rtimes (X)$ plays    different roles in  the geometry of $X$.

  \vspace {1mm}

  Let $V$ be a Riemannian manifold  homeomorphic to the product 
$X\times [-1,+1]$, where $X$ is a {\it basic {\color {red!50!black}$\mathbf\nexists$}PSC
$n$-manifold}, i.e. where the underlying reason  for non-existence of a metric with $Sc>0$ on $X$ is {\it  of the  same kind  as what is    presented in section 16.1}.\footnote {Conjecturally, all {\color {red!50!black}$\mathbf\nexists$}PSC manifolds will do, at least for $n\neq 4$} For instance $X$ is diffeomorphic to  the product of the  {\it torus by  Hitchin's sphere}.

{\bf $2\pi\over \sqrt{Sc^\rtimes}$-Inequality}.\footnote {See [Gro 2024] for more about it.}
  Let $V$ be a Riemannian manifold  homeomorphic to the product 
$X\times [-1,+1]$, where $X$ is a {\it basic {\color {red!50!black}$\mathbf\nexists$}PSC
$n$-manifold}, i.e. where the underlying reason  for non-existence of a metric with $Sc>0$ on $X$ is {\it  of the  same kind  as what is    presented in section 16.1}.\footnote {Conjecturally, all {\color {red!50!black}$\mathbf\nexists$}PSC manifolds will do, at least for $n\neq 4$} For instance $X$ is diffeomorphic to  the product of the  {\it torus by  Hitchin's sphere}.

Then the  distance between the two boundary components of $V$ is {\it bounded as follows:
{\color {blue}$$dist( X\times \{-1\}, X\times \{+1\})\leq 2\pi\sqrt\frac { n}{Sc^\rtimes (V) (n+1)}.$$}

{\it Examples of Evaluation of $Sc^\rtimes$.}} 
The rectangular solids satisfy
{\color {blue}$$Sc^{\rtimes}\left( \bigtimes_1^n[-a_i,b_i]\right)=4\sum_1^n \lambda_1[a_i,b_i]=
\sum_1^n \frac {4\pi^2}{(b_i-a_i)^2},$$}
the unit hemispheres  satisfy:
{\color {blue}$$Sc^\rtimes \left(S^n_+\right )=n(n-1)+4n=n(n+3),$$ }
 the unit balls  satisfy 
 $${\color {blue}Sc^\rtimes (B^n)=4j_\nu^2,}$$
 for the first zero of the Bessel function $J_\nu$,   $\nu=\frac {n}{2}-1$, where 
 $j_{-1/2}=\frac {\pi}{2}$, $j_0=2.4042...,$  $j_{1/2}=\pi$ and if $\nu>1/2$, then

   { \color{blue} {$$\nu+\frac{a \nu^\frac{1}{3}}{2^\frac{1}{3}}<j_\nu<
\nu+\frac{a \nu^\frac{1}{3}}{2^\frac{1}{3}}+\frac{3}{20}\frac{2^\frac{2}{3}a^2}{\nu^\frac{1}{2}}$$}}
   where $a=\left(\frac {9\pi}{8} \right)^\frac{2}{3}(1+\varepsilon)\approx 2.32$  with
  $\varepsilon< 0.13\left(\frac {8}{2.847\pi}\right )^{2}<0.1.$    
\vspace {1mm}

{\bf Corollary.} Let $X$ be  a basic {\color {red!50!black}$\mathbf\nexists$}PSC$^\rtimes$  manifold  of dimension $n-1$, e.g. $X=\mathbb T^{n-1}$, and $f: X\to B^{n}(r)$ be a smooth immersion. 
 
 Then the  focal radii of immersions $X\hookrightarrow  B^{n}(r)$  satisfy:
{\color {blue} $$ rad^\perp\big( X\hookrightarrow B^{n} (r)\big)\leq  {\pi r\over 2j_\nu} \sqrt\frac { n}{ n+1}\leqno {[foc.rad]_{j_\nu}}$$}
and 
{\color {blue} $$ curv^\perp\big( X\hookrightarrow B^{n} (r)\big)\geq  \left({2j_\nu\over \pi r}\sqrt {n+1\over n}\right)-r
 \leqno {[curv^\perp]_{j_\nu}} $$}
where $${2j_\nu\over \pi r}\geq 
{n-1/2+3.68(n/2-1)^{1/3}\over \pi r}$$







{\color {blue!22!black}This implies, in particular,  the {\it low curvature bounds from the {\color {blue}$\mathbb T^\rtimes$-remark} in section 15}. 

Also this can be used along with the following.

{\bf Mean Curvature/Ricci $4j_\nu^2$--Inequality.} Let $Y$ be a compact connected  Riemannian $n$-manifold with a non-empty boundary, such that {\it the Ricci curvature of $Y$ is nonnegative}, e.g.  $Y$ is a bounded  Euclidean domain, and {\it the mean curvature of the boundary of $W$ is bounded from below} by that of the unit ball,
$$mean.curv(\partial Y)\geq n-1= mean.curv(\partial B^n).$$
Then 
{\color {blue}$$Sc^\rtimes(Y)\geq  Sc^\rtimes(B^n)=4j_\nu^2.$$}

Thus, the above inequalities

{\color {blue}$ {[foc.rad]_{j_\nu}}$ and $ {[curv\perp]_{j_\nu}}$
remain valid for immersions $X\hookrightarrow Y_r$ for {\it all compact connected  Riemannian $n$-manifolds $Y_r$ with non-empty  boundaries, such that $Ricci(Y_r)\geq 0$ and  $mean.curv(\partial Y_t)\geq {n-1 \over r}$}}.

 {\it Remark/Question.} Let $V\subset \mathbb R^n$ be a bounded domain with two boundary components,    
let $d(V)$ be the distance between these componetns and let $\lambda_1(V)$  the first eigenvalue  
of the Dirchlet problem in $V$. 

 The above shows that \vspace {1mm}
 
  \hspace{-1mm}{\it topology of $V$ may impose a non-trivial bound  on  the product $d^2(V)\lambda_1(V)$}.\vspace {1mm}

{\color {magenta}What are} other  cases of a similar role of the topology of a $V\subset \mathbb R^n$ on metric invariants of $V$?

\subsection  {Curvatures of  Regular Homotopies  of Immersions}


 Due to the Atiyah-Singer index theorem for families of Dirac operators, the index  theoretic obstructions to  $Sc>0$ apply to families of metrics with $Sc>0$, which imply the following
 (see [Hit 1974])

{\bf 16.3.A. The spheres $S^{n-1}$,  $n = 8k+1,8k+2$, $k=1,2,...$.} These  admit (Smale/Milnor) 
$$\mbox{  diffeomorphisms } \mu:S^{n-1}\to S^{n-1},$$
such that the  usual spherical  metric $g_o$  ($sect.curv(g_o)=1$) and   the induced 
metric $g_o^\ast=\mu^\ast(g_o)$  (also $sect.curv(g^\ast_o)=1$)  {\it    can't be joined by a $C^2$-continuous  homotopy $g_t$, such that  $Sc(g_t)>0$}. 

(The  diffeomorphism $\mu$ establishes an {\it isometry} of $(S^{n-1} ,g_o^\ast)$ with the usual sphere $(S^{n-1},g_o)$, where Milnor's  theorem doesn't allow 
a homotopy $g_t$ between $g_o$ and $g_o^\ast$, such that  the metrics $g_t$ have {\it constant sectional} curvatures.)

{\bf 16.3.B.    $O(\sqrt n)$-Curvature Corollary.}} 
Let $f_o: S^{n-1}\to S^{n}(1)$,  be the standard equatorial embedding of the sphere and let $f_t: S^{n-1}\to S^{n}(1)$, $t\in [0,1]$, be a $C^2$-continuous regular homotopy, (a family of $C^2$-immersions\footnote {Such a family does exist by the Smale immersion theorem.}) between   $f_o$ and 
$f^\ast _o=f_o\circ \mu : S^{n-1}\to S^{n}(1)$.
Then {\it there exists} $t_0\in [0,1]$, such that  the normal  curvature of the immersion $f_{t_0}$  satisfies  : 
$$curv^\perp \big(  S^{n-1}\overset{ f_{t_0}} \hookrightarrow S^{n}(1)\big)\geq \sqrt {n-2}. $$

Indeed, if  $curv^\perp \big(  S^{n-1}\overset{ f_{t_0}} \hookrightarrow S^{n}(1)\big)< \sqrt {n-2}. $ for all $t$ then, by Petrunin's Gauss formula from section14.1, the $f_t$-induced metrics $g_t$  on $S^{n-1}$ would  have $Sc>0$ in contradiction with 16.3.A.

{\bf 16.3.C. $O( n)$-Curvature  {\color {magenta}Conjectural} Corollary.} Let $f_o: S^{n-1}\to B^{n}(1)\subset \mathbb R^{n}$   be the standard embedding of the sphere and let $f_t: S^{n-1}\to B^{n}(1)$, $t\in [0,1]$, be a $C^2$-continuous regular homotopy, (a family of $C^2$-immersions\footnote {Such a family does exist by the Smale immersion theorem.}) between   $f_o$ and 
$f^\ast _o=f_o\circ \mu : S^{n-1}\to B^{n}(1)$.
Then {\it there exists} $t_0\in [0,1]$, such that  the normal  curvature of the immersion $f_{t_0}$  satisfies  : 
$$curv^\perp \big(  S^{n-1}\overset{ f_{t_0}} \hookrightarrow  B^{n}(1)\big)\geq j_\nu/\pi>{n+1\over \pi}-1.  $$

To show this one needs an index theorem for families of Callias operators on Riemannian bands.

Milnor's diffeomorphsms seem very different in this regard from the following; 
   
 {{\bf 16.3.D. }   Let $S^n\subset _- B^{n+1}(1)$ be the embedding  obtained from the standard one    $S^n\subset  B^{n+1}(1)$ by an orientation reversing transformation  from $O(n+1).$  
   
    If $n=2,6 \mod 8$, then the two can be joined by a regular homotopy of immersions  
    $S^n\hookrightarrow B^{+1}(1)$.
    ("Turning a sphere inside out" .)

    {\bf Conjecture} A regular homotopy  between these two embeddings  can be achieved for all $n=2,6 \mod 8$ with immersions $f_t$, where  $curv^\perp(f_t(S^n))\leq C$, where $C$ is  a universal constant (probably, $C\leq 100$).

{\bf 16.3.E.  Higher Homotopy Remark.}  There is a body of  results on higher homotopy groups 
of the space $\mathcal G_{Sc>0}(S^n)$ of metrics $g$  with $Sc(g)>0$ on $S^n$, but it is unclear what to do with   (the homotopy   structure of) the map from the space of  immersions $S^n\to B^{n+k}(1)$ (and/or $S^n\to S^{n+k}(1)$) with sufficiently small curvatures  to $\mathcal G_{Sc>0}(S^n)$.
\vspace {1mm}

Not only  Hitchin's spheres but  all  {\color {red!50!black}$\mathbf\nexists$}PSC manifolds $X$ of dimension $n\geq 5$ contain hypersurfaces  $H\subset X$, which support pairs of Riemannian metrics $g_0$ and $g_1$, such that   $Sc(g_i)>0$, $i=0,1$, and where these metrics  {\it    can't be joined by a $C^2$-continuous  homotopies $g_t$, such that  $Sc(g_t)>0$, $0\leq t\leq 1$}. 

To see that,  let $\psi:X\to\mathbb R $ be a Morse function and let $Z=\psi^{-1}(r_0)\subset X$, for some $r_0\in \mathbb R$     be a level of $\psi$, such that all critical point $x\in X$ of $\psi$ with indices $\leq m$ lie below $Z$, i.e. $\psi(x)(x)<r_0$. 

Then $Z$ serves as the common  boundary of the regions $X_0\subset X$ and $X_1\subset  X$, where 
$$\mbox {$X_0=\{x \in X\}_{\psi(x)\leq r_0}$ and $X_1=\{x \in X\}_{psi(x)\geq r_0}$.}$$

 Since $X_0$ represents a regular neighbourhood of a ($\psi$-cellular) $m$-skeleton of $X$
the manifold $X_0$ carries a natural Riemannian metric $g_0$ with $Sc(g_0)>0$, provided  $n-m \geq 3$ and    since  $X_1$ represents a regular neighbourhood of a  $n-m-1$-skeleton of $X$ there is another "natural"metric $g_1$ on $Z$ with $Sc(g_1)>0$ for $m\leq 2$.

Also on knows that if $g_0$ and $g_1$  lie in the same connected component of  $\mathcal G_{Sc>0}(Z)$, then $ X$ admits a metric with $Sc>0$. 

Similarly, 
 if $g_0$ and $g_1$  lie in the same connected component of  $\mathcal G_{Sc^\rtimes>0}(Z)$, 
 then $X$ admits a metric with $Sc^\rtimes >0$. 

{{\bf 16.3.F.   Higher Homotopy \color {magenta} Problem.}  Is there a development of this construction in the spirit of {\it 3.5.D. Higher Homotopy Remark}, e.g. something about the fundamental group  of the space 
$\mathcal G_{Sc^\rtimes>0}(Z')$ for some hypersurface $Z'\subset Z$? 

{\bf 16.3.G.   Toral Example/Question}. Let $X=\mathbb T^n$ and $ 2\leq m \leq n-3$. 
Then, one can show that $Z$ admits an immersion $f_0: Z\to B^n(1)$ 
with 
$$curv^\perp\big (Z\overset {f_0}\hookrightarrow   S^n(1)\big )\leq c_m. $$

It follows that if $n>>m$, then the induced metric  $g_{f_0}$. on $Z$ has $Sc>0$; moreover, one can find  an $f_0$ such that $g_{f_0}$ is homotopic to $g_0$ in $\mathcal G_{Sc>0}(Z)$.

{\color {magenta}When does} $Z$ also admits a similar  immersion $f_1$ to $S^n$ with a sufficiently  small curvature and  a homotopy between $g_{f_1}$ and  $g_1$?

When do manifolds like  $Z$ admit pairs of regularly homotopic immersion $f_0, f_1: Z\hookrightarrow B^n(1)$ with curvatures $\leq c$, yet not regularly homotopic by immersions with curvatures $\leq C$ for some costants $c $ and $C>>c$?


    \section  {Overtwisted Immersions.}
    
Riemannian

    Let $Y=(Y,g=g_Y)$ be a Riemannian manifolds, such as a bounded Euclidean domain, e.g.  the unit ball $B^N(1)\subset\mathbb R^N$.

     An {\it overtwisted}  immersion from $X$ to $Y$  is  a $C^1$ continuous family of 
    smooth immersions $f_t: X\to Y$, $0\leq t<\infty$, such that  the curvature 
    remains bounded 
    $$curv^\perp (f_t(X)\leq C<\infty,$$ such that the induced Riemannian metrics 
    $g_t = f^\ast_t (g_Y)$
    "tends to infinity"  which (to allow non-compact X) is understood  as 
    $g_{t+\delta }\geq \geq 2g_t$ h\ all $t$  and  $\delta >1$.
    
    Here are two  motivating examples of families  of immersions from  the circle  to the  disc of radius $3$ in the $(x,y)$ plane
    $$f_t: S^1\to B^2(3), \mbox  { with } curv^\perp (f_t(X)\leq 1,$$ 
    
  (1)  {\it Winding on a Circle.} This $f_t$ is a family of  immersions from  $S^1$ to  the annulus $A^2(1,3)$ between the unit  circle  $S^1(1)\subset  B^2(3)$  and $S^1(3)=\partial B^2(3)$, which are compositions of two maps that are: 
  
  $\bullet$ embeddings $\phi_t$ from $S^1$ to the band of width 2, i.e. to
  $ \mathbb R \times[-1,1]$, where the image of $\phi_t$  is  the 1-encircling (boundary of the 1-neighbourhood) of the  straight segment of length $t$  in the central line of the band, 
  $[0,t]\subset   \mathbb R\times 0\subset \mathbb R \times[-1,1]$; 
  

 $\bullet$ the covering map   $\mathbb R\times[-1,1]\to A^2(1,3) =S^1(1)\times [-1,1].$ 
    
    \vspace {1mm} 
    
    (2) {\it $\Circle$\hspace {-0.3mm}$\Circle$-Construction.} Another family is  obtained by  repetitively using (compare with 1.C)   regular homotopy from the $1$-encircling $f_0(S^1)$ of the  interval  $[-2,2]$ in the $x$-line   to an immersion $f_1$ with the image  in the union of three  unit circles  with the centres at  
    $-2,0,2$ on the $x$-line,
    $$f_2: S^1\to \mbox { $\Circle$\hspace {-0.3mm}$\Circle$\hspace {-0.3mm}$\Circle$}\subset B^2(3),$$ 
     where the two parallel horizontal bars (of lengths = 4) in  the (convex curve) $f_0(S^1)\subset B^2(2)$  are moved to the two arcs of the  central unit circle, where the upper bar is pushed down to the lower arc and the lower bar is pushed  up to the  arc on the top of this circle.
    
  (This is achieved by moving the the third  disc $\Circle_3$ (positioned on the   right) 
  to  the  position  of 
  $\Circle_1$ by "rolling" $\Circle_3$ along $\Circle_2.$)

  {\it Remark.}   The resulting immersion  $f_1$ from $S^1$ to $\Circle$\hspace {-0.3mm}$\Circle$\hspace {-0.3mm}$\Circle$ can be regularly homotoped with $curv^\perp\leq 1$  inside $B^2(3)$ to $f_2 :S^1\to \Circle$\hspace {-0.3mm}$\Circle$, which is contained in the smaller disc $B^2(2)\subset B^2(3)$.
    
   {\bf  \color {magenta} Conjecture.} There exist   no  regular homotopy  with curvature $\leq 1$ from $f_1$ 
    to $f_2$ within a disc or radius $r<3$. 
    
    \vspace{1mm}

 A version of the  $\Circle$\hspace {-0.3mm}$\Circle$-construction  can be adapted 
   to families of maps, and also to 1-dimensional foliations. Here is a potential application.
   
   \textbf  {16.A. Parametric $1D$-Approximation  {\color {magenta} Conjecture.}}\footnote {This is announced in [Gro 2023]. The proof  I had in mind is technical,  I don't intend 
writing it and would be only happy   if somebody else does it. }   
{\sf Let $\tau $ be a smooth non-vanishing vector field   on a manifold $X$
Then there exists  a smooth map from $X$ to the disc $B^2(4+\varepsilon)$ for a given $\varepsilon>0$, such that the orbits of $\tau$ are sent to  smooth (immersed) curves with curvatures 
$curv^\perp\leq 1$
 in the disc.}

\footnote{Probably, the minimal  possible radius of the receiving disc is $3+\varepsilon$, 
 where "extremal members"  in 
  families of immersions $S^1\to B^2(3-\varepsilon)$  with curvatures $\leq 1$ must be    associated with
 certain patterns comprised    of  unit circles in $B^2(3-\varepsilon)$  tangent to  the unit circle  cantered at zero.
 and where such a pattern  associated with  a regular isotopy from the unit circle to an immersion with the image {\sf 8}
 is comprised of three mutually tangent unit circles inside $B^2(1_2/\sqrt 3$ (which  seems to imply that $\Xi_1> 1+2/\sqrt 3\approx 2.1547...$).}


\textbf  {Example.} There exists a smooth map $f: S^{2n+1} \to  
B^2(4+\varepsilon )$, for all $n=1,2,...$ and $\varepsilon>0$,   such that  the $f$-images of the Hopf circles 
are smooth immersed circles  with  curvatures 
 $curv^\perp \leq 1.$

\vspace{1mm}

The above (1) and (2) generalize to  immersions of $n$-manifolds  for $n>1$, yet  only in a limited way.

{\bf 16.B.} {\color {magenta} $n$-Dimensional Overtwisting  Conjecture.}   Smooth   immersions  of compact $n$-manifolds to the unit balls, $f_0:X\to B^N(1)$, must admit {\it overtwisted} regular homotopies     
$f_t:X\to B^N(1)$, i.e. where the induced metrics $g_t$ in $X$ tend to infinity,  and such that the curvatures of $f_t$ are bounded as follows:
$$curv^\perp (f_t(X))\leq C_n curv^\perp (f_0(X)) \mbox  { for  } t\leq  1
\mbox {  and 
$curv^\perp f_t(X)\leq C_n$ for   $t\geq 1 $}, $$
where, ideally,    $C_n\leq 100n$.

Overtwisting immersion  $f_t$ (with no  control of the curvature for the initial values of $t$)  of special $n$-manifolds $X$
can be obtained with 
"winding  $X$ on an $n$-torus". 

The existence of these. {\sl overtwisted $f_t:X\to B^N(1)$,
such that
 $$curv^\perp( f_t(X))\leq 1000 n^{3/2}\mbox { for $t>1$},
 $$
 is not hard to show in two cases. 
\footnote {See  section 4.4. in [Gro 2022]. There is 
  also   an attempt  in this paper  
 to prove a version of the above n-dimensional overtwisting  conjecture,  but  there is  an error  at  Step 2 in 4.3.A of the "proof".}

(i) $X$ is an orientable  $n$-manifold, which admits an immersion to  $\mathbb R^{n+1}$

(ii) $X$ is an $n$-manifold, which admits an immersion to $\mathbb R^{N+1}$.}

{\bf 16.C.} {\color {magenta}   Curvature  Stable  Flexibility  Conjecture.}  Let $Y$ be a compact Riemannian manifold, e.g the unit ball $B^N(1)$ or the unit $n$-sphere $S^N$

Loosely speaking, the conjecture claims the existence of  a constant $C(Y)<\infty $, such that the all homotopy theoretic invariants 
of immersions $f$ and subspaces  $Im_C(X,Y)\subset Im_C(X,Y)=Im_\infty C(X,Y)$  of immersions $f:X\hookrightarrow Y$ with curvatures $curv^\perp (f(X))\leq C$  do not depend on $C$ for all smooth manifolds $X$ and  $C\geq C(Y)$.

One should  be aware of possible existence of "rigid immersion" 
$f_{\bullet}: X\to Y$ with large  (large)  $curv^\perp(f)=C$, where no small  deformation of  such an  $f_{\bullet}$   decreases the curvature.  (I suspect these exist,  except for $n=1$).
 
To be  safe,  we  {\color {magenta}  conjecture}  that  every compact subset $K\subset Im_(C_1)(X,Y)$ can be brought to $Im_{C_2}(X,Y)$ by a homotopy of $K$ in $Im_{C_3}(X,Y)$  for all 
    $C_1\geq C_2\geq C(Y)$ and  $C_3 \leq C_1(1+C(Y))$.

    {Critical Curvature  \color {magenta} Problem.}  What is the set 
    $\mathcal C^n _{crit}=\mathcal C^n_{crit} (Y)\subset \mathbb R$ of  "critical  values" $C_{crit}$, e.g. where the "homotopy content"  of the subspace $Im_{C}(X,Y)$ increases for some $n$-manifold $X$
   at the point $C=C_{crit}$.
   
   For instance a number $C_\circ$ is critical, if there exits a manifold $X_\circ$ and an 
   immersion $f_\circ :X_\circ \to Y$,  such that $curv^\perp(f_\circ (X_\circ)=C_\circ$
    and there is {\it no immersion} $f$ regularly homotopic to $f_\circ$, such that  
    $curv(f(X_\circ))<C_\circ$.
    
    What is the topology of the set $\mathcal C^n _{crit}(Y)$ e.g. for $Y=B^N(1)$,   $Y=S^N$ and $n=2$?
   
   Is this set finite? discrete?
   
 \vspace {1mm}

\hspace {13 mm}{\it \sc On Immersions between Manifolds with Boundaries} 

One of the problems in proving a  topological Smale-Hirsch type $h$-principle for overtwisted immersions is  typical non-extendability of  immersions from  $X=X_0$ with controlled curvature  to such immersions from the (wide) band 
    $X\times [0,1]\subset X_0=\{0\}\times X$.
   
   On the other hand, much of what we know and what we don't know  about closed manifolds applies to the curvatures of  immersions between pairs of manifolds, especially  to immersions
   $ f: (X,\partial X)\hookrightarrow (Y,\partial Y) $,
   (with $f(X)$ normal to $\partial Y$?), where   the mean curvature of $\partial X\subset X$  may play a similar role  to that of the scalar curvature of $X$.

  (A promising $X$ is  the complement of a small neighbourhood of  the  2-skeleton  of the $n$-torus.)
   \vspace{1mm}

   {\sc {\color {magenta} Final Question.}} Are immersions   with $curv^\perp\leq const$  rigid  or mainly flexible? Is this 
    "hard" or "soft"mathematcs?   


  .


    \section {References}


[A-C-C-D-D-L-O-P 2008] T.    Zachary Abel David Charlton Sébastien Collette Erik D. Demaine
Martin L. Demaine Stefan Langermank
Joseph O’Rourke Val Pinciu
Godfried Toussaint   {\sl Cauchy’s Arm Lemma on a Growing Sphere} 
    https://arxiv.org › pdf › 0804Arxiv

   \vspace{2mm}

[Cecc-Zeid  2023]  Simone Cecchini and Rudolf Zeidler, {\sl Scalar Curvature and Generalized Callias Operators}, in 
Perspectives in Scalar Curvature,
World Scientific Publishing 2023.
\url{https://doi.org/10.1142/9789811273223_0002}

[C-E-M] K. Cieliebak, Y. Eliashberg, N. Mishachev. {\sl Introduction to
the h-Principle}
AMS, 2024.

[Chern, 1961]  Chern, Shiing-Shen  {\sl Studies in global geometry and analysis}
Publication date 1967
Topics Geometry, Differential, Global analysis  \vspace{2mm}

  [Gho-Li 2020] Otis Chodosh, Chao Li {\sl Generalized soap bubbles and the topology of manifolds with positive scalar curvature}
arXiv:2008.11888 [math.DG] \vspace{2mm}

  [Cho-Ma-Sch 2023] Otis Chodosh, Christos Mantoulidis, Felix Schulze {\sl Generic regularity for minimizing hypersurfaces in dimensions 9 and 10}
	arXiv:2302.02253 [math.DG]\vspace{2mm}

 [Conn 1982]  {\sl Connely Rigidity and Energy}    
Inventiones mathematicae (1982)
Volume: 66, page 11-34 \vspace{2mm}  

[DW 1971] M. P. do Carmo and N. R. Wallach, {\sl Minimal immersions of
spheres into spheres}. Ann. of Math. (2)93, 43–62 (1971).

      [Esch  1986] J.-H. Eschenburg,
    {\sl Local convexity and nonnegative curvature -- Gromov's proof of the sphere theorem.}
Inventiones mathematicae (1986)
Volume: 84, page 507-522.\vspace{2mm}  
   
   [Feld  1965]  A. A. Feldbaum, { Optimal. Control  Systems}. Academic Press, 1965 
 
\vspace{2mm}

 [Ge 2021] J. Ge, {\sl Gehring Link Problem, Focal Radius and Over-torical width}
arXiv
https://arxiv.org › math \vspace{2mm}  

[GH 2024] M. Gromov, B.Hanke {\sl Torsion Obstructions to Positive Scalar Curvature
}\url {https://www.emis.de/journals/SIGMA/2024/069/}\vspace{2mm} 
      
      [GL 1983] M.Gromov, B Lawson, {\sl Positive scalar curvature and
the Dirac operator on complete Riemannian manifolds}, Inst. Hautes Etudes
Sci. Publ. Math.58 (1983), 83-196.\vspace{2mm}

      [Gr-Gr 1987]  D. Gromoll and K. Grove, {\sl A generalization of Berger’s rigidity theorem for positively curved
manifolds,} Ann. Scient. Ec. Norm. Sup., Vol. 20 (1987), 227-239.\vspace{2mm}

   [Gro 1991]   M. Gromov,  {\sl Sign and geometric meaning of curvature}, 
     Conferenza tenuta il 14 giugno 1990.    \vspace{2mm}

    [Gro  2017] M. Gromov, {\sl Metric Inequalities with Scalar Curvature}
arXiv:1710.04655 [math.DG]\vspace{2mm}

     [Gro 2021] M. Gromov,   {\sl Four Lectures on Scalar Curvature}, 
arXiv:1908.10612. \vspace{2mm}

  [Gro 2022] M. Gromov,   {\sl Curvature, Kolmogorov Diameter, Hilbert Rational Designs and Overtwisted Immersions}
https://arxiv.org/pdf/2210.13256v1. \vspace{2mm}

[Gro 2023]  M. Gromov, {\sl Isometric Immersions with Controlled Curvatures}
https://arxiv.org/pdf/2212.06122.\vspace{2mm}
   
  [Gro 2024]  M. Gromov,    { \sl Product   Inequalities for $\mathbb T^\rtimes$-Stabilized Scalar Curvature}
https://arxiv.org/pdf/2306.02932 \vspace{2mm}

  [Hit 1974]    N .Hitchin,  {\sl  Harmonic Spinors. }Advances in Mathematics. September 1974, Pages 1-55. 
    \vspace{2mm}

[Guo-Xie-Yu 2022] Hao Guo, Zhizhang Xie, Guoliang Yu,
{\sl Quantitative K-theory, positive scalar curvature, and band width}.

arXiv:2010.01749v4 [math.KT] 
   \vspace{2mm}

[Hopf 1946] H. Hopf,  {\sl Differential Geometry
in the Large}
Seminar Lectures New York University 1946
and Stanford University 1956
With a Preface by S.S. Chern. \vspace{2mm}

[Kon 1995] H. Konig, {\sl Isometric imbeddings of Euclidean spaces into finite di-
mensional $l_p$ -spaces,} Banach Center Publications (1995) Volume: 34, Issue: 1,
page 79-87.
\vspace{2mm}

  [Mend 2025] R. A. E. Mendes,   Diameter and focal radius of submanifoldsz.}
Springer
https://link.springer.com › content › pdf. \vspace{2mm}

     [Nad 1996] N. Nadirashvili, {\sl  Hadamard’s and Calabi-Yau’s conjectures on negatively curved and minimal
surfaces.} Invent. Math. 126 (1996), 457-465.  \vspace{2mm}

[Ni 2023] Lei Ni,    {\sl A Schur’s theorem via a monotonicity and
the expansion module.}
J. reine angew. Math., Ahead of Print DOI 10.1515/crelle-2023-0061 Journal für die reine und angewandte Mathematik. \vspace{2mm}

   [O'R 2000]   Joseph O’Rourke, \vspace{2mm}   An Extension of Cauchy's Arm Lemma with Application to Curve Development} \url {https://link.springer.com/chapter/10.1007/3-540-47738-1_27} \vspace{2mm}

   [Petr 2023]  Anton Petrunin,   {\sl Gromov's tori are optimal }     https://arxiv.org › abs › 2304Arxiv[2304.00886] \vspace{2mm}
  
 [Petr 2024]  Anton Petrunin, {\sl Veronese minimizes normal curvatures} https://arxiv.org › abs › 2408Arxiv[2408.05909] 
   
    \vspace{2mm}

  [Roz 1961]  E. R. Rozendorn, {\sl The construction of a bounded, complete `
surface of nonpositive curvature,}
Uspekhi Mat. Nauk, 1961, Volume 16, Issue 2, 149–156
    \vspace{2mm}

[Sab 2004] I.K Sabitov, {\sl  Around the proof of the Legendre-Cauchy lemma on convex polygons}. Sabitov, I.Kh. Sibirskij Matematicheskij Zhurnal (2004). Volume: 45. 
 \vspace{2mm}

[Sch-Z 1967] 
Schoenberg, I.J., Zaremba, S.K. {\sl On Cauchy's Lemma Concerning Convex Polygons.} Canad. J. of Math.19, 1062–1071 (1967). Published online by Cambridge University Press:  20 November 2018 \vspace{2mm}

     [Sch 1921]  Axel Schur, 
{\sl Über die Schwarzsche Extremaleigenschaft des Kreises unter den Kurven konstanter Krümmung}, Mathematische Annalen 
Volume 83, pages 143–148, (1921)\vspace{2mm}

 [Sul 2007]  John M. Sullivan  {\sl Curves of Finite Total Curvature} arXiv:math/0606007v2 [math.GT] 24 Oct 2007.\vspace{2mm}

 [SY 1979] R. Schoen and S. T. Yau, {\sl On the structure of manifolds
with positive scalar curvature}, Manuscripta Math. 28 (1979), 159-183.

  \vspace{2mm}

  [SY 2017] R. Schoen and S. T. Yau {\sl Positive Scalar Curvature
and Minimal Hypersurface Singularities}. arXiv:1704.05490 ¯

      \end {document}